\documentclass[11pt]{article}
\usepackage[margin=1in]{geometry}
\usepackage{amsmath} 
\usepackage{amssymb} 
\usepackage{amsbsy} 
\usepackage{amsthm} 
\usepackage{amsfonts}
\usepackage{thmtools}
\usepackage{thm-restate}
\usepackage{enumitem}
\usepackage{bbm}
\usepackage{url}
\usepackage{booktabs}
\usepackage{mathrsfs}
\usepackage{multirow}
\usepackage{chngpage}
\usepackage[percent]{overpic}
\usepackage{subcaption}
\usepackage{float}
\usepackage{subfiles}
\usepackage[math]{cellspace}
\usepackage[font={small}]{caption}
\usepackage{tikz, pgfplots}
\usepackage{graphicx}
\usepackage{chngcntr}
\usepackage[square,sort,comma,numbers]{natbib}
\usepackage{soul}
\usepackage{mathtools}
\usepackage{array}
\usepackage{arydshln}
\usepackage[letterspace=-85]{microtype}
\usepackage{nicefrac}
\usepackage{palatino}
\usepackage{hyperref}

\usetikzlibrary{shapes.geometric,arrows,positioning,calc}
\usetikzlibrary{patterns}
\usetikzlibrary{decorations.pathmorphing,backgrounds,fit,matrix,decorations.pathreplacing}

\newcommand{\E}{\mathbb E}

\newcommand*{\pfstart}{\begin{proof}}
\newcommand*{\pfend}{\end{proof}}

\tikzstyle{vertex} = [circle, minimum size=0.1cm, inner sep=0pt, draw=black, fill=black]
\tikzstyle{circ} = [circle, minimum width=0.5mm, inner sep=0pt,draw,fill]
\tikzstyle{hcirc} = [circle,minimum width=5mm, inner sep=0pt,draw]
\tikzstyle{bcirc} = [circle, minimum width=1.5mm, inner sep=0pt,draw,fill]
\tikzstyle{bhcirc} = [circle, minimum width=1.5mm, inner sep=0pt,draw, dotted ]
\tikzstyle{ept} = [circle,minimum width=0mm, inner sep=0pt, white]
\tikzstyle{txt} = [text width=1.3cm,draw,rounded corners=3pt]
\tikzstyle{ncirc} = [circle,draw=black, inner sep=1pt]

\newcommand{\iprod}[1]{\left\langle #1 \right\rangle}

\theoremstyle{definition}
\newtheorem{definition}{Definition}

\theoremstyle{remark}

\usepackage{bm}
\theoremstyle{remark}
\newtheorem*{claim*}{Claim}
\theoremstyle{remark}
\newtheorem*{remark*}{Remark}
\theoremstyle{remark}
\newtheorem{remark}{Remark}

\theoremstyle{plain}
\newtheorem{proposition}{Proposition}

\theoremstyle{plain}

\theoremstyle{plain}
\newtheorem*{lemma*}{Lemma}

\theoremstyle{definition}
\newtheorem{assumption}{Assumption}

\theoremstyle{definition}
\newtheorem*{assumption*}{Assumption}
\declaretheorem[name=Example, style=definition]{example}
\renewcommand\thmcontinues[1]{Cont.}
\theoremstyle{plain}

\theoremstyle{plain}

\DeclareMathOperator{\mspan}{span}
\DeclareMathOperator{\conv}{conv}
\DeclareMathOperator{\lin}{lin}
\DeclareMathOperator{\rec}{rec}
\DeclareMathOperator{\epi}{epi}
\DeclareMathOperator{\id}{id}
\DeclareMathOperator{\ess}{ess}

\DeclareMathOperator{\cone}{cone}
\DeclareMathOperator{\cl}{cl}

\hypersetup{
    colorlinks=true,
    citecolor=blue,
    linkcolor=blue,
    filecolor=magenta,      
    urlcolor=cyan,
}

\pgfplotsset{compat=1.15}
\tikzset{snake it/.style={decorate, decoration={snake, amplitude=.3mm, segment length=2mm}}}

\allowdisplaybreaks[4]
\begin{document}

\title{A Gauge Set Framework for Flexible Robustness Design}

\author{
    Ningji Wei
    \thanks{
    Department of Industrial, Manufacturing \& Systems Engineering, Texas Tech University, Lubbock, TX 79409 (email: {\tt ningji.wei@ttu.edu}).}
    \and
    Xian Yu
    \thanks{
      Department of Integrated Systems Engineering, The Ohio State University, Columbus, OH 43210 (email: {\tt yu.3610@osu.edu}).}
      \and
      Peter Zhang
    \thanks{
      Heinz College of Information Systems and Public Policy, Carnegie Mellon University, Pittsburgh, PA 15213 (email: {\tt pyzhang@cmu.edu}).}
    }
	\date{}
	\maketitle

\begin{abstract}
  This paper proposes a unified framework for designing robustness in optimization under uncertainty using \emph{gauge sets}, convex sets that generalize distance and capture how distributions may deviate from a nominal reference.
Representing robustness through a \emph{gauge set reweighting formulation} brings many classical robustness paradigms under a single convex-analytic perspective. The corresponding dual problem, the \emph{upper approximator regularization model}, reveals a direct connection between distributional perturbations and objective regularization via polar gauge sets.
This framework decouples the design of the nominal distribution, distance metric, and reformulation method, components often entangled in classical approaches, thus enabling modular and composable robustness modeling.
We further provide a gauge set algebra toolkit that supports intersection, summation, convex combination, and composition, enabling complex ambiguity structures to be assembled from simpler components.
For computational tractability under continuously supported uncertainty, we introduce two general finite-dimensional reformulation methods. 
The \emph{functional parameterization} approach guarantees any prescribed gauge-based robustness through flexible selection of function bases, 
while the \emph{envelope representation} approach yields exact reformulations under empirical nominal distributions and is asymptotically exact for arbitrary nominal choices.
A detailed case study demonstrates how the framework accommodates diverse robustness requirements while admitting multiple tractable reformulations. \\

\end{abstract}

\section{Introduction}
Designing decisions that remain effective under uncertainty is increasingly central to modern optimization and learning applications~\citep{shapiro2021lectures}, a concept we refer to as \textit{robustness design}. Over the past decades, a variety of paradigms have emerged, most notably stochastic programming (SP), robust optimization (RO), and distributionally robust optimization (DRO), each developed to hedge against different types of uncertainty and equipped with its own modeling methods and solution techniques.
Although unifying frameworks have been proposed within individual paradigms \citep{bennouna2022holistic,blanchet2023unifying}, this work aims to offer a cross-paradigm perspective on robustness design that reveals common structural principles and enables modular modeling choices. To illustrate the benefits of our framework, we begin with the following example.

\paragraph{Illustrative Example.}  
Consider a city $\Xi\subseteq \mathbb{R}^2$ partitioned into regions where emergency incidents may occur. We model the random location of the incident as our uncertainty $\xi\in\Xi$.
A planner must determine the location of a response center to minimize the expected travel distance from the response center to the incident. Although the spatial uncertainty $\xi \in \Xi$ is continuous, the regional partition follows established administrative divisions that govern how data are collected and organized.
With only limited observations and seasonally varying incident patterns, the true incident distribution is uncertain.
The planner therefore seeks to construct a nominal distribution that blends empirical observations with prior knowledge, and to design a model that
(i) guards against shifts in incident frequencies, e.g., via $\phi$-divergence, since future incident patterns may differ from historical observations,
(ii) hedges against uneven regional data quality, e.g., via a region-specific Wasserstein metric, since some regions have rich historical records while others have very limited data, and
(iii) ensures robust performance under high-impact, tail events, e.g., via Conditional Value-at-Risk (CVaR), since the response plan should still work well for the most remote or hard-to-reach locations.
\hfill$\triangle$
\ \\

Existing approaches, such as $\phi$-divergence DRO \citep{ben2013robust,hu2013kullback}, Wasserstein-based DRO \citep{mohajerin2018data, blanchet2019quantifying, gao2023distributionally}, and CVaR optimization \citep{rockafellar2000optimization}, capture one of these important aspects of uncertainty or risk in this setting. 
There has been an emerging stream of research that focuses on integrating multiple modeling perspectives but often requires ad-hoc or ambiguity-set-specific constructions, such as merging moment- and distance-based ambiguity sets \citep{cheramin2022computationally}, integrating divergence- and Wasserstein-based formulations \citep{blanchet2023unifying}, or composing Wasserstein DRO with CVaR-type objectives \citep{hanasusanto2015distributionally,yu2022multistage}. It is unclear what the \emph{general principles} are in terms of combining and analyzing different robustness modeling approaches. We answer this question positively in this paper.

To enable this flexible robustness design, derive a clean algebraic interpretation of different modeling perspectives, and provide unified reformulation approaches, we introduce a framework grounded in the classical concept of \emph{gauge sets}, that is, convex, zero-containing sets that serve as generalized ``unit balls'' for measuring deviation. Consider a minimization problem $\min_{x\in\mathcal{X}}f(x,\xi)$, where $x\in\mathcal{X}$ is our decision, $\xi\in\Xi$ is the uncertainty, and $f(x,\xi)$ is the associated cost. Let $f_x := f(x, \cdot)$ denote the random variable induced by fixing a solution $x \in \mathcal X$. Our goal is to hedge against the uncertainty $\xi$.
Given any gauge set $\mathcal V$, we define the following \emph{optimal reweighting problem}
\begin{align}
  \sup_{\nu\geq 0 \in L^2(\mathbb P),\ \mathbb E_{\mathbb P}[\nu]=1} \left\{\mathbb E_{\mathbb P}[\nu\cdot f_x] ~\middle|~ \|\nu - 1\|_{\mathcal V} \leq \epsilon\right\}, \label{eq:gauge-set-reweighting}
\end{align}
where $\nu$ denotes a \emph{distribution reweighting function} from the space of square-integrable random variables $L^2(\mathbb P)$, and $\|\cdot\|_{\mathcal{V}}$ is the gauge function induced by $\mathcal{V}$. 
Every gauge function is convex and positively homogeneous, widely employed in convex analysis~\citep{pryce1973r,drusvyatskiy2020convex} and gauge optimization~\citep{freund1987dual,friedlander2014gauge,aravkin2018foundations} to generalize the notion of distance. 
Accordingly, problem~\eqref{eq:gauge-set-reweighting} admits an intuitive interpretation: it constrains the ``distance'' between the reweighting function~$\nu$ and the nominal distribution (represented by the trivial reweighting~$1$) within a prescribed radius~$\epsilon$. When $\nu$ is set to 1, \eqref{eq:gauge-set-reweighting} becomes the expected cost under the nominal distribution $\mathbb{P}$, which reduces to the risk-neutral SP: $\min_{x\in\mathcal{X}}\mathbb{E}_{\mathbb{P}}[f_x]$. 
Under mild conditions, we show that this primal problem admits the following dual formulation, termed the \emph{upper-approximator regularization problem},
\begin{align}
  \inf_{\alpha \in \mathbb R, w \in L^2(\mathbb P)} \left\{ \alpha + \mathbb E_{\mathbb P}[w] + \epsilon\|w\|_{\mathcal V^\circ} ~\middle|~ \alpha + w \geq f_x\right\}, \label{eq:dual-gauge-set-reweighting}
\end{align}
where $\mathcal{V}^\circ$ denotes the polar set of~$\mathcal{V}$, characterizing the form of regularization imposed. 
This dual perspective provides a transparent interpretation of robustness: non-constant upper-approximators $w$ are penalized by the gauge function induced by the polar $\mathcal{V}^\circ$, thereby extending the intuition of \citet{gao2024wasserstein} from Wasserstein-based formulations to general gauge-induced variations. Beyond providing an interpretable and unifying perspective, this framework further enables flexible and composable robustness design, as illustrated below.

\paragraph{Illustrative Example (Continued).}  
Within this framework, the planner can take any distribution~$\mathbb P$ as the nominal distribution to reflect the knowledge on empirical data and prior belief, and express the three robustness requirements using the $\phi$-divergence gauge~$\mathcal V_{\phi}$, the region-wise Wasserstein gauge~$\mathcal V_{\mathrm{Wass}}$ (Example~\ref{eg:hdro}), and the CVaR gauge~$\mathcal V_{\mathrm{CVaR}}$, respectively. 
The first two gauge sets may be combined either by intersection $\mathcal V_{\text{Comb}}:=\delta_1\mathcal{V}_{\phi} \cap \delta_2\mathcal{V}_{\text{Wass}}$ when the distributional shift must satisfy both conditions simultaneously, or by Minkowski sum
$\mathcal V_{\text{Comb}}:=\delta_1\mathcal{V}_{\phi} + \delta_2\mathcal{V}_{\text{Wass}}$ for maximal robustness. The scalings (gauge set radii) $\delta_i$ control the confidence assigned to each modeling component.
Composing this with CVaR realizes the primal problem \eqref{eq:gauge-set-reweighting} as
$$
\sup_{\substack{\nu_1 \geq 0 \\ \mathbb E_{\mathbb P}[\nu_1]=1 \\ \|\nu_1 - 1\|_{\mathcal V_{\text{Comb}}} \leq \epsilon_1}}
\sup_{\substack{\nu_2 \geq 0 \\ \mathbb E_{\nu_1\mathbb P}[\nu_2]=1 \\ \|\nu_2 - 1\|_{\mathcal V_{\text{CVaR}}} \leq \epsilon_2}}
\mathbb E_{\mathbb P}[\nu_2 \nu_1 f_x]
$$
to obtain the worst-case CVaR tail performance over distributions in $\mathcal V_{\text{Comb}}$.
The associated dual \eqref{eq:dual-gauge-set-reweighting} follows immediately from the gauge algebra (Theorem~\ref{thm:alg}) for polar set computation and Theorem~\ref{thm:gcomp} for gauge composition. For instance, if $\mathcal V_{\text{Comb}}$ is defined as the Minkowski sum, according to Corollary~\ref{coro:multisum} and Theorem~\ref{thm:gcomp}, the dual becomes
\begin{equation}
  \label{eq:caseprob}
\begin{aligned}
  \inf_{\alpha_1,\alpha_2 \in \mathbb R, w_{1}, w_2 \in L^2(\mathbb P)} &~\alpha_1 + \alpha_2 + \mathbb E_{\mathbb P}[w_{1}]+ \epsilon_1\delta_1 \|w_{1}\|_{\mathcal V_{\phi}^\circ} + \epsilon_1\delta_2 \|w_{1}\|_{\mathcal V_{\text{Wass}}^\circ} + \epsilon_2 \|w_2\|_{\mathcal V_{\text{CVaR}}^\circ} \\
  \text{s.t.} &~ \alpha_1 + w_{1} \geq w_2\\
            &~ \alpha_2 + w_2 \geq f_x.
\end{aligned}
\end{equation}
For tractable computation, if one of the polar sets is a Lipschitz gauge (Definition~\ref{def:lipgauge}), the \emph{envelope representation} reformulation (Theorem~\ref{thm:saa}) yields a finite convex program that is exact when~$\mathbb P$ is the empirical distribution and is asymptotically exact for other nominal choices (Corollary~\ref{coro:saa}). 
Otherwise, when only regional moments are relevant, a \emph{piecewise moment parameterization} (Example~\ref{eg:affparam}) projects~$\mathcal V^\circ$ onto the corresponding feature space with preserved robustness (Theorem~\ref{thm:paradual}). Both approaches lead to tractable finite-dimensional reformulations. Section~\ref{sec:case} presents a case study of this problem with multiple reformulations.
\hfill$\triangle$
\ \\

\begin{table}[t]
  \centering
  \small
  \begin{tabular}{cccc}
    \toprule
    Method & Corresponding gauge set $\mathcal V$ & Polar set $\mathcal V^\circ$ & Section \\ \midrule
    CRM   & shifted risk envelope $\mathcal Q$ & functions with bounded $\mathcal Q$-induced penalty& \ref{sec:crm} \\
     CVaR   & shifted non-positive cone  & functions with bounded expectation   & \ref{sec:cvar} \\
     Risk-neutral SP   & bounded set & absorbing set  & \ref{sec:spro}\\
    RO   & absorbing set  & bounded set  & \ref{sec:spro} \\
    MDRO   & moment ball  & polynomials with bounded coefficients  & \ref{sec:mdro}\\
    Type-1 WDRO   & shifted $W_1$ ball  & Lipschitz-1 functions & \ref{sec:wdro1}\\
    Type-$p$ WDRO   & shifted $W_p$ ball  & functions with bounded type-$p$ smoothness & \ref{sec:wdrop}\\
    $\phi$-Divergence & $\phi$-divergence ball  & functions with bounded $\phi^*$-penalty & \ref{sec:divergence}\\ 
    Total variation & Total variation ball  & functions with bounded oscillation & \ref{sec:sum}\\ 
    $\chi^2$-Divergence & 2-norm ball  & 2-norm ball & \ref{sec:case}\\ \bottomrule
  \end{tabular}
  \caption{High-level description of gauge sets to represent various robustness designs. CRM refers to the \emph{coherent risk measure}, $W_p$ ball is the type-$p$ Wasserstein ball, and $\phi^*$ is the convex conjugate of $\phi$. The algebraic rules for integrating these gauge sets in flexible ways are presented in Theorem~\ref{thm:alg}.}
  \label{tab:gauge_sets}
\end{table}

From the example above, the proposed framework admits the following benefits when compared to the existing approaches.
\begin{itemize}
  \item \textbf{Separation of design elements.}
  The dual formulation \eqref{eq:dual-gauge-set-reweighting} decouples the nominal distribution from the distance metric, two elements often intertwined in the reformulations in existing DRO paradigms. 
  In~\eqref{eq:dual-gauge-set-reweighting}, the nominal measure~$\mathbb P$ and the polar gauge~$\mathcal V^\circ$ independently evaluate and regularize~$w$ through the expectation and the gauge function, respectively.
This separation enables a principled design of robustness. In contrast to existing data-driven DRO methods that typically build a nominal distribution based on empirical data,
our approach permits the incorporation of both data and model-based information into a hybrid nominal
that encodes the best available belief about the underlying distribution. On the other hand, the gauge~$\mathcal V$ specifies the form of uncertainty to be guarded against.
  
  \item \textbf{Modular composition of robustness.}
  Expressing diverse robustness criteria as their corresponding gauge sets (Table~\ref{tab:gauge_sets}) enables a modular design through algebraic operations, including intersection, summation, convex combination, and gauge composition (Section~\ref{sec:toolkit}). 
  Moreover, multiple robustness measures can be incorporated either from the \emph{distributional-deviation} perspective~\eqref{eq:gauge-set-reweighting} or the \emph{objective-regularization} perspective~\eqref{eq:dual-gauge-set-reweighting}, with the dual interpretations immediately derived through polar-gauge computation.
  \item \textbf{Tractable reformulations for continuously supported uncertainty.}
  When the uncertainty support is continuous, problem~\eqref{eq:dual-gauge-set-reweighting} becomes infinite-dimensional. 
  We develop two finite-dimensional reformulations (Section~\ref{sec:comp}) for general problems.
The \emph{functional parameterization} approach enforces prescribed gauge-based robustness through flexible selection of function bases, while the \emph{envelope representation} approach yields exact finite reformulations under empirical distribution and is asymptotically exact for arbitrary nominal distributions (Theorem~\ref{thm:saa}).
Together, these two methods generalize existing reformulation techniques and further decouple reformulation choices from gauge set design, thereby providing enhanced flexibility for tractable, application-tailored computation.
\end{itemize}

Collectively, these developments support a flexible and modular design of the distributional center, distance metric (gauge set), and reformulation method, enabling tailored and composable robustness modeling that aligns more naturally with data geometry and decision priorities.

\subsection{Related Work}
Multiple optimization paradigms have been established in the literature to enhance solution robustness based on available distributional information, including SP, RO, and DRO. We review each of these topics next and end with a connection to gauge optimization.

\paragraph{Stochastic Programming (SP).} This paradigm aims to optimize a certain risk measure of a random outcome (e.g., the expectation or CVaR of the random cost) given a fully known distribution of the uncertain parameters. When the expectation is used as the performance metric, we have a ``risk-neutral'' framework and aim to find a solution that performs well on average. However, focusing solely on minimizing expected costs does not inherently prevent rare instances of exceptionally high costs. In many real-world scenarios, a ``risk-averse'' framework is preferable to a risk-neutral one to ensure reliable performance under extreme situations, where different risk measures can be used as the objective function. Risk-averse optimization has thus been extensively studied in widespread applications such as portfolio optimization \citep{chen2008two}, energy management \citep{tavakoli2018cvar}, and inventory problems \citep{ahmed2007coherent}. We refer to \citet{birge2011introduction, shapiro2021lectures} for detailed discussions about model formulations, solution algorithms, and applications in risk-neutral and risk-averse SP. In particular, \emph{coherent risk measures} (CRMs) have been widely used in the literature since they satisfy several natural and desirable properties. It has been understood that each CRM, in its dual representation, corresponds to the expectation of the reweighted objective function with respect to the worst-case reweighting probability density function chosen from a candidate density set, referred to as a \emph{risk envelope} \citep{shapiro2021lectures}. Each CRM can be uniquely identified by its risk envelope. However, the risk envelope remains in an abstract form in general, and it has an explicit definition only for particular risk measures. For instance, if the decision maker focuses on the tail performance, CVaR can be used to quantify the tail risk, which corresponds to a fairly simple box-constrained risk envelope \citep{rockafellar2000optimization,pavlikov2014cvar}. This paper adopts the worst-case reweighting perspective and aims to present a general yet explicit formulation of the risk envelope by imposing constraints on the reweighting function using \emph{gauge sets}. It turns out that a variety of existing methods in the literature can be linked to this framework, and their corresponding primal and dual gauge sets offer an intuitive interpretation from a distance-regularization perspective.

\paragraph{Robust Optimization (RO).} When we do not have any information on the underlying distribution except for the support set and the worst-case performance over this support set is a primary concern, RO has proven to be beneficial, ensuring that solutions remain effective even under the most adversarial conditions \citep{ben2009robust,ben2002robust}. Significant efforts have been dedicated to deriving duality results and tractable reformulations under various uncertainty sets \citep{bertsimas-sim2004,bertsimas2003robust,wei2024adjustability,han2023finite,hanasusanto2015k,el2021optimality}, yielding impactful results across various application domains in transportation, supply chain management, power system, and operation management \citep{ardestani2021linearized,STZ2019,bandi2019robust,shi2020integration,subramanyam2021robust}.

\paragraph{Distributionally Robust Optimization (DRO).} As a middle ground between SP and RO, when only partial distributional information is available, DRO can be employed to hedge against distributional ambiguity by constructing ambiguity sets containing all plausible distributions. We refer interested readers to \cite{rahimian2019distributionally} for an extensive survey on DRO. Traditional forms of ambiguity sets include (i) moment-based ambiguity sets \citep[see, e.g.,][]{wagner2008stochastic, delage2010distributionally, mehrotra2014cutting, zhang2018ambiguous,yu2022multistage}, and (ii) distance-based ambiguity sets, such as norm-based distance \citep[see][]{jiang2018risk}), $\phi$-divergence \citep[see][]{jiang2016data,bayraksan2015data}, and Wasserstein metric \citep[see, e.g.,][]{mohajerin2018data, blanchet2019quantifying, gao2023distributionally}. Moment-based ambiguity sets consider different moments of the underlying probability distributions, ensuring the optimal decision remains robust against a family of distributions whose moments are within a certain range from the empirical ones \citep{delage2010distributionally}. Distance-based ambiguity sets impose restrictions on the distance between the candidate distribution and the reference one. Leveraging concentration theorems, it has been shown that Wasserstein distance-based ambiguity sets can achieve effective out-of-sample performance \citep{mohajerin2018data,gao2023distributionally}. 
Meanwhile, $\phi$-divergence metrics, extensively applied in statistical inferences, quantify the ``ratio'' between probability measures \citep{pardo2018statistical,jager2007goodness,cressie1984multinomial}. Based on this, divergence-based DRO has been developed to tackle distributional ambiguity by providing a divergence budget to an adversarial opponent \citep{ben2013robust}. Due to its tractability, this method has been applied in various areas such as data-driven SP \citep{bayraksan2015data} and network design \citep{wang2023globalized}.

\paragraph{Connection to Gauge Optimization.}
The concepts of gauge sets and gauge functions have been extensively studied in convex analysis \citep{pryce1973r,drusvyatskiy2020convex}. Their associated duality theory has been developed in the gauge optimization literature \citep{aravkin2018foundations,freund1987dual,friedlander2014gauge}, where gauge functions are used to evaluate the objective function or constraint violations. Building on this perspective and extending gauge optimization to functional spaces, this paper employs gauge sets to measure the distributional distance between the reweighting function and the nominal one, establishing a unified framework for solution robustness. \ \\

We organize the rest of the paper as follows. Section \ref{sec:setting1} establishes the main assumptions and the strong duality of the gauge reweighting problem. Section \ref{sec:reweighting} investigates gauge set designs in existing robustness paradigms. In Section \ref{sec:toolkit}, we develop several technical tools for manipulating and designing gauge sets, demonstrating their utility using several examples. Section \ref{sec:comp} discusses two tractable reformulation strategies to solve the reweighting problem. Section~\ref{sec:case} presents a detailed case study of the illustrative example, providing two finite-dimensional reformulations and their computational analysis. Finally, Section \ref{sec:conclusion} concludes the paper with discussions on future directions. To streamline the presentation, we discuss potential applications to other robustness frameworks in Appendix~\ref{sec:app2} and defer all the proofs to Appendix~\ref{sec:app1}.

\paragraph{Notation.}  
Let $(\Xi,\mathcal F,\mathbb P)$ denote the nominal probability space and $(\tilde{\Xi},\tilde{\mathcal F},\tilde{\mathbb P})$ the true one. Let $\mathcal M(\Xi)$ be the space of finite signed Borel measures on $\Xi$, endowed with the weak$^\ast$ topology, and let $\mathcal P(\Xi) \subseteq \mathcal M(\Xi)$ denote the subset of probability measures. For any $\mu \in \mathcal M(\Xi)$, $\langle f,\mu\rangle$ denotes the integration of $f$ with respect to $\mu$. The space $L^2(\mathbb P)$ consists of square-integrable random variables equipped with the inner product $\langle \nu,w\rangle_{\mathbb P}=\mathbb E_{\mathbb P}[\nu w]$, where the subscript $\mathbb P$ is omitted when clear. For any $f:L^2(\mathbb P)\to\mathbb R$, its convex conjugate is $f^*(w)=\sup_\nu \langle w,\nu\rangle - f(\nu)$.  Given a convex subset $\mathcal V\subseteq\mathcal M(\Xi)$ or $L^2(\mathbb P)$ that contains the origin, we write $\mathcal V^\circ$, $\conv(\mathcal V)$, $\cone(\mathcal V)$, $\mathcal V^\perp$, $\rec(\mathcal V)$, $\lin(\mathcal V)$, $\mathrm{int}(\mathcal V)$, and $\overline{\mathcal V}$ (or $\cl \mathcal V$) for its polar, convex hull, conic hull, orthogonal space, recession cone (i.e., $\{w \mid \gamma w \in \mathcal V,~\forall \gamma \geq 0\}$), lineality subspace (i.e., $\{w \mid \gamma w \in \mathcal V, ~\forall \gamma \in \mathbb R\}$), interior, and closure (under the ambient topology). Given some closed subspace $\mathcal U$, $\mathrm{int}_{\mathcal U}(\mathcal V \cap \mathcal U)$ and $\cl_{\mathcal U}(\mathcal V \cap \mathcal U)$ denote the interior and closure of $\mathcal V \cap \mathcal U$ relative to the subspace topology of $\mathcal U$. If $\mathcal V$ is convex and contains zero, its gauge function is defined as $\|\nu\|_{\mathcal V}:=\inf\{t>0:\nu\in t\mathcal V\}$. 
For a family $\{\mathcal V_i\}_{i\in I}$, we define $\bigoplus_{i\in I}\mathcal V_i$ to be the closure of $\{\sum_{i\in I}\nu_i \mid \nu_i\in\mathcal V_i,\ \nu_i=0\ \text{for all but finitely many }i\}$. We use $\mathrm{id}(\cdot)$ to denote the identity function. A function is called closed if its epigraph is closed, i.e., it is lower-semicontinuous. Given two vectors $x, y$, we use $x\otimes y:=xy^\intercal$ to denote the tensor product and $x^{\otimes k}$ for the $k$-tensor product using the same $x$. Given any matrix $A$, $\mathrm{vec}(A)$ denotes the vectorization.

\section{Optimal Reweighting Problem}
\label{sec:setting1}
We focus on the optimization problem $\min_{x \in \mathcal X} f(x,\xi)$ where $\mathcal X$ is the solution space and $\xi$ is a random vector from some underlying probability space $(\tilde\Xi, \tilde{\mathcal F}, \tilde{\mathbb P})$. We use $\mathbb P$ with a support $\Xi$, called the \emph{nominal measure}, to denote some empirical probability measure of the unknown true measure $\tilde{\mathbb P}$, and use $f_x$ to denote the random variable $f_x(\xi) = f(x, \xi)$, termed the cost distribution. 

\begin{assumption}[Space Regularity]
  \label{asmp}
 Throughout the paper, we assume the following 
 \begin{enumerate}
   \item $\Xi \subseteq \mathbb R^n$ is Polish and closed, and $\Xi \supseteq \tilde \Xi$;
   \item $\mathbb P$ is fully supported on $\Xi$.
 \end{enumerate}
\end{assumption}

We do not pose other restrictions on the type of $\Xi$, which can be continuous, discrete, or mixed.
The modeling choice \( \Xi \supseteq \tilde{\Xi} \) is standard: robust formulations typically posit a design support that covers all plausible realizations. 
Although Assumption~\ref{asmp}.2 differs from those used in data-driven DRO, it primarily serves as a technical device for analytical convenience. 
To the best of the authors’ knowledge, a wide range of reformulations in existing robustness paradigms can be recovered under this setup (see Section~\ref{sec:reweighting}), including those based on discrete nominal distributions (Corollary~\ref{coro:saa}).

\subsection{Gauge Set}
This subsection provides the basic definition and properties of gauge sets.
We begin with the following definition.

\begin{definition}[Gauge Set]
  \label{defi:l2gauge}
  A \emph{gauge set} is any convex subset $\mathcal V \subseteq L^2(\mathbb P)$ that contains $0$ as a relative interior in the subspace $\mathcal R_0:=\{w \in L^2(\mathbb P) \mid \iprod{1, w} = 0\}$, i.e., $0 \in \mathrm{int}_{\mathcal R_0}(\mathcal V \cap \mathcal R_0)$.
For any $\nu \in L^2(\mathbb P)$, the gauge function induced by $\mathcal V$ is defined as
$\|\nu\|_{\mathcal V}:=\inf\{t > 0 \mid \nu \in t \mathcal V\}$.
We define the set of \emph{reweighting functions} as
$\mathcal R(\mathbb P):=\{\nu \in L^2(\mathbb P) \mid \nu \geq 0, \iprod{1, \nu}= \mathbb E[\nu] = 1\}$.
We further define the set of induced probability measures with a variable center $w$ as
$$\mathcal P_{\epsilon\mathcal V, w}:=\{\nu \mathbb P \in \mathcal M(\Xi) \mid  \nu \in w + \epsilon\mathcal V, \iprod{1, \nu}=1, \nu \geq 0\}$$
and denote by $\overline{\mathcal P}_{\epsilon\mathcal V, w}$ its weak$^\ast$ closure in
$\mathcal M(\Xi)$. 
When the center is $1$, we write
$\overline{\mathcal P}_{\epsilon\mathcal V}:=\overline{\mathcal P}_{\epsilon\mathcal V,1}$.
\end{definition}
The constraint in the above problem defines a ``$\mathcal V$-shaped $\epsilon$-ball'' around the nominal reweighting $1$ (see Statement 4 in Proposition~\ref{prop:bipolar}). The requirement of containing $0$ as a relative interior in $\mathcal R_0$ has two implications: the gauge function is continuous when restricted to $\mathcal R_0$, and it allows the nominal $1$ to be perturbed in every direction within the subspace $\mathcal R_0$ corresponding to probability reweightings. When the gauge set $\mathcal V$ is symmetric ($\nu \in \mathcal V \Longleftrightarrow -\nu \in \mathcal V$), full-dimensional, and bounded, then $\|\cdot\|_{\mathcal V}$ is equivalent to a norm. Thus, the gauge function introduces a more liberal notion of length, using $\mathcal V$ as the ``unit ball'' for measurement, a concept commonly introduced and applied in convex analysis \citep{pryce1973r,drusvyatskiy2020convex} and the gauge optimization literature \citep{aravkin2018foundations,freund1987dual,friedlander2014gauge}. In particular, without the boundedness, we call the gauge function a \emph{seminorm}, allowing nonzero elements to have zero length; without the full-dimensionality, we term the gauge function a \emph{pseudonorm}, allowing elements to have infinite length.
The following proposition summarizes basic properties of the gauge sets used throughout the paper. 
We include the derivations for completeness, as our notion of gauge sets differs slightly from the classical definition: 
we do not assume closedness of $\mathcal V$, and instead require that $0$ lies in the relative interior of $\mathcal V$ 
with respect to $\mathcal R_0$.

\begin{restatable}{proposition}{bipolar}
  \label{prop:bipolar}
  The following relations hold for any given gauge set $\mathcal V \subseteq L^2(\mathbb P)$:
  \begin{enumerate}
    \item $\|\epsilon \nu\|_{\mathcal V} = 0$ when $\epsilon = 0$, and equals  $\epsilon \|\nu \|_{\mathcal V}$ for every $\epsilon > 0$.
    \item $\|\nu + w\|_{\mathcal V} \leq \|\nu\|_{\mathcal V} + \|w\|_{\mathcal V}$.
    \item Given any closed subspace $\mathcal U$, if $0 \in \mathrm{int}_{\mathcal U}(\mathcal V \cap \mathcal U)$, then $\|\cdot\|_{\mathcal V}$ is Lipschitz continuous on $\mathcal U$, and $\|\nu\|_{\mathcal V} = \|\nu\|_{\cl_{\mathcal U}(\mathcal V \cap \mathcal U)}$ for every $\nu \in \mathcal U$.
    \item $\epsilon \mathcal V \subseteq \{\nu \in L^2(\mathbb P) \mid \|\nu\|_{\mathcal V} \leq \epsilon\} \subseteq \epsilon \overline{\mathcal V}$ for every $\epsilon > 0$.
    \item $\|\cdot\|_{\mathcal V}$ is convex, and is closed if $\mathcal V$ is.
  \item $\overline{\mathcal V} = \mathcal V^{\circ\circ}$.
    \item $\ker(\|\cdot\|_{\mathcal V}) = \rec(\mathcal V)$.
  \item $\|w\|_{{\mathcal V}^\circ} = \sup_{\nu \in \mathcal V} \iprod{w, \nu}$.
    \item For every $w \neq 0$, $\|w\|_{\mathcal V} = 0$ implies $\|w\|_{\mathcal V^\circ} = \infty$.
    \item If $w \in \mathcal V^\perp$, $\|w\|_{\mathcal V^\circ} = 0$, and $\|w\|_{\mathcal V} = \infty$ if $w \neq 0$.
  \end{enumerate}
\end{restatable}

\subsection{Optimal Reweighting Problem}

Given a gauge set $\mathcal V$ and radius $\epsilon \ge 0$, we define the associated \emph{optimal reweighting problem} as
\begin{subequations}
  \label{eq:gdist}
\begin{align}
  z_{\epsilon\mathcal V}:=\sup_{\nu(\cdot) \in \mathcal R(\mathbb P)}~& \iprod{f_x, \nu}\\
  \text{s.t.} ~& \|\nu -1\|_{\mathcal V} \leq \epsilon.\label{eq:gdist02}
\end{align}
\end{subequations}
To establish the duality results in the subsequent section, we need some technical pieces to be established. The following lemma introduces the construction of \emph{extended gauge set} that preserves the optimal value.

\begin{restatable}{lemma}{extball}
  \label{lem:extball}
  For every gauge set $\mathcal V$, the \emph{extended gauge} defined as $\tilde{\mathcal V}:=(\mathcal V \cap \mathcal R_0) + \mathcal R_0^\perp$ satisfies (i) $\tilde{\mathcal V}$ contains $0$ as an interior in $L^2(\mathbb P)$, (ii) $z_{\epsilon \mathcal V} = z_{\epsilon \tilde{\mathcal V}}$ for every $\epsilon \geq 0$.
\end{restatable}

Throughout, either of the following two types of regularity conditions will be imposed on a given optimal reweighting problem to serve as a light-tail assumption in the DRO literature \cite{gao2023distributionally,kuhn2019wasserstein}.

\begin{assumption}[Gauge Regularity]
  \label{asm:reg}
  For a given optimal reweighting problem with nominal $\mathbb P$, gauge $\mathcal V$, radius $\epsilon$, and a cost function $f_x$, we assume (i) $f_x$ is continuous on $\Xi$, and (ii) there exists some $\epsilon' > \epsilon$ and a finite-valued, closed, and coercive function $\Phi \geq 0$ (i.e., $\Phi(\xi) \to \infty$ as $\|\xi\|\to \infty$) such that one of the following two conditions is satisfied:
  \begin{itemize}
    \item Type-I: $|f_x| \leq \alpha$ for some $\alpha \in (0, \infty)$, and $\sup_{\mathbb Q \in \overline{\mathcal P}_{\epsilon'\mathcal V}} \mathbb E_{\mathbb Q}[\Phi] < \infty$; or
    \item Type-II: $|f_x| \leq \alpha + \beta \Phi$ for some $\alpha, \beta \in (0, \infty)$, and $\sup_{\mathbb Q \in \overline{\mathcal P}_{\epsilon'\mathcal V}} \mathbb E_{\mathbb Q}[\Phi^{1+\eta}] < \infty$ for some $\eta > 0$.
  \end{itemize}
\end{assumption}
Either form of gauge regularity is sufficient for our subsequent analysis, reflecting a trade-off between enforcing bounds on $f_x$ and controlling tails via the auxiliary function $\Phi$. In particular, when $\Xi$ is compact, both types are satisfied.
Type-I assumes that $f_x$ is bounded, implying uniform tightness of $\overline{\mathcal P}_{\epsilon'\mathcal V}$ (Lemma~\ref{lem:unitight}). 
Type-II drops the boundedness requirement on $f_x$, allowing heavier tails, and instead imposes a stronger light-tail condition on $\overline{\mathcal P}_{\epsilon' \mathcal V}$.
The next lemma summarizes the consequences of gauge regularity.

\begin{restatable}{lemma}{unitight}
  \label{lem:unitight}
  For the given $\epsilon'$ from Assumption~\ref{asm:reg}, Type-I regularity entails that for every $\epsilon \leq \epsilon'$:
  \begin{enumerate}[label=(\roman*)]
    \item $\sup_{\mathbb Q \in \overline{\mathcal P}_{\epsilon\mathcal V}}\mathbb E_{\mathbb Q}[\Phi]< \infty$,
    \item $\overline{\mathcal P}_{\epsilon\mathcal V}$ is uniformly tight,
    \item $\sup_{\mathbb Q \in \overline{\mathcal P}_{\epsilon\mathcal V}}\mathbb E_{\mathbb Q}[f_x]< \infty$,
    \item Given $\epsilon < \epsilon'$, under the extended gauge $\tilde{\mathcal V}$, there exists some $\delta > 0$ such that $\|w - 1\|\leq \delta$ implies $\overline{\mathcal P}_{\epsilon \tilde{\mathcal V}, w}$ is uniformly tight.
  \end{enumerate}
  Type-II regularity additionally entails that (v) for every $\epsilon \leq \epsilon'$, $\sup_{\mathbb Q \in \overline{\mathcal P}_{\epsilon\mathcal V}}\mathbb E_{\mathbb Q}[\Phi\mathbb I_{\Phi > M}]\to 0$ as $M \to \infty$.
\end{restatable}

Robustness models that are based on distributions absolutely continuous with respect to $\mathbb P$, such as coherent risk measures and $\phi$-divergence DRO, are clearly included in the gauge set framework. 
For ambiguity sets defined directly in $\mathcal M(\Xi)$, such as moment-based and Wasserstein DRO, the following proposition shows that \eqref{eq:gdist} is value-equivalent to the corresponding optimization over its probability measure closure under Assumption~\ref{asmp} and \ref{asm:reg}.

\begin{restatable}{lemma}{valeq}
  \label{lem:valeq}
  The functional $\mathbb Q \mapsto \iprod{f_x, \mathbb Q}$ is weak$^\ast$-continuous on $\overline{\mathcal P}_{\epsilon \mathcal V}$. Consequently, let $z_{\epsilon\mathcal V}$ be the optimal value of \eqref{eq:gdist}, the following identity is satisfied,
  $$z_{\epsilon\mathcal V} = \sup_{\mathbb Q \in \overline{\mathcal P}_{\epsilon \mathcal V}} \iprod{f_x, \mathbb Q}.$$
\end{restatable}

This proposition guarantees that as long as worst-case measures can be
approximated in the weak$^\ast$ sense by distributions from
$\mathcal P_{\epsilon\mathcal V}$, the corresponding optimization problems are
value-equivalent under Assumption~\ref{asmp} and \ref{asm:reg}. In particular, the full support Assumption~\ref{asmp}.2 allows for the weak$^\ast$ approximation of measures that are not absolutely continuous with respect to $\mathbb P$. 

\subsection{Dual of the Optimal Reweighting Problem}
\label{sec:grp}
To derive the dual of \eqref{eq:gdist}, we follow the conjugate duality framework introduced in \cite{rockafellar1974conjugate,bot2009conjugate} for generating and analyzing dual problems. Given a convex primal problem $\inf_{x}f(x)$ with a properly constructed convex perturbation function $F(x,u)$ satisfying $F(x,0) = f(x)$, the dual problem can be produced as $\sup_{y} - F^*(0, -y)$ where $F^*$ is the convex conjugate of $F$. A comprehensive list of regularity conditions for strong duality can be found in the paper \cite{bot2009conjugate}. Most of these conditions are designed to guarantee two aspects simultaneously: (i) the primal and dual problems share the same optimal value; (ii) both problems can attain optimality. Since our main interest is to enforce (i) for solution robustness, the following definition and proposition will be used for duality derivation.

\begin{definition}[Quasi-Strong Duality]
  Given a primal problem $\inf_{x} f(x)$ and its dual $\sup_{y} g(y)$, we say the quasi-strong duality holds if $-\infty < \inf_{x} f(x) = \sup_{y} g(y) < +\infty,$ while both optimal solutions may not exist.
\end{definition}

\begin{proposition}[{\cite[p.~11, Theorem 1.4]{bot2009conjugate}}]
  \label{prop:qduality}
  Given that the perturbation function $F: \mathcal X\times \mathcal U \rightarrow \mathbb R \cup \{\pm \infty\}$ is proper and convex, the quasi-strong duality holds if and only if the infimal value function $\phi(u):=\inf_{x \in \mathcal X} F(x,u)$ is finite at $0$ and lower-semicontinuous at $0$.
\end{proposition}

Using conjugate duality, the next theorem derives the dual problem of \eqref{eq:gdist}. The main technical challenge is that the strong duality may not hold, which means none of the strong duality conditions can be directly applied. Instead, we need to prove the quasi-strong duality using Proposition \ref{prop:qduality}. 

\begin{restatable}{theorem}{dist}
  \label{thm:dist}
  The quasi-strong duality holds for the following dual problem of \eqref{eq:gdist}
\begin{subequations}
  \label{eq:distdual}
\begin{align}
  \inf_{\alpha \in \mathbb R, w(\cdot)\in L^2(\mathbb P)}~& \alpha + \mathbb E_{\mathbb P}[w] + \epsilon\|w\|_{\mathcal V^\circ}\label{eq:distdual01}\\
  \text{s.t.} ~& \alpha + w \geq f_x.\label{eq:distdual02}
\end{align}
\end{subequations}
\end{restatable}

\begin{remark}
  The dual problem always provides an upper bound on the primal value by weak duality, irrespective of any continuity or semicontinuity of $f_x$. Hence, the dual formulation can be used as a tractable upper-bounding relaxation of the optimal reweighting problem; additional regularity of $f_x$ in Assumption~\ref{asm:reg} is only needed to ensure the exact value equivalence.
\end{remark}

This reformulation provides an intuitive dual interpretation. The objective function evaluates the expected value of the upper approximation \(\alpha + w\), alongside a penalty on the magnitude of \(w\) gauged by the polar set \(\mathcal V^{\circ}\).
Thus, this result explicitly links the distance and regularization perspectives, enabling robustness to be designed from one side while yielding a dual interpretation via gauge set computation. 

\section{Gauge Set Design in Existing Frameworks}
\label{sec:reweighting}
 Through the gauge set reweighting perspective, this section explores existing robustness paradigms, including general CRM, CVaR, risk-neutral SP, RO, MDRO, WDRO, and $\phi$-divergence DRO, to gain insights into gauge set design patterns. Some results presented here rely on technical tools for gauge set manipulation, which will be fully developed in Section~\ref{sec:toolkit} and are referenced throughout this section as needed.

\subsection{Gauge Set Design in Coherent Risk Measures}
\label{sec:crm}

A CRM is a function $\rho:L^2(\mathbb P) \rightarrow \mathbb R$ that satisfies several axioms to quantify a certain type of risk on cost distributions \citep{artzner1999coherent,shapiro2021lectures}. In this section, we prove that any CRM can be equivalently recast as a gauge set reweighting problem \eqref{eq:gdist}.

According to \cite{artzner1999coherent} and \cite{rockafellar2007coherent}, every CRM adopts a dual representation
$\rho(f_x)=\sup_{\nu \in \mathcal Q} \iprod{f_x, \nu}$ for some convex-closed subset $\mathcal Q \subseteq \mathcal R(\mathbb P)$, where $\mathcal Q$ is called the risk envelope of $\rho$. Moreover, \cite{artzner1999coherent} has shown a one-to-one correspondence between risk envelopes $\mathcal Q$ and CRMs. 
Let $\mathcal Q := \tilde{\mathcal Q} \cap \mathcal R(\mathbb P)$ for some convex and closed set $\tilde{\mathcal Q} \subseteq L^2(\mathbb P)$. Without loss of generality, we may assume that the nominal reweighting $1 \in \mathcal Q$, or equivalently, redefine the reference measure as $\mathbb P := \nu_0 \mathbb P$ for some $\nu_0 \in \mathcal Q$.
Then, the following theorem proves that every risk envelope $\mathcal Q$ can be equivalently described by some gauge set $\mathcal V$.
\begin{restatable}{proposition}{crm}
  \label{prop:crm}
  Every CRM with a risk envelope $\mathcal Q:=\tilde{\mathcal Q}\cap \mathcal R(\mathbb P)$ is equivalent to \eqref{eq:gdist} under the gauge set $\mathcal V=\tilde{\mathcal Q} - 1$ with a radius $\epsilon = 1$. In particular, when $\tilde{\mathcal Q}$ is represented as $\{\nu \in L^2(\mathbb P) \mid g(\nu)\leq 0\}$ for some convex-closed function $g:L^2(\mathbb P) \rightarrow \mathbb R^m$, the polar gauge set is
  $$\mathcal V^\circ=(\tilde{\mathcal Q} - 1)^\circ = \left\{w \in L^2(\mathbb P) ~\middle|~ \inf_{\gamma \geq 0}\iprod{\gamma, g(\cdot)}^*(w) - \iprod{1, w} \leq 1\right\}$$
  where $\iprod{\gamma, g(\cdot)}^*$ is the convex conjugate of the map $\nu \mapsto \iprod{\gamma, g(\nu)}$.
\end{restatable}

An immediate implication is the following explicit form for a general CRM.

\begin{restatable}{corollary}{cmgauge}
  Given a CRM $\rho$ with the risk envelope $\mathcal Q:=\tilde{\mathcal Q}\cap \mathcal R(\mathbb P)$ such that $\tilde{\mathcal Q}:=\{\nu \mid g(\nu) \leq 0\}$ from some convex-closed $g$ satisfying $g(1) \leq 0$, we have
  \begin{align*}
    \rho(f_x) = \inf_{\gamma \geq 0, \alpha, w(\cdot)} \left\{\alpha + \iprod{\gamma, g(\cdot)}^*(w) ~\middle|~ \alpha + w \geq f_x \right\},
  \end{align*}
where $\iprod{\gamma, g(\cdot)}^*$ is the convex conjugate of the map $\nu \mapsto \iprod{\gamma, g(\nu)}$.
\end{restatable}

\subsection{Gauge Set Design in CVaR}
\label{sec:cvar}
For general CRMs, the primal and dual gauge sets are defined abstractly through the representation function $g$. For specific CRMs such as CVaR, the resulting gauge set is more geometrically intuitive.

In CVaR optimization \citep{rockafellar2000optimization}, the $\beta$-CVaR is the conditional expectation of the upper $(1-\beta)$-tail of the cost distribution.  Constraint \eqref{eq:gdist02} can then be written as $\nu \leq (1-\beta)^{-1}$, implying that the reweighting function can increase the original distribution by a factor of at most $(1-\beta)^{-1}$. In this design, the worst-case distribution will move all the probability mass to the upper $(1-\beta)$-percentile, which recovers the CVaR interpretation. The following proposition investigates this constraint under the gauge set perspective.

\begin{restatable}{proposition}{cvar}
  \label{prop:cvar}
 CVaR constraint $\nu \leq 1/(1-\beta)$ is equivalent to 
 $\|\nu - 1\|_{\mathcal V_\beta} \leq 1 \text{ with } \mathcal V_\beta:=\{\nu \mid \nu \leq \beta(1-\beta)^{-1}\}.$
 The corresponding polar gauge set is
 $\mathcal V^\circ_\beta = \{w \geq 0 \mid \beta(1-\beta)^{-1}\mathbb E[w] \leq 1\}.$
 Then, the gauge function is defined as
 $\|w\|_{\mathcal V^\circ_\beta} = \beta(1-\beta)^{-1}\mathbb E[w]$ if $w \geq 0$ and equals $+\infty$ otherwise. This recovers the standard objective function for CVaR optimization as 
  $\inf_{\alpha} \alpha + (1-\beta)^{-1}\mathbb E[(f_x - \alpha)_+]$.
\end{restatable}

In this case, the primal gauge set $\mathcal V_\beta$ is designed as a shifted non-negative cone. The upper bound is deliberately designed to ensure the cut-off point is exactly at the $(1-\beta)$-percentile.

\subsection{Gauge Set Design in Risk-Neutral SP and RO}
\label{sec:spro}
Risk-neutral SP and RO represent opposite ends of robustness: the former optimizes average performance, while the latter guards against the worst case.
In their corresponding gauge-set formulations, this contrast appears through the $\mathcal V$-ball radius $\epsilon$: setting $\epsilon = 0$ in SP restricts reweighting to the nominal function~$1$, whereas a sufficiently large $\epsilon$ in RO renders constraint~\eqref{eq:gdist02} inactive.
Achieving this requires certain basic properties in the design of~$\mathcal V$.

\begin{definition}[Bounded \& Absorbing Set]
 $\mathcal V\subseteq L^2(\mathbb P)$ is bounded if there exists some $L < +\infty$ such that $\|\nu\| \leq L$ for every $\nu \in \mathcal V$; it is absorbing if the origin is an interior point.
\end{definition}

\begin{restatable}{proposition}{bdab}
  If $\mathcal V$ is bounded, then $\ker\|\cdot\|_{\mathcal V}$ is zero; if $\mathcal V$ is absorbing, then $\cone(\mathcal V) = L^2(\mathbb P)$. Therefore, when $\mathcal V$ is bounded, \eqref{eq:gdist} reduces to SP with $\epsilon = 0$; when $\mathcal V$ is absorbing, \eqref{eq:gdist02} becomes redundant when $\epsilon \rightarrow \infty$ and the problem \eqref{eq:gdist} reduces to RO.
\end{restatable}

Such effects are also carried over to the dual problem through the polar gauge set $\mathcal V^\circ$. The following proposition reveals the dual relationship between bounded and absorbing sets.

\begin{restatable}{proposition}{absorb}
  \label{thm:absorb}
$\mathcal V^\circ$ is absorbing if and only if $\mathcal V$ is bounded.
\end{restatable}

\subsection{Gauge Set Design in DRO with Moment-based Ambiguity Sets}
\label{sec:mdro}
Moment-based ambiguity sets have been introduced in the DRO literature to hedge against ambiguity around different moment functions of the nominal distribution, termed MDRO \citep{delage2010distributionally}. For a given nominal distribution $\mathbb P$, the main idea is to construct certain deviation ranges for different moment functions, e.g., the expectation and covariance matrix of $\mathbb P$. Intuitively, these ranges can also be interpreted as some gauge on the distance between the reweighting function $\nu$ and the nominal weight $1$. The following definition generalizes this idea to arbitrary degrees of moment.

\begin{definition}[Generalized Moment Gauge Sets]
  Let $\Omega:\mathbb R^n \rightarrow \mathbb R^n$ be some injective affine transformation and $T_m: \mathbb R^n \rightarrow (\mathbb R^n)^{\otimes m}$ be the $m$-th order tensor product defined as $T_m(\xi) = \xi^{\otimes m}$ with the $(i_1, i_2, \dots, i_m)$-th entry equal to $\xi_{i_1}\xi_{i_2}\cdots\xi_{i_m}$, then the $m$-th moment gauge set can be defined as $\mathcal V_m := \left\{\nu \mid \left\|\mathbb E_{\nu \mathbb P}[T_m \circ \Omega]\right\|_{\mathcal N} \leq 1\right\}$, where $T_m\circ \Omega$ is a random tensor that can be realized at each scenario $\xi_0$ with $T_m\circ \Omega(\xi_0)$, and $\|\cdot\|_{\mathcal N}$ is some compatible norm in the tensor space $(\mathbb R^n)^{\otimes m}$ with $\mathcal N$ being the corresponding unit norm ball.
\end{definition}

The following proposition shows that the classic MDRO constraints can indeed be expressed in terms of these moment gauge sets.
\begin{restatable}{proposition}{mdro}
  \label{prop:mdro}
  Denoting $\mu = \mathbb E[\xi]$ and $\Sigma=\mathbb E[(\xi - \mu)(\xi-\mu)^\intercal]$ as the expectation and covariance matrix of the nominal distribution $\mathbb P$ and $id(\cdot)$ as the identity function, we have the following equivalence,
  \begin{align*}
    (\mathbb E_{\nu \mathbb P}[\xi] - \mu)^\intercal \Sigma^{-1}(\mathbb E_{\nu\mathbb P}[\xi] - \mu) \leq \gamma_1 &\Longleftrightarrow \|\nu - 1\|_{\mathcal V_1} \leq \sqrt{\gamma_1},\\
    \mathbb E_{\nu \mathbb P}[(\xi - \mu)(\xi - \mu)^\intercal]\preceq \gamma_2 \Sigma &\Longleftrightarrow \|\nu - 1\|_{\mathcal V_2} \leq \gamma_2 - 1,
  \end{align*}
  where the affine operator $\Omega_1$ for $\mathcal V_1$ is defined as $\Omega_1:=\Sigma^{-1/2}=\Lambda^{-1/2}Q$ for the eigenvalue decomposition $\Sigma = Q^\intercal \Lambda Q$ with 2-norm on $\mathbb R^n$ as the compatible norm; $\Omega_2$ for $\mathcal V_2$ is defined as $\Omega_2:=\Sigma^{-1/2}(id - \mu)$ with spectral norm $\|A\| = \sigma_{\max}(A)$ extracting the largest singular value as the compatible norm.
\end{restatable}

Therefore, MDRO also falls into the gauge set reweighting problem where \eqref{eq:gdist02} is realized with multiple moment gauge sets. Using Corollary~\ref{coro:multicon} (in Section~\ref{sec:toolkit}) for gauge set intersection, we can directly obtain the dual formulation. The following theorem shows that these moment gauge sets are quite convenient to analyze. We use $\mathfrak J := [n]^{[m]}$ to denote the set of multi-indices of the tensor space $(\mathbb R^n)^{\otimes m}$.

\begin{restatable}{theorem}{mdroi}
  \label{thm:mdro1}
 For every moment gauge set $\mathcal V_m$, the polar set $\mathcal V_m^\circ$ induces a pseudonorm and can be written as $\mathcal V^\circ_m = \{\iprod{X, T_m \circ \Omega} \mid X \in \mathcal N^\circ\}$, where $\iprod{X, T_m \circ \Omega} \in L^2(\mathbb P)$ is defined as $\iprod{X, T_m \circ \Omega}(\xi) = \sum_{J \in \mathfrak J}X_J [T_m \circ \Omega(\xi)]_J$. The corresponding gauge of $w \in L^2(\mathbb P)$ can be explicitly computed as
  $$ \|w\|_{\mathcal V^{\circ}_m} = \begin{cases}
    \|[w]_{T_m\circ \Omega}\|_{\mathcal N^\circ}, & \text{if } w \in \mspan(T_m \circ \Omega)\\
    +\infty, & \text{otherwise},
  \end{cases}
  $$
  where $[w]_{T_m\circ \Omega}:=\arg\min_A\{\|A\|_{\mathcal N^\circ} \mid \iprod{A, T_m \circ \Omega} = w\}$ is the coefficient tensor with respect to $T_m\circ \Omega$, and $\|\cdot\|_{\mathcal N^\circ}$ is the dual norm of $\|\cdot\|_{\mathcal N}$. Moreover, $\mathcal V_m$ induces a seminorm and can be decomposed as $\mathcal V'_m + (\mathcal V'_m)^{\perp}$ with $\mathcal V'_m=\{\iprod{X, T_m \circ \Omega} \mid X \in \mathfrak C^{-1} \mathcal N\}$ and $(\mathcal V'_m)^\perp$ the largest subspace in $\mathcal V_m$ orthogonal to $\mathcal V'_m$,
where $\mathfrak C$ is the symmetric 2-tensor on $(\mathbb R^n)^{\otimes m}$ defined by $[\mathfrak C]_{JJ'}= \iprod{[T_m \circ \Omega]_J, [T_m \circ \Omega]_{J'}}_{\mathbb P}$ for every index $(J,J') \in \mathfrak J^2$. In particular, $\mathfrak C$ is the identity tensor if entries in $T_m \circ \Omega$ form an orthonormal set.
\end{restatable}

This theorem indicates that the polar gauge set $\mathcal{V}_m^\circ$ is obtained by lifting the polar norm ball $\mathcal N^\circ$ into $L^2(\mathbb P)$ through the polynomials from $T_m \circ\Omega$. In particular, $\mathcal{V}_m^\circ$ associated with the classic first-moment constraint is an $L_2$-ellipsoid within the subspace of linear functions, and the one associated with the second-moment constraint induces a spectral-norm-ellipsoid within the subspace spanned by some second-degree polynomials. We also note that polynomials in $T_m \circ\Omega$ are not necessarily linearly independent, thus we define the coefficient tensors to be the ones with the smallest size under $\|\cdot\|_{\mathcal N^\circ}$.
With these pseudonorms used in \eqref{eq:distdual}, only polynomial functions are allowed for upper approximation, which leads to the following corollary.

\begin{restatable}{corollary}{mdroform}
  With the first $m$-th moment constraints $\|\nu -1\|_{\mathcal V_i}\leq \epsilon_i$ for $i \in [m]$ in \eqref{eq:gdist02}, the dual problem \eqref{eq:distdual} is a degree-$m$ polynomial programming
\begin{align}
  \label{eq:mdro}
  \inf_{w(\cdot) \in \mathcal P_m} \left\{ \mathbb E[w] + \sum_{i \in [m]}\epsilon_i\left\|[w]_{T_i \circ \Omega_i}\right\|_{\mathcal N_i^\circ} ~\middle|~ w \geq f_x\right\},
\end{align}
where $\mathcal P_m$ is the space of polynomials of degree less than or equal to $m$.
\end{restatable}

This result echoes the equivalence between MDRO and polynomial programming discovered by \citet{nie2023distributionally}. 
We note that if certain lower-order moment constraints are omitted prior to the $m$-th moment constraint, then the chosen functional basis does not span the entire space $\mathcal P_m$. Instead, it spans only the subspace generated by the functional elements $\{T_i \circ \Omega_i\}$. Another interesting design of MDRO is that $\mathcal V^\circ_i$'s induce pseudonorms so that only specific types of functions (polynomials in this case) can be used for upper approximation, which has the potential to be generalized for other function bases (see Example~\ref{eg:indi} and \ref{eg:affparam}).

\subsection{Gauge Set Design in DRO with Wasserstein Distance-based Ambiguity Sets}
\label{sec:wdro}
Another popular method is to perturb nominal distributions within a certain range measured by Wasserstein distance for gaining solution robustness \citep{gao2023distributionally,mohajerin2018data}. Essentially, the Wasserstein $p$-distance $W_p(\mu,\nu) := \inf_{\pi \in \Pi(\mu,\nu)} (\mathbb E_\pi[d(\xi, \xi')^p])^{1/p}$ defines a metric on the space of probability measures, where $d(\cdot, \cdot)$ is a non-negative and closed cost function and $\Pi(\mu,\nu)$ contains all the joint distributions with marginals $\mu$ and $\nu$. It is well known that the following duality holds for sufficiently general spaces, which we adapt to the Euclidean spaces.
\begin{proposition}[{\citet[p.~19, Theorem 1.3]{villani2021topics}}]
  \label{eq:mdual}
  Suppose $c: \Xi \times \Xi \rightarrow \mathbb R_+ \cup \{+\infty\}$ is a closed function, then we have $W_p(\mu,\nu) = \sup_{\phi(\xi) + \psi(\xi') \leq d(\xi,\xi')} \{\mathbb E_{\mu}[\phi] + \mathbb E_{\nu}[\psi]\}$, where $\phi$ and $\psi$ are continuous and bounded functions.
\end{proposition}

\subsubsection{Gauge Sets of Wasserstein \texorpdfstring{$1$-Distance}{1-Distance}}
\label{sec:wdro1}
When restricting to the probability measures in $\mathcal R(\mathbb P)$, the Wasserstein $1$-distance directly provides the gauge set interpretation according to the following proposition. We omit the proof as it directly follows the Kantorovich-Rubinstein theorem \citep{edwards2011kantorovich}.
\begin{proposition}
  \label{prop:lip}
 Given $\mathbb P$ and $\nu \mathbb P$, the associated Wasserstein $1$-distance is equal to
 $$W_1(\mathbb P, \nu \mathbb P) = \sup_{w \in \text{Lip}_1} \iprod{w, \nu -1} = \|\nu - 1\|_{\text{Lip}_1^\circ},$$ where $\text{Lip}_1$ is the set of non-expanding functions.
\end{proposition}

Hence, for $W_1$ metric, the distance constraint \eqref{eq:gdist02} becomes $\|\nu - 1\|_{\text{Lip}_1^\circ} \leq \epsilon$, and the dual problem \eqref{eq:distdual} uses the $\text{Lip}_1$ gauge set to penalize the upper approximator $w$. We note that $\text{Lip}_1^\circ$ is not a ball defined on the probability measures anymore; instead, it is the original Wasserstein $\epsilon$-ball centered at $1$ translated to the center $0$. The following theorem provides more detailed information.

\begin{restatable}{proposition}{wdroi}
  \label{prop:w1}
The gauge set $\mathcal V_1 = \text{Lip}_1^\circ$ can be written as $\{\nu \mid (\nu+1) \in \mathcal R(\mathbb P), W_1((\nu+1) \mathbb P, \mathbb P) \leq 1\}$. It induces a pseudonorm with $\mspan(1)$ as its orthogonal space. The polar gauge set $\text{Lip}_1$ induces a seminorm with $\mspan(1)$ as its kernel. In particular, $\|w + \alpha\|_{\text{Lip}_1} = \|w\|_{\text{Lip}_1}$ for every $\alpha \in \mathbb R$. 
\end{restatable}

This analysis on the $W_1$ distance will later enable the derivation of the general $W_p$ distance. One immediate result is the following dual problem with respect to the $W_1$ distance constraint.
\begin{restatable}{corollary}{droiform}
  Given the constraint $\|\nu - 1\|_{\text{Lip}_1^\circ} \leq \epsilon$, the dual problem \eqref{eq:distdual} becomes
\begin{align}
  \label{eq:w1dro}
  \inf_{w(\cdot)} \left\{ \mathbb E[w] + \epsilon\|w\|_{\text{Lip}_1} ~\middle|~ w \geq f_x\right\}.
\end{align}
\end{restatable}
This formulation is an infinite-dimensional problem. Two different finite-dimensional solution methods will be introduced in Section \ref{sec:comp}.

\subsubsection{Gauge Sets of Wasserstein \texorpdfstring{$p$-Distance}{p-Distance}}
\label{sec:wdrop}
Using a similar idea as in the $W_1$ distance, we define the gauge set for the $W_p$ distance as follows.
\begin{definition}
  \label{defi:wball}
Let $\mathcal V_{p, \epsilon}:=\{\nu \in L^2(\mathbb P) \mid (\nu + 1) \in \mathcal R(\mathbb P), W_p((\nu + 1)\mathbb P, \mathbb P) \leq \epsilon\}$, the constraint \eqref{eq:gdist02} under the $W_p$ distance can be realized as $\|\nu - 1\|_{\mathcal V_{p,\epsilon}} \leq 1$.
\end{definition}

Since $W_p$ also defines a metric on the probability simplex, the gauge set $\mathcal V_{p,\epsilon}$ also shares the same properties as $\text{Lip}^\circ_1$ in Proposition~\ref{prop:w1}: it is the shifted Wasserstein $p$-ball centered at the origin and is orthogonal to $1$. The following theorem provides an exact description of $\mathcal V_{p,\epsilon}^\circ$.

\begin{restatable}{proposition}{wpdro}
  \label{prop:wp}
  The polar set $\mathcal V^\circ_{p,\epsilon}$ is the following
  \begin{align*}
  \mathcal V_{p,\epsilon}^\circ &= \left\{w \in L^2(\mathbb P) ~\middle|~ \left\{\inf\limits_{\beta \geq 0} \iprod{1,-w(\cdot) - \inf_{\xi} \left\{\beta (d(\xi, \cdot)^p-\epsilon^p) - w(\xi)\right\}}\right\} \leq 1 \right\}.
  \end{align*}
\end{restatable}
We note that the term inside the inner product is the difference between $-w(\xi')$ and its smoothed version $\inf_{\xi}{\beta(d(\xi, \xi')^p-\epsilon^p) - w(\xi)}$, i.e., the infimum convolution of $-w(\cdot)$ with the smoothing term $\beta(d(\cdot,\xi')^p -\epsilon^p)$. Then, the expectation of this difference measures a certain type of smoothness of $w$. We call this quantity the type-$p$ smoothness of $w$. Hence, $\mathcal V_{p,\epsilon}^\circ$ contains functions with their type-$p$ smoothness bounded by one. The corresponding dual problem \eqref{eq:distdual} can be derived in the next corollary, which recovers the general results obtained by \cite{gao2023distributionally}.

\begin{restatable}{corollary}{wpdroform}
  Given Wasserstein $p$-distance $\|\nu -1 \|_{\mathcal V_{p,\epsilon}} \leq 1$, the dual problem \eqref{eq:distdual} becomes
$\inf_{\beta \geq 0} \epsilon^p \beta - \iprod{ 1,\inf_{\xi} \left\{\beta d(\xi, \cdot)^p - f_x(\xi)\right\}}$.
\end{restatable}

\subsection{Gauge Set Design in DRO with \texorpdfstring{$\phi$}{phi}-Divergence-based Ambiguity Sets}
\label{sec:divergence}
Given some convex-closed function $\phi:[0, \infty) \rightarrow \mathbb R$ with additional properties: (i) $\phi(1) = 0$, (ii) $0 \phi(a/0) = a \lim_{t \rightarrow \infty} \phi(t)/t$ for $a > 0$, and (iii) $0 \phi(0/0) = 0$, the corresponding $\phi$-divergence-based worst reweighting problem is defined by realizing \eqref{eq:gdist02} as $\mathbb E[\phi(\nu)] \leq \epsilon$, where $\phi$ acts on $\nu$ in an entry-wise manner by $\phi(\nu)(\xi) = \phi(\nu(\xi))$ \citep{ben2013robust}. The following theorem provides the gauge sets design with respect to $\phi$-divergence.

\begin{restatable}{proposition}{phigauge}
  \label{prop:phigauge}
  Given $\phi$-divergence-based constraint $\mathbb E[\phi(\nu)] \leq \epsilon$, the associated constraint \eqref{eq:gdist02} can be written as $\|\nu - 1\|_{\mathcal V_{\phi,\epsilon}} \leq 1$ for the primal gauge set $\mathcal V_{\phi,\epsilon} = \{\nu \mid \mathbb E[\phi(\nu + 1)] \leq \epsilon\}$. The associated polar set in \eqref{eq:distdual} is 
  $\mathcal V^\circ_{\phi,\epsilon} = \left\{w ~\middle|~ \inf_{\gamma \geq 0}\iprod{1, \gamma (\phi^*(w/\gamma) + \epsilon) - w} \leq 1\right\}$
  where $\phi^*$ is the convex conjugate of $\phi$ and $0 \phi^*(w/0)$ denotes the convex indicator function $\delta_0(w)$.
\end{restatable}

Thus, for any given $w$, we consider the value 
$\inf_{\gamma \geq 0} \iprod{1, \gamma(\phi^*(w/\gamma)+\epsilon) - w}$
as a specific type of penalty on $w$, which we call the $\phi^*$-penalty of $w$. Then, the following corollary provides the dual formulation \eqref{eq:distdual} with respect to $\phi$-divergence.
\begin{restatable}{corollary}{divdual}
  \label{thm:divdual}
  Given $\mathcal V_{\phi, \epsilon}$ as the gauge set in \eqref{eq:gdist02}, the dual problem \eqref{eq:distdual} becomes the following
\begin{align}
  \label{eq:phidual}
  \inf_{\alpha, \gamma \geq 0, w(\cdot)} \left\{ \alpha + \mathbb E[\gamma \phi^*(w/\gamma)]  + \epsilon\gamma ~\middle|~ \alpha + w \geq f_x\right\},
\end{align}
where $\phi^*$ is the convex conjugate of $\phi$ and $0\phi^*(w/0) = \delta_0(w)$. In particular, when $\phi$ is strictly convex and continuously differentiable, $\phi^*$ can be directly computed as $\phi^*(w) = w \cdot (\phi')^{-1}(w)-\phi\circ(\phi')^{-1}(w)$.
\end{restatable}

This corollary provides an intuitive interpretation for DRO with $\phi$-divergence-based ambiguity sets. In the primal problem, the function $\phi$ is designed to measure the divergence of $\nu$ relative to the nominal reweighting function $1$; in the dual problem, it induces the conjugate penalty $\phi^*$ and uses its perspective function to penalize the upper approximation functional $w$ in an entry-wise fashion. In the next section, we will develop technical tools for manipulating multiple gauge sets, facilitating a more flexible approach to robustness design.

\section{Gauge Set Design Methods}
\label{sec:toolkit}
From Section \ref{sec:reweighting}, we observe that various existing robustness solution schemes can be imposed by carefully designing the associated gauge sets. To enable systematic and flexible robustness design, in this section, we develop three technical tools for manipulating gauge sets: (i) an algebraic framework for combining gauge sets, (ii) a decomposition theorem that enables fine-grained design of penalty schemes using selected functional bases, and (iii) a gauge composition theorem for recursively applying multiple robustness requirements.

\subsection{Operations on Gauge Sets and Gauge Functions}
\label{sec:alg}
We begin with some basic properties of gauge sets and gauge functions, providing a convenient toolset for designing gauge sets, as will be shown in later examples. We present the main results in the following two theorems.

\begin{restatable}[Algebra of Gauge Sets and Functions]{theorem}{alg}
  \label{thm:alg}
  Let $\{\mathcal V_i\}_{i \in I}$ be a (possibly infinite) family of convex-closed sets, each of which contains the origin, and let $I_n \subseteq I$ be an arbitrary finite index subset. We define the generalized simplex as
  $\Delta:=\left\{\lambda \in \bigoplus_{i \in I}\mathbb R_+ ~\middle|~ \iprod{1, \lambda} = 1\right\}.$
Then, we have the following results.
  \begin{enumerate}
    \item $(\epsilon \mathcal V)^\circ = \mathcal V^\circ/\epsilon$ for every $\epsilon > 0$.
    \item $\left(\bigcap_{i \in I} \mathcal V_i\right)^\circ = \cl\conv\left(\bigcup_{i \in I} \mathcal V_i^\circ\right)$.
    \item $(\bigoplus_{i \in I} \mathcal V_i)^\circ = \cl\left(\bigcup_{\lambda \in \Delta}\bigcap_{i \in I}\lambda_i \mathcal V_i^\circ\right)$.
    \item $\epsilon \|\nu\|_{\mathcal V} = \|\epsilon \nu\|_{\mathcal V} = \|\nu\|_{\mathcal V/\epsilon}$ for every $\epsilon > 0$.
    \item $\|\nu\|_{\bigcap_{i \in I} \mathcal V_i} = \sup_{i \in I}\|\nu\|_{\mathcal V_i}$.
    \item $\|\nu\|_{\bigcup_{i \in I} \mathcal V_i} = \inf_{i \in I}\|\nu\|_{\mathcal V_i}$.
    \item $\|\nu\|_{\conv\left(\bigcup_{i \in I} \mathcal V_i\right)} = \inf\limits_{I_n\subseteq I, \nu=\sum_{i \in I_n}\nu_i}\sum_{i \in I_n}\|\nu_i\|_{\mathcal V_i}$.
    \item $\|\nu\|_{\bigoplus_{i \in I} \mathcal V_i} = \inf\limits_{I_n\subseteq I, \nu=\sum_{i \in I_n}\nu_i}\max_{i \in I_n}\|\nu_i\|_{\mathcal V_i}$.
    \item $\|w\|_{\bigcup_{\lambda \in \Delta} \bigcap_{i \in I}\lambda_i \mathcal V_i} = \sum_{i \in I}\|w\|_{\mathcal V_i}$, when $I$ is finite.
  \end{enumerate}
\end{restatable}

This theorem enables the computation of the polar gauge set from any compounded primal gauge set and simplifies the polar gauge function representation. Similarly, the following theorem provides a method to express gauge functions in a more specific form.

\begin{restatable}{theorem}{gauge}
  \label{thm:gauge}
  Given any function $g$ that satisfies (i) Non-negativity: $g(w) \geq 0$ for all $w \in L^2(\mathbb P)$ and (ii) Positive homogeneity: $g(\alpha w) = \alpha g(w)$ for every $\alpha \geq 0$, and any gauge set $\mathcal V:= \{w \mid g(w) \leq \epsilon\}$ with $\epsilon > 0$, we have $\|w\|_{\mathcal V} = g(w)/\epsilon$.
\end{restatable}

\subsubsection{Application of Gauge Algebra I: Intersection}
\label{sec:multi}
Combining multiple distance metrics through intersection could reduce the solution conservativeness. The following corollary demonstrates how this intersection in the primal problem \eqref{eq:gdist} influences the dual penalization scheme.

\begin{restatable}{corollary}{multicon}
  \label{coro:multicon}
  Given constraint \eqref{eq:gdist02} as $\|\nu -1 \|_{\mathcal V_i} \leq \epsilon_i$ for all $i \in [m]$, the dual problem becomes
\begin{align}
  \label{eq:pweighted3}
  \inf_{\alpha, w_i(\cdot)} \left\{ \alpha + \sum_{i\in [m]}\mathbb E_{\mathbb P}[w_i] + \sum_{i \in [m]}\epsilon_i\|w_i\|_{\mathcal V^\circ_i} ~\middle|~ \alpha + \sum_{i \in [m]}w_i \geq f_x\right\}.
\end{align}
Moreover, the quasi-strong duality holds if $\mathcal V_i$'s are convex-closed and contain the origin.
\end{restatable}
According to this corollary, using the intersection of multiple distance constraints in the primal problem equips the dual problem with multiple functional components for upper approximating \(f_x\). Then, the objective function measures the expectation of the approximation and applies a size penalty on each component functional \(w_i\) via the gauge set \(\mathcal{V}_i^\circ\). We can consider that each component $w_i$ encodes a certain feature of $w$. Hence, using intersection, we can penalize multiple aspects of the upper approximator $w$. We use the following example for illustration.

\begin{example}[Combination of Multiple Ambiguity Sets I]\label{exp:combination1}
When the underlying distributional ambiguity arises from multiple sources, we may want to combine multiple distributional distance metrics, such as WDRO with $\phi$-divergence \citep{blanchet2023unifying,jin2024constructing} or multiple Wasserstein ambiguity sets, to achieve solution robustness against various sources. For instance, the following reweighting problem
$$\sup_{\nu(\cdot) \in \mathcal R(\mathbb P)} \left\{\iprod{f_x, \nu} ~\middle|~ \|\nu -1\|_{\text{Lip}_1^\circ} \leq \epsilon_1, \|\nu -1\|_{\mathcal V_{\phi,\epsilon_2}} \leq 1 \right\}$$
imposes that the reweighting function $\nu$ should not be too far away from the nominal reweighting function $1$ under both the Wasserstein $1$-distance and $\phi$-divergence metrics. Then, the primal gauge set is the intersection $\epsilon_1 \text{Lip}_1^\circ \cap \mathcal V_{\phi,\epsilon_2}$ with a radius of one. Using Corollary~\ref{coro:multicon}, we immediately obtain the following dual problem
$$\inf_{\alpha, w_1(\cdot), w_2(\cdot)} \left\{ \alpha + \mathbb E[w_1 + w_2] + \epsilon_1\|w_1\|_{\text{Lip}_1} + \|w_2\|_{\mathcal V^\circ_{\phi,\epsilon_2}} ~\middle|~ \alpha + w_1 + w_2 \geq f_x\right\},$$
where two parts $w_1$ and $w_2$ are under distinct penalties. Moreover, due to the generality of our framework, the above duality result remains valid for a broad range of ambiguity sets that can be described using gauge sets.
\hfill$\triangle$
\end{example}

\subsubsection{Application of Gauge Algebra II: Summation}
\label{sec:sum}
Combining gauge sets through summation provides protection against multiple uncertainty sources. For example, an optimal solution obtained under the sum of Wasserstein and KL-divergence gauge sets is simultaneously certified to be robust against both types of distributional perturbations.
The following corollary reveals the effect of this operation on the gauge set design.
\begin{restatable}{corollary}{multisum}
  \label{coro:multisum}
  Given $\mathcal V = \sum_{i \in [m]}\beta_i\mathcal V_i$ in  \eqref{eq:gdist02} for some scalar $\beta_i \geq 0$, the dual problem becomes
\begin{align}
  \label{eq:sumgauge}
  \inf_{\alpha, w(\cdot)} \left\{ \alpha + \mathbb E_{\mathbb P}[w] + \epsilon\sum_{i \in [m]}\beta_i \|w\|_{\mathcal V^\circ_i} ~\middle|~ \alpha + w \geq f_x\right\}.
\end{align}
Moreover, the quasi-strong duality holds if $\mathcal V_i$'s are convex-closed and contain the origin.
\end{restatable}
According to Corollary \ref{coro:multisum}, when adding multiple primal gauge sets, we are also adding their penalty in the dual problem \eqref{eq:distdual}. Thus, it is possible to design multiple gauge sets $\mathcal V_i$ with distinct weights $\beta_i$ to enable a sophisticated robustness solution scheme. In particular, the convex combination of reweighting problems can be seen as a special case of gauge set summation. 
We use the following examples to illustrate the utility of gauge set summation for different robust design purposes.

\begin{example}[Combination of Multiple Ambiguity Sets II]
As an alternative to Example \ref{exp:combination1}, we can also combine multiple ambiguity sets from the dual perspective:
$$\inf_{\alpha, w(\cdot)} \left\{ \alpha + \mathbb E[w] + \epsilon_1\|w\|_{\text{Lip}_1} + \|w\|_{\mathcal V^\circ_{\phi,\epsilon_2}} ~\middle|~ \alpha + w \geq f_x\right\},$$
which penalizes the upper approximator $w$ based on its Lipschitz constant as well as the $\phi^*$-penalty. Applying Corollary~\ref{coro:multisum}, we get the following primal problem
$$
\sup_{\nu(\cdot) \in \mathcal R(\mathbb P)} \left\{\iprod{f_x, \nu} ~\middle|~ \|\nu -1\|_{\epsilon_1 \text{Lip}_1^\circ + \mathcal V_{\phi,\epsilon_2}} \leq 1 \right\}
$$
where the corresponding primal gauge set is the sum $\epsilon_1 \text{Lip}_1^\circ + \mathcal V_{\phi,\epsilon_2}$. Hence, the distance interpretation is that the reweighting function $\nu$ should be near $1$ under this summed gauge set. This method provides a more robust solution than the gauge set intersection as shown in Example \ref{exp:combination1}, since the summation is a superset of the intersection. Again, this duality result also holds for other gauge sets, such as multiple Wasserstein balls \citep{rychener2024wasserstein}.
\hfill$\triangle$
\end{example}

\begin{example}[Flexible Tail-Behavior Selection \& Total Variation Gauge]
  \label{eg:rng}
Utilizing multiple gauge sets, we can extend the idea of CVaR to design flexible tail-behavior selectors as follows.
\begin{align*}
  \inf_{x\in \mathcal X, \alpha}~& \alpha + \sum_{i \in [m]}\epsilon_i \|(f_x - \alpha)_+\|_{\mathcal V_i^\circ}.
\end{align*}
 For instance, when some $\mathcal V_i^\circ$ is $\text{Lip}_1$, the optimal $f_x$ also concerns the Lipschitz constant at the tail part. 
 In contrast, when the polar gauge set is defined as
 $$\text{Osc}_1 := \left\{w ~\middle|~ \frac{1}{2}\left((\sup_{\xi \in \Xi}w(\xi) - \inf_{\xi \in \Xi}w(\xi)\right) \leq 1\right\},$$
the optimal solution \(f_x\) seeks to minimize the oscillation of the objective, while the parameters \(\epsilon_i\) govern the trade-off between tail expectation and tail variation, enforcing smaller dispersion when risks materialize.
A direct computation shows that \(\mathrm{Osc}_1\) is the polar of the total variation gauge
\[
\mathcal V_{\mathrm{TV}}
:=
\left\{
\nu \in L^2(\mathbb P)
\;\middle|\;
\iprod{1,\nu}=0,\;
\iprod{1,|\nu|} \le 1
\right\}.
\]
Indeed, any \(\nu \in \mathcal V_{\mathrm{TV}}\) admits the decomposition 
\(\nu = \nu_+ - \nu_-\) with 
\(\iprod{1,\nu_+} = \iprod{1,\nu_-} = 0.5\).
Maximizing \(\iprod{w,\nu}\) therefore assigns half the mass to \(\sup w\) and half to \(\inf w\), which confirms that \(\mathrm{Osc}_1 = \mathcal V_{\mathrm{TV}}^\circ\) (up to closure).
\hfill$\triangle$
\end{example}

\subsection{Gauge Set Decomposition}
\label{sec:decomp}
In MDRO and WDRO, we observe the dual relationship between seminorms and pseudonorms, reflected through the associated primal and polar gauge sets. The following decomposition theorem offers a more detailed characterization of this relationship, enabling the use of the function basis enforcement technique to address other types of ambiguities.

\begin{restatable}[Gauge Set Decomposition]{theorem}{decomp}
  \label{thm:decomp}
For a closed gauge set $\mathcal V$, we have the following decomposition
 $$L^2(\mathbb P) =  \lin(\mathcal V) \oplus\mathcal V^\perp \oplus \ess(\mathcal V),$$
 where $\ess(\mathcal V):=(\lin(\mathcal V) \oplus \mathcal V^\perp)^\perp$ is termed the essential subspace induced by $\mathcal V$. Define $\mathcal V^\dagger:=\ess(\mathcal V) \cap \mathcal V$ to be the essential gauge set of $\mathcal V$, then $\mathcal V$ and $\mathcal V^\circ$ can be decomposed as
 $$
 \begin{aligned}
   \mathcal V &= \lin(\mathcal V) + \mathcal V^\dagger\\
   \mathcal V^\circ &= \mathcal V^\perp + \left(\mathcal V^\dagger\right)^\circ_{\ess(\mathcal V)},
 \end{aligned}
 $$
where $(\mathcal W)^\circ_{\ess(\mathcal V)}:=\{w \in \ess(\mathcal V) \mid \iprod{w, v} \leq 1 ~~ \forall v \in \mathcal W\}$ is the polar set relative to the essential subspace $\ess(\mathcal V)$ for any $\mathcal W \subseteq \ess(\mathcal V)$. In particular, we have
$$\lin(\mathcal V^\circ) = \mathcal V^\perp, \quad (\mathcal V^\circ)^\dagger = \left(\mathcal V^\dagger\right)^\circ_{\ess(\mathcal V)}, \quad \ess(\mathcal V) = \ess(\mathcal V^\circ).$$
Moreover, $\mathcal V^\dagger$ is convex-closed and contains $0$.
\end{restatable}

This depicts a more intuitive picture regarding the primal and polar gauge sets: the lineality subspace and the orthogonal subspace associated with $\mathcal V$ will swap in its polar set $\mathcal V^\circ$, and the ``essential'' part of the gauge set $\mathcal V$ will be converted to its relative polar set in the essential subspace $\ess(\mathcal V)$. The following example illustrates one usage of this result.

\begin{example}[Indicator Function Basis for Spatial Uncertainty]
  \label{eg:indi}
In this case, $\Xi = \bigcup_{i \in I} \Xi_i$ represents a region that is partitioned into multiple districts. Based on historical data, different districts may have different types of ambiguity. A simple scheme is to define the following polar gauge sets based on the indicator functions basis 
$\mathcal V^\circ_i:=\{r_i \mathbb I_{\Xi_i} \mid |r_i| \leq 1\}$.
Then, the dual problem becomes
$$
\inf_{w(\cdot)=\sum_{i \in I} r_i\mathbb I_{\Xi_i}} \left\{ \alpha + \mathbb E[w] + \sum_{i \in I}\epsilon_i |r_i| ~\middle|~ \alpha + w \geq f_x \right\}.
$$
Hence, every $w \in \mathcal V^\circ$ is a piecewise function with each piece having a coefficient $r_i$. Each piece also has a distinct penalty $\epsilon_i$. From the primal perspective, constraint \eqref{eq:gdist02} becomes 
$$\left|\iprod{\nu - 1, \mathbb I_{\Xi_i}}\right| = \left|\mathbb E[(\nu -1) \mathbb I_{\Xi_i}]\right| \leq \epsilon_i \Longleftrightarrow \nu \mathbb P(\Xi_i) \in [\mathbb P(\Xi_i)-\epsilon_i, \mathbb P(\Xi_i)+\epsilon_i],\ \forall i\in I.$$
That is, the spatial distributional ambiguity at each region $i$ is modeled by the probability variation $\epsilon_i$ from the nominal probability, providing an intuitive distance interpretation.
\hfill$\triangle$
\end{example}

We can further combine this indicator function basis with other penalty methods, as illustrated in the next example.

\begin{example}[Heterogeneous DRO]
  \label{eg:hdro}
Let $\{\Xi_1, \Xi_2\}$ be a partition of the uncertainty space $\Xi$, and suppose the data associated with $\Xi_1$ are more sufficient than in $\Xi_2$. Then, the user may want to mitigate more distributional uncertainty over $\Xi_2$ than $\Xi_1$. Define $\text{Lip}_1^1 := \{w \cdot \mathbb I_{\Xi_1} \mid w \in \text{Lip}_1\}$ and $\text{Lip}_1^2 := \{w \cdot \mathbb I_{\Xi_2} \mid w \in \text{Lip}_1\}$, we can set up the following dual problem
$$\inf_{w(\cdot)} \left\{ \mathbb E[w] + \epsilon_1\|w \cdot \mathbb I_{\Xi_1}\|_{\text{Lip}^1_{1}} + \epsilon_2\|w \cdot \mathbb I_{\Xi_2}\|_{\text{Lip}^2_{1}} ~\middle|~ w \cdot \mathbb I_{\Xi_1} + w \cdot \mathbb I_{\Xi_2} \geq f_x\right\}.$$
This function combines two polar gauge sets according to Corollary \ref{coro:multicon}. Since $\text{Lip}_1^1$ does not contain any functions that have nonzero values on $\Xi_2$, these functions are prevented from usage since their gauge would be infinity. Hence, $w_1$ is a function that has zero values on $\Xi_2$ and has a Lipschitz penalty $\epsilon_1$ on the $\Xi_1$ part. From the primal perspective, the associated distance constraints are
$\|\nu -1\|_{(\text{Lip}^i_1)^\circ} \leq \epsilon_i \text{ for } i \in \{1,2\}.$
Thus, it first projects $\nu-1$ onto the $\Xi_i$ part, then ensures that its $W_1$ distance is less than $\epsilon_i$, realizing a heterogeneous penalty. Although such a modification does not guarantee the global Lipschitz (the changing rate between points in $\Xi_1$ and $\Xi_2$ is not penalized), we can add an additional term $\epsilon\|w_1 + w_2\|_{\text{Lip}_1}$ to fine-tune the global Lipschitz if needed.
\hfill$\triangle$
\end{example}

\subsection{Gauge Set Composition}
When additional regularization is imposed on the worst-case distribution, the resulting formulation involves a composition of gauge sets. The next theorem formalizes this recursive construction and presents its explicit dual representation.

\begin{restatable}{theorem}{gcomp}
  \label{thm:gcomp}
  Given $m$ gauge sets satisfying Type-II regularity in Assumption~\ref{asm:reg}, the \emph{composed primal problem} is
  \begin{equation}
    \label{eq:gcomp}
  \sup_{\substack{\nu_1 \geq 0 \\ \iprod{1, \nu_1}_{\mathbb P}=1 \\ \|\nu_1 - 1\|_{\mathcal V_1} \leq \epsilon_1}}\sup_{\substack{\nu_2 \geq 0 \\ \iprod{1, \nu_2}_{\nu_1\mathbb P}=1 \\ \|\nu_2 - 1\|_{\mathcal V_2} \leq \epsilon_2}} \cdots \sup_{\substack{\nu_m \geq 0 \\ \iprod{1, \nu_m}_{\nu_1 \nu_2 \cdots \nu_{m-1}\mathbb P}=1 \\ \|\nu_m - 1\|_{\mathcal V_m} \leq \epsilon_m}} \iprod{\prod_{i \in [m]}\nu_i, f_x}_{\mathbb P}.
  \end{equation}
  Define $C_{\Phi}(\Xi):= \{w \in C(\Xi) \mid \sup_{\xi\in \Xi} |w(\xi)|/(1+\Phi(\xi)) < \infty \}$ where $C(\Xi)$ is the set of continuous functions over $\Xi$. The associated dual problem is 
  $$
  \begin{aligned}
    \inf_{\{\alpha_i, w_i(\cdot)\}_{i \in [m]}} &~ \sum_{i \in [m]} \left(\alpha_i + \epsilon_i\|w_i\|_{\mathcal V_i^\circ}\right) + \mathbb E_{\mathbb P}[w_1]\\
    \text{s.t.} &~ \alpha_i + w_i \geq w_{i+1}, \quad \forall i \in [m],
  \end{aligned}
  $$
  where $w_{m+1} = f_x$ and $w_i \in C_\Phi(\Xi)$ for every $i \in [m]$.
\end{restatable}

Given the optimal reweighting~$\nu_i$ at level~$i$, the next stage applies a new reweighting~$\nu_{i+1}$ to the updated distribution~$\big(\prod_{k \le i} \nu_k\big)\mathbb{P}$, yielding the composed reweighting problem above. An illustrative example is provided below.

\begin{example}[Tail Performance under Worst-Case Distribution]
  \label{eg:gcomp}
  In distributionally robust risk optimization, a decision-maker facing uncertain outcomes seeks to hedge against tail risk by optimizing performance with respect to the worst-case distribution within a Wasserstein ambiguity set, thereby ensuring reliability under potential model misspecification. 
  This risk attitude can be represented as a two-level composition of gauge sets: a Wasserstein gauge~$\mathcal V_1:=\text{Lip}_1$ with radius $\epsilon$ capturing distributional perturbations, and a CVaR gauge~$\mathcal V_2:=\mathcal V_{\beta}$ with radius $1$ modeling tail sensitivity, where both gauge sets $\text{Lip}_1$ and $\mathcal V_\beta$ are introduced in Section~\ref{sec:reweighting}. 
  According to Theorem~\ref{thm:gcomp}, the associated dual problem is given by
  $$
  \begin{aligned}
    \inf_{\alpha,\, w_1(\cdot),\, w_2(\cdot)} 
      &~ \alpha + \mathbb E_{\mathbb P}[w_1] 
        + \epsilon \|w_1\|_{\text{Lip}_1} 
        + \|w_2\|_{\mathcal V^\circ_\beta}\\
    \text{s.t.}\quad
      &~ w_1 \ge w_2,\\
      &~ \alpha + w_2 \ge f_x.
  \end{aligned}
  $$
  By substituting the definitions of gauge sets, the formulation can be equivalently expressed as
  $$
  \begin{aligned}
    \inf_{\alpha,\, w(\cdot)} 
      &~ \alpha 
        + \mathbb E_{\mathbb P}\!\left[w 
          + \tfrac{\beta}{1-\beta}(f_x - \alpha)_+\right]
        + \epsilon \|w\|_{\text{Lip}_1}\\
    \text{s.t.}\quad 
      &~ w \ge (f_x - \alpha)_+.
  \end{aligned}
  $$
  Finally, a finite-dimensional program can be obtained using either of the reformulation methods introduced in the next section of computational approaches.
\end{example}

\section{Computational Approaches}
\label{sec:comp}

Both the decision variable $w$ and the constraint set \eqref{eq:distdual02} of the dual reweighting problem \eqref{eq:distdual} are indexed by the elements in $\Xi$. When $\Xi$ contains only a finite number of scenarios, the problem is often tractable with a finite number of variables and constraints. Otherwise, \eqref{eq:distdual} is an infinite-dimensional optimization. This section introduces two finite-dimensional approaches, namely \emph{functional parameterization} and \emph{envelope representation}, to handle this challenge, generalizing solution methods in the DRO literature.

\subsection{Functional Parameterization}
\label{sec:comp_param}
In practice, one often focuses on robustness with respect to a finite number of features (e.g., covariances, moments, or probabilities over selected regions). To formalize this, we introduce the functional parameterization method.

\begin{definition}[Functional Parameterization]
Let $\phi = (\phi_i)_{i \in [\ell]}$ be a collection of basis functions with each $\phi_i \in L^2(\mathbb P)$.
For any coefficient vector $\lambda \in \Lambda$ where $\Lambda \subseteq \mathbb{R}^\ell$ is a convex-closed cone, denote the linear combination by $\iprod{\lambda,\phi}:= \sum_{i=1}^\ell \lambda_i \phi_i$.
The associated parametric functional subspace is denoted by $\Lambda_\phi := \{\iprod{\lambda, \phi} \mid \lambda \in \Lambda\}$. Then, we define the \emph{parametric primal problem} under $\phi$ as
\begin{subequations}
  \label{eq:paraprime}
\begin{align}
  \sup_{\nu(\cdot) \in \mathcal R(\mathbb P)} &~ \iprod{f_x, \nu}\\
  \text{s.t.} &~ \|\nu -1\|_{\conv\left(\mathcal V \cup \Lambda_\phi^\circ\right)} \leq \epsilon,
\end{align}
\end{subequations}
where $\Lambda_\phi^\circ$ is the polar cone of the induced cone $\Lambda_\phi$ in  $L^2(\mathbb P)$.
\end{definition}

We observe that~\eqref{eq:paraprime} is always at least as robust as the original problem under the gauge $\mathcal V$, as it corresponds to a superset of the original primal gauge.
In practice, explicitly constructing this parametric primal gauge is unnecessary, since its semi-infinite dual admits a simple and tractable representation, as shown in the following theorem.

\begin{restatable}{theorem}{paradual}
  \label{thm:paradual}
  When $\Lambda_\phi$ is closed in $\mathrm{span}(\phi)$, the dual associated with \eqref{eq:paraprime} is
\begin{subequations}
  \label{eq:paradual}
\begin{align}
  z_{\phi,\Lambda}:=\inf_{\lambda \in \Lambda, \alpha} &~ \alpha + \iprod{\lambda, \mathbb E_{\mathbb P}[\phi]} + \epsilon\|\iprod{\lambda, \phi}\|_{\mathcal V^\circ}\\
  \text{s.t.} &~ \alpha + \iprod{\lambda, \phi} \geq f_x.
\end{align}
\end{subequations}
Otherwise, $z_{\phi, \Lambda}$ is an upper bound of the dual of \eqref{eq:paraprime}.
\end{restatable}

\begin{remark}
Since $\Lambda$ is closed in $\mathbb R^\ell$ and $\mathrm{span}(\phi)$ is a finite-dimensional (hence closed) subspace of $L^2(\mathbb P)$, the set $\Lambda_\phi$ can be viewed as the image $T(\Lambda)$ under the linear map $T(\lambda):=\iprod{\lambda,\phi}$. Consequently, $\Lambda_\phi$ is closed in $\mathrm{span}(\phi)$ if and only if $T$ preserves closedness on $\Lambda$. This property holds, for example, if $\Lambda$ is polyhedral, or $T$ is injective, or $\ker T\cap\Lambda=\{0\}$. When $\Lambda_\phi$ is not closed in $\mathrm{span}(\phi)$, the above reformulation still provides a conservative (i.e., potentially larger) evaluation of $f_x$.
\end{remark}

This theorem enables flexible finite-dimensional parameterizations while preserving robustness. For example, existing moment-based WDRO reformulations (e.g., elliptical reformulation in \cite{kuhn2019wasserstein}) require specific nominal distributions and Wasserstein metrics for tractability, whereas our result allows arbitrary choices of functional bases, nominal distributions, and robustness metrics. Moreover, when each $\phi_i$ is piecewise convex and $f_x$ is piecewise concave, the semi-infinite constraints often admit a finite-dimensional dual reformulation.

\begin{example}[Moment-Based Parameterization]
  Both the classical MDRO model \citep{delage2010distributionally} and the WDRO model with an elliptical nominal distribution \citep{kuhn2019wasserstein} employ variants of moment-based parameterizations with $\phi(\xi) = (\xi,\, \xi^{\otimes 2})$
to extract first- and second-moment information. The distinction between these approaches lies in their choices of the nominal distribution~$\mathbb P$ and the gauge set~$\mathcal V$. In MDRO, $\mathbb P$ is interpreted as an arbitrary distribution characterized only by its mean~$\mu$ and covariance~$\Sigma$, and $\mathcal V$ is the moment-based uncertainty set described in Proposition~\ref{prop:mdro}.  
In contrast, WDRO with an elliptical nominal assumes $\mathbb P$ to be elliptical and uses a type-2 Wasserstein ball to construct the primal gauge.
The theorem above, however, provides a more general perspective where nominal distribution and gauge sets can be independently and flexibly chosen.
For example, one may take $\phi(\xi) = \bigl(\cos \xi_i,\, \sin \xi_i\bigr)_{i\in[n]}$
as the functional basis with guaranteed robustness. Moreover, each expected feature value $\mathbb E_{\mathbb P}[\phi_i]$ can be computed analytically when available, or estimated via sampling when closed-form expressions are unavailable.
\end{example}

\begin{example}[Region-Based (Piecewise-Constant) Parameterization]
The indicator-function basis introduced in Example~\ref{eg:indi} induces a partition of the uncertainty space~$\Xi$ into regions~$\{\Xi_i\}_{i \in [\ell]}$ with corresponding indicator functions~$\mathbb I_{\Xi_i}$. 
The resulting region-based parameterization is given by $\phi(\xi) = (\mathbb I_{\Xi_i}(\xi))_{i \in [\ell]}$.
Overall, this parameterization offers a principled approach to discretizing the support of the ambiguity set and remains fully compatible with different choices of~$\mathcal V$ and~$\mathbb P$.
\end{example}

\begin{example}[Piecewise-Affine Parameterization]
  \label{eg:affparam}
The region-based parameterization can be overly coarse, as it captures only the distributional distance of the zeroth moment within each region.  
To achieve finer control, we introduce the \emph{piecewise-affine parameterization} defined as  
$$\phi(\xi) = \left(\mathbb I_{\Xi_i}(\xi), \xi_j\mathbb I_{\Xi_i}(\xi)\right)_{i \in [\ell], j \in [n]}$$
In addition to the constant basis functions used previously, each new basis functional $\xi_j\mathbb I_{\Xi_i}(\xi)$ encodes the first-moment information within region~$\Xi_i$.  
\end{example}

\subsection{Lipschitz Gauge and Envelope Representation}
\label{sec:saa}
When the gauge set is to measure some type of Lipschitz property, every function adopts an envelope representation, enabling finite-dimensional reformulation.
We begin with the following definition based on a \emph{quasimetric}, a relaxed notion of a metric that may fail to satisfy symmetry.

\begin{definition}[Lipschitz Gauge]
  \label{def:lipgauge}
  A \emph{quasimetric} is a function $c: \Xi \times \Xi \rightarrow \mathbb R$ that satisfies
  \begin{itemize}
    \item Zero-Diagonal: $c(\xi,\xi) = 0$ for all $\xi \in \Xi$.
    \item Nonnegativity: $c(\xi, \xi') > 0$ for all $\xi \neq \xi'$.
    \item Triangle inequality: $c(\xi_1, \xi_2) + c(\xi_2, \xi_3) \geq c(\xi_1, \xi_3)$ for all $\xi_1, \xi_2, \xi_3 \in \Xi$.
  \end{itemize}
The associated $c$-Lipschitz gauge set is defined as 
  $$\mathcal V_c:=\{w \mid w(\xi) - w(\xi') \leq c(\xi, \xi'),\ \forall \xi, \xi' \in \Xi\} = \left\{w ~\middle|~ \sup_{\xi \neq \xi' \in \Xi} \frac{w(\xi) - w(\xi')}{c(\xi, \xi')} \leq 1\right\}.$$
  For a given $(\gamma, s_i, \xi_i)$, we call $\theta_{\gamma, s_i, \xi_i}(\xi):= s_i + \gamma c(\xi, \xi_i)$ an \emph{atomic envelop} associated with $\mathcal V_c$. Given a finite number of atomic envelops $\{\theta_{\gamma, s_i, \xi_i}\}_{i \in [m]}$ that share the same $\gamma$, we define $\hat w_{\gamma, s}:=\min_{i \in [m]} s_i + \gamma c(\xi, \xi_i)$ the associated envelope function. 
  For every $\xi \in \Xi$,  $\theta_{\gamma, s_i, \xi_i}$ is \emph{active} at $\xi$ if $\hat w_{\gamma, s}(\xi) = \theta_{\gamma, s_i, \xi_i}(\xi)$. We say $\theta_{\gamma, s_i, \xi_i}$ is active if it is active at some $\xi \in \Xi$.
\end{definition}

\begin{figure}[tbp]
\centering

\begin{tikzpicture}
\begin{axis}[
    axis lines = middle,
    axis line style={thick},
    xlabel = {},
    ylabel = {},
    xtick=\empty, ytick=\empty,
    clip=false,
    xmin=0, xmax=500,
    ymin=0, ymax=12.5,
    domain=20:480,
    width=11cm, height=5cm,
    samples=200,
    restrict y to domain=0:20
]

\addplot[smooth, thick]
    {sin(x)+2 + 0.002*x};
    \node at (axis cs:490,3.5) [anchor=west] {$w(\xi)=\inf_{\xi'}\theta_{\gamma, w(\xi'), \xi'}(\xi)$};

\pgfmathsetmacro{\gval}{0.02}

\def\xiA{100}
\def\xiB{230}
\def\xiC{420}

\pgfmathsetmacro{\wA}{sin(\xiA) + 2 + 0.002*(\xiA)}
\pgfmathsetmacro{\wB}{sin(\xiB) + 2 + 0.002*(\xiB)}
\pgfmathsetmacro{\wC}{sin(\xiC) + 2 + 0.002*(\xiC)}

\addplot[
    dotted,
    thick,
    domain=20:480,
    samples=200
]
({x},{\wA + \gval*abs(x - \xiA)});

\addplot[
    dotted,
    thick,
    domain=20:480,
    samples=200
]
({x},{\wB + \gval*abs(x - \xiB)});

\addplot[
    dotted,
    thick,
    domain=20:480,
    samples=200
]
({x},{\wC + \gval*abs(x - \xiC)});

\node[circle,fill,inner sep=1pt] at (axis cs:\xiA,0) {};
\node[anchor=north] at (axis cs:\xiA,0) {$\xi_1$};
\node[circle,fill,inner sep=1pt] at (axis cs:\xiB,0) {};
\node[anchor=north] at (axis cs:\xiB,0) {$\xi_2$};
\node[circle,fill,inner sep=1pt] at (axis cs:\xiC,0) {};
\node[anchor=north] at (axis cs:\xiC,0) {$\xi_3$};

\addplot[
    very thick,
    dashed,
    black
]
({x},{
    min(
        min(
            \wA + \gval*abs(x - \xiA),
            \wB + \gval*abs(x - \xiB)
        ),
        \wC + \gval*abs(x - \xiC)
    )
});

\node at (axis cs:490,5.5) [anchor=west]
{$\hat w(\xi) = \min_i \theta_{\gamma, w(\xi_i), \xi_i}(\xi)$};


\end{axis}
\end{tikzpicture}

%
%
%
%
%
%
%
%

\caption{Illustration of the envelope functions under the setting $\mathcal V^\circ = \mathrm{Lip}_1$.
Given a function $w$ with Lipschitz constant $\gamma$, each atomic envelope is defined as
$\theta_{\gamma, w(\xi_i), \xi_i}(\xi) = w(\xi_i) + \gamma\|\xi - \xi_i\|$, shown as dotted lines centered at each sample $\xi_i$.
Their envelope $\min_{i} \theta_{\gamma, w(\xi_i), \xi_i}(\xi)$ forms an upper approximation of $w$.}
\label{fig:thmsaa}
\end{figure}

Figure~\ref{fig:thmsaa} illustrates the envelope functions.
The following proposition provides basic properties of quasimetrics and the induced Lipschitz gauge sets.

\begin{restatable}{proposition}{lipg}
  \label{prop:lipg}
  For any quasimetric $c$, let $\mathcal V_c$ be the associated Lipschitz gauge. The following holds
  \begin{enumerate}
    \item $\mathcal V_c$ is convex and contains the origin.
    \item Every $w \in L^2(\mathbb P)$ adopts the representation $w(\xi)= \inf_{\xi' \in \Xi} \theta_{\gamma, w(\xi'), \xi'}(\xi)$ for every $\gamma \geq \|w\|_{\mathcal V_c}$.
    \item $\|\alpha + w\|_{\mathcal V_c} = \|w\|_{\mathcal V_c}$ for every constant $\alpha \in L^2(\mathbb P)$.
    \item $\|c(\cdot, \xi)\|_{\mathcal V_c} =1$ for every $\xi \in \Xi$.
    \item If $\theta_{\gamma, s_i, \xi_i}(\xi_j) \leq \theta_{\gamma, s_j, \xi_j}(\xi_j)$, then $\theta_{\gamma, s_i, \xi_i} \leq \theta_{\gamma, s_j, \xi_j}$ pointwise.
    \item If $\theta_{\gamma, s_i, \xi_i}(\xi_j) < \theta_{\gamma, s_j, \xi_j}(\xi_j)$, then $\theta_{\gamma, s_i, \xi_i} < \theta_{\gamma, s_j, \xi_j}$ pointwise.
    \item $\theta_{\gamma, s_i, \xi_i}$ is active if and only if it is active at $\xi_i$.
    \item \text{[SAA Compatibility]} $\hat w_{\gamma, s}(\xi_i) \leq s_i$. Equality holds if $\theta_{\gamma, s_i, \xi_i}$ is active.
    \item \text{[Gauge Compatibility]} $\|\hat w_{\gamma, s}\|_{\mathcal V_c} \leq \gamma$. Equality holds if some $\theta_{\gamma, s_i, \xi_i}$ is active at multiple points. In particular, suppose the cardinality of $\Xi$ is strictly larger than the sample size $m$, then $\|\hat w_{\gamma, s}\| = \gamma$.
  \end{enumerate}
\end{restatable}

For the analysis in the remainder of this section, we consider the following problem that generalizes \eqref{eq:distdual}, and always assume $\mathcal V^\circ$ is a Lipschitz gauge generated by a quasimetric.
\begin{subequations}
  \label{eq:dualgen}
\begin{align}
  \inf_{\alpha \in \mathbb R, w(\cdot)\in L^2(\mathbb P)}~& \alpha + g(\mathbb E_{\mathbb P}[h\circ w]) + \epsilon\|w\|_{\mathcal V^\circ}\label{eq:dualgen1}\\
  \text{s.t.} ~& \alpha + w \geq f_x,\label{eq:dualgen2}
\end{align}
\end{subequations}
where $g$ and $h$ satisfy the following properties:
\begin{itemize}
  \item Both $g$ and $h$ are Lipschitz continuous with Lipschitz constants $L_g$ and $L_h$, respectively.
  \item The composite functional $g \circ h_{\mathbb Q}(w) := g\left(\mathbb E_{\mathbb Q}[h\circ w]\right)$
is nondecreasing with respect to the pointwise ordering of $w$, for every probability measure $\mathbb Q$.
\end{itemize}
Such functions $g$ and $h$ naturally occur when incorporating multiple gauge sets in the design (see case study in Section~\ref{sec:case}). In particular, the original problem \eqref{eq:distdual} is recovered as a special case in which both $g$ and $h$ are the identity function. 
To obtain a tractable reformulation of \eqref{eq:distdual}, we consider the following \emph{envelope reformulation} of \eqref{eq:dualgen} with respect to a given sample set $S=\{\xi_i\}_{i \in [m]}$ sampled from the nominal or the true distributions.
\begin{subequations}
  \label{eq:distdual_sp}
\begin{align}
  \inf_{\gamma \geq 0, \alpha, s}~& \alpha + g\left(\sum_{i \in [m]}h(s_i)/m\right) + \epsilon \gamma \label{eq:distdual_sp00}\\
  \text{s.t.} ~& \theta_{\gamma, s_i, \xi_i} \geq f_x - \alpha, \quad \forall i \in [m]. \label{eq:distdual_sp02}
\end{align}
\end{subequations}
This reformulation approximates $w$ by its finite envelope representation $\hat w_{\gamma, s}$ to obtain a semi-infinite program. 
  Suppose $\Xi$ is convex, each $\theta_{\gamma, s_i, \xi_i}$ is piecewise-convex, and $f_x$ is piecewise-concave, then \eqref{eq:distdual_sp02} can be equivalently written as $\inf_{\xi \in \Xi}\{\theta^k_{\gamma, s_i, \xi_i}(\xi) - f^j_x(\xi)\} \geq -\alpha$ for every piece $k$ for $\theta$ and every piece $j$ for $f_x$, allowing the entire problem to be reformulated into a convex optimization problem with finite decision variables and constraints. 
The following lemma provides some results of this formulation when $g$ and $h$ satisfy the required properties.

\begin{restatable}{lemma}{saacomp}
  \label{lem:saacomp}
  Given a feasible solution $(\gamma, \alpha, s)$ of \eqref{eq:distdual_sp} under samples $\{\xi_i\}_{i \in [m]}$, let $\hat w_{\gamma, s}(\xi)$ be the associated envelope function. Then, $\alpha + \hat w_{\gamma, s}$ is feasible to \eqref{eq:dualgen}. Moreover, if an optimal solution exists, there must be some optimal $(\gamma, \alpha, s)$ such that $\hat w_{\gamma, s}(\xi_i) = s_i$ for all $i \in [m]$.
\end{restatable}

We call this type of optimal solution \emph{non-redundant}. The following theorem characterizes the approximation gap between \eqref{eq:dualgen} and \eqref{eq:distdual_sp}.

\begin{restatable}{theorem}{saa}
  \label{thm:saa}
Suppose $\mathcal V^\circ$ is a Lipschitz gauge induced by a quasimetric $c$.
Let $(\alpha^\star,w^\star)$ and $z^\star$ denote an optimal solution and the optimal value of~\eqref{eq:dualgen}.
For a given set of i.i.d.\ samples $\{\xi_i\}_{i=1}^m$, let $(\gamma_m,\alpha,s)$ and $z_m$ denote an optimal solution and the optimal value of~\eqref{eq:distdual_sp}.
Let $\bar{\mathbb P}_m := \frac{1}{m}\sum_{i=1}^m \delta_{\xi_i}$ denote the empirical measure, and let 
$W_1^{\bar c}$
denote the type-1 Wasserstein distance induced by the transport cost $\bar c(\xi, \xi') = \max\{c(\xi,\xi'), c(\xi',\xi)\}$.
Then the following bound holds:
\[
-L_g\left|\iprod{h\circ w^\star,\bar{\mathbb P}_m-\mathbb P}\right|
\le
z^\star - z_m
\le
L_g L_h \gamma_m W_1^{\bar c}(\bar{\mathbb P}_m,\mathbb P).
\]
In particular, $\limsup_{m \to \infty}z_m \leq z^\star$ almost surely. Moreover, if $W_1^{\bar c}(\bar{\mathbb P}_m,\mathbb P)\to 0$ almost surely as $m\to\infty$ and $\limsup_{m \to \infty} \gamma_m < \infty$, then $z_m \xrightarrow{\mathrm{a.s.}} z^\star$.
\end{restatable}

\begin{remark}
Since $c$ is a quasimetric, the symmetrized cost $\bar c$ defines a metric.
Under mild regularity conditions, both the SAA estimation error (lower bound) and the Wasserstein distance
$W^{\bar c}_1(\bar{\mathbb P}_m,\mathbb P)$ (upper bound) vanishes almost surely as the empirical measure
$\bar{\mathbb P}_m$ converges to $\mathbb P$.
Theorem~\ref{thm:saa} therefore implies that the envelope reformulation \eqref{eq:distdual_sp}
is a consistent approximation of \eqref{eq:dualgen} when scenarios are sampled from the nominal
model, and likewise consistently approximates the optimal reweighting problem centered at the
true distribution when scenarios are drawn from data.
From this perspective, sampling from the nominal model versus from observed data induces two
distinct robustness mechanisms: in the former, the ambiguity set is centered at a prescribed
nominal and robustness guards against model misspecification (with SAA used only to approximate
its value); in the latter, the true (but unknown) distribution is the conceptual center, samples provide its SAA proxy, and the robustness radius compensates for statistical estimation error.
\end{remark}

The following corollary shows that, when the empirical distribution is taken as the nominal, the reformulation~\eqref{eq:distdual_sp} is exact.
Although this conclusion follows directly from the bound in Theorem~\ref{thm:saa}, we provide a constructive proof in the appendix for completeness.

\begin{restatable}{corollary}{discssa}
  \label{coro:saa}
Let $z^\star$ and $z_m$ denote the optimal values of \eqref{eq:dualgen} and \eqref{eq:distdual_sp}, respectively. If the nominal $\mathbb P$ is taken as the empirical measure 
$\bar{\mathbb P}_m := \tfrac{1}{m}\sum_{i \in [m]} \delta_{\xi_i}$, then $z_m = z^\star$.
\end{restatable}

\begin{example}[$\mathcal V^\circ = \text{Lip}_1$]
  In \citep{mohajerin2018data}, the nominal distribution is the empirical distribution $\bar {\mathbb P}_m$ and the polar gauge set is $\text{Lip}_1$. The atomic envelope is $\theta_{\gamma, s, \xi'}(\xi):= s + \gamma \|\xi - \xi'\|$. Hence, the corresponding tractable convex reformulation is obtained by dualizing the constraints in \eqref{eq:distdual_sp02}. By Corollary~\ref{coro:saa}, this reformulation is exact when the empirical distribution is taken as the nominal.
\end{example}

The following example shows that $\mathrm{Osc}_1$ (see Example~\ref{eg:rng}) is also a Lipschitz gauge with respect to the discrete metric $
c(\xi,\xi') := \mathbb I(\xi \neq \xi')$,
which assigns unit distance to any pair of distinct points.

\begin{example}[$\mathcal V^\circ = \text{Osc}_1$]
  \label{eg:rng1}
  In this case, the primal gauge is $\mathcal V_{\mathrm{TV}}$ and the associated dual problem is  
\[
\inf_{w(\cdot)} \Big\{\, \mathbb{E}[w] + \epsilon \|w\|_{\mathrm{Osc}_1} \ \big|\  w \ge f_x \,\Big\},
\]
where $\mathrm{Osc}_1 = \mathcal V_{\mathrm{TV}}^\circ := \{w \mid \sup_{\xi \in \Xi} w(\xi) - \inf_{\xi \in \Xi} w(\xi) \le 2\}$
denotes the unit oscillation ball.  
Every function \(w \in L^2(\mathbb P)\) with oscillation \(\gamma\) can be expressed via the envelope representation $w(\xi) = \inf_{\xi'} [w(\xi') + \gamma \mathbb{I}(\xi \neq \xi')]$,
where the binary metric \(\mathbb{I}(\xi \neq \xi')\) equals \(0\) when \(\xi = \xi'\) and \(1\) otherwise. Clearly, the function $\mathbb I$ satisfies all three properties of a quasimetric, thus $\mathrm{Osc}_1$ is indeed a Lipschitz gauge.
Accordingly, each atomic envelope takes the form $\theta_{\gamma,s,\xi'}(\xi) = s + \gamma\, \mathbb{I}(\xi \neq \xi')$.  
Substituting this envelope into \eqref{eq:distdual_sp} yields the following SAA reformulation:
$$
\begin{aligned}
  \inf_{\gamma \geq 0, s_i}~& \sum_{i \in [n]}\frac{s_i}{n} + \frac{\epsilon}{2} \gamma \\
  \text{s.t.} ~& s_i \geq f_x(\xi_i), \quad \forall i \in [n] \\
  ~& s_i + \gamma \geq \sup_{\xi \in \Xi} f_x(\xi), \quad \forall i \in [n].
\end{aligned}
$$
As discussed earlier, when \(f_x\) is piecewise-concave, the supremum term in the last constraint admits a dual representation, yielding a tractable convex reformulation. By Corollary~\ref{coro:saa}, this reformulation is also exact when the empirical distribution is taken as the nominal.
\end{example}

\section{Case Study}
\label{sec:case}
This section illustrates the proposed framework using the illustrative example in the introduction. We derive two tractable reformulations under multiple combined gauge sets and conduct a computational study to demonstrate the resulting formulations. 
We stress that the purpose of this case study is not to benchmark robustness paradigms. Rather, it is designed to illustrate the flexibility of the proposed framework, including its ability to accommodate customized robustness specifications and to support multiple reformulation strategies.

\subsection{Tractable Reformulations}
\label{sec:casereform}
For simplicity, the city region is assumed to be a two-dimensional box $\Xi = [l,u] \subseteq \mathbb R^2$, 
partitioned into finite box-shaped districts 
\(\Xi_k = \{\xi \in \Xi \mid l_k \le \xi \le u_k\}\) for $k \in K$
that may share boundaries but have no overlapping interiors.
The objective is to determine the location of an emergency response center within 
$x \in \Xi$ to minimize the expected distance 
$\mathbb{E}[\|x - \xi\|_1]$ to a random incident $\xi$, where distance is measured using 
the Manhattan metric $\|\cdot\|_1$. Following the same requirement as introduced in the example, the planner aims to (i) hedge against sampling noise using $\phi$-divergence, (ii) guard against region-wise ambiguity via Wasserstein metric, and (iii) ensure robust performance under tail events via CVaR. 

For maximal robustness, we define $\mathcal V_{\text{Comb}}$ as the Minkowski sum of the 
divergence and the region-aware Wasserstein gauge sets. By Corollary~\ref{coro:multisum} and 
Theorem~\ref{thm:gcomp}, this leads to the reformulation~\eqref{eq:caseprob}. 
We now further simplify this expression using results from the previous sections.
First, by Proposition~\ref{prop:cvar}, the minimizer of $w_2$ satisfies 
$w_2 = (f_x - \alpha_2)_+$, and therefore
\[
\|w_2\|^\circ_{\mathcal V_{\text{CVaR}}}
= \frac{\beta}{1-\beta}\,\mathbb E_{\mathbb P}\big[(f_x - \alpha_2)_+\big].
\]
Next, utilizing the $\chi^2$-divergence with $\phi(\nu) = (\nu - 1)^2$ \cite{ben2013robust}, we have
\[
\mathbb E[\phi(\nu)] = \langle \nu - 1, \nu - 1\rangle = \|\nu - 1\|_2^2.
\]
The associated gauge set is therefore the $L^2(\mathbb P)$ unit ball, $\mathcal V_{\phi} = \{\nu \mid \|\nu\|_2 \le 1\}$, which is self-dual. Consequently, $\|w_1\|_{\mathcal V_{\phi}^\circ} = \|w_1\|_2 = \sqrt{\mathbb E_{\mathbb P}[w_1^2]}$. To achieve the region-wise Wasserstein metric, we adopt the design in Example~\ref{eg:hdro} with $\epsilon_i$ as the type-1 Wasserstein radius over region $\Xi_k$. Combining these components yields the following reformulation, where we denote $w := w_1$ for simplicity.
\begin{equation}
  \label{eq:casemain}
\begin{aligned}
  \inf_{\alpha_1,\,\alpha_2,\, w(\cdot)} \quad &
  \alpha_1 + \alpha_2 
  + \mathbb E_{\mathbb P}[w]
  + \delta \sqrt{\mathbb E_{\mathbb P}[w^2]}
  + \sum_{k \in K} \epsilon_k \|w \cdot \mathbb I_{\Xi_k}\|_{\mathrm{Lip}_1^k}
  + \frac{\beta}{1-\beta}\,\mathbb E_{\mathbb P}[(f_x - \alpha_2)_+] \\[4pt]
  \text{s.t.} \quad &
  \alpha_1 + w \ge (f_x - \alpha_2)_+.
\end{aligned}
\end{equation}
From here, we apply two different tractable reformulation methods introduced in Section~\ref{sec:comp}.

\subsubsection{Reformulation via Functional Parameterization}
According to Theorem~\ref{thm:paradual}, we can parameterize the functional space over each $\Xi_k$ using a distinct basis $\phi$ with preserved robustness. In particular, we adopt the moment basis with the conic parameter space $\mathbb R \times \mathbb R^n \times \mathbb S_+^n$ so that the polynomial $h_k(\xi) := q_{0k}+ \iprod{q_k, \xi} + \iprod{Q_k, \xi\xi^\intercal}$ is the functional used for upper approximation over each $\Xi_k$ with some $Q_k \succeq 0$. Due to this explicit format, its Lipschitz is $\sup_{\xi \in \Xi_k}\|q_k + 2Q_k \xi\|_*$ where $\|\cdot\|_*$ is the dual norm of the given norm for Lipschitz. Let $p_k, \mu_k, \Sigma_k$ be the probability mass, conditional expectation, and conditional covariance matrix over each $\Xi_k$, define 
$$\bar \mu := (p_k, p_k\mu_k, p_k\mathrm{vec}(\Sigma_k + \mu_k\mu_k^\intercal))_{k \in K},\ \bar q_k := (q_{0k}, q_k, \mathrm{vec}(Q_k)),\ \bar q:= (\bar q_k)_{k \in K}$$ 
as the stacked vectors, where $\mathrm{vec}(\cdot)$ is the vectorization of a matrix. We can express $\mathbb E_{\mathbb P}[w]$ as
$$
  \mathbb E_{\mathbb P}[w] = \sum_{k \in K} p_k(q_{0k} + \iprod{\mu_k, q_k} + \iprod{\Sigma_k + \mu_k\mu_k^\intercal, Q_k}) = \iprod{\bar \mu, \bar q}.
$$
For $\mathbb E_{\mathbb P}[w^2]$, we define $\Lambda_k$ as the conditional expectation of the matrix $(1, \xi, \mathrm{vec}(\xi\xi^\intercal))^{\otimes 2}$ and $\Lambda := \mathrm{diag}\left([\sqrt{p_k}\Lambda_k^{1/2}]_{k \in K}\right)$ be the associated diagonally stacked matrix, leading to
$$
  \mathbb E_{\mathbb P}[w^2] = \sum_{k \in K} p_k (\bar q_k^\intercal \Lambda_k \bar q_k) =\sum_{k \in K} p_k \|\Lambda_k^{1/2} \bar q_k\|_2^2 = \sum_{k \in K}\|\sqrt{p_k} \Lambda_k^{1/2}\bar q_k\|_2^2 = \|\Lambda \bar q\|_2^2,
$$
Then, we obtain the following reformulation with $f_x(\xi) = \|x-\xi\|_1$ expressed as the piecewise-affine function $\max_{d \in \{\pm 1\}^2} \iprod{d, x - \xi}$.
\begin{align*}
  \inf_{\substack{\alpha, \gamma, \eta \\ \bar q=(q_{k0}, q_k, \mathrm{vec}(Q_k))_{k \in K}}} \quad &
  \alpha_1 + \alpha_2 
  +  \iprod{\bar \mu, \bar q}
  + \delta \|\Lambda \bar q\|_2
  + \sum_{k \in K} \epsilon_k \gamma_k
  + \frac{\beta}{(1-\beta)m} \sum_{j \in [m]} \eta_j \\[4pt]
  \text{s.t.} \quad &
  \alpha_1 + q_{0k} + \iprod{q_k, \xi} + \iprod{Q_k, \xi\xi^\intercal} \ge \iprod{d, x-\xi} - \alpha_2, \ \forall k \in K, \xi \in \Xi_k, d \in \{\pm 1\}^2\\
  \quad & \alpha_1 + q_{0k} + \iprod{q_k, \xi} + \iprod{Q_k, \xi\xi^\intercal} \ge 0, \ \forall k \in K, \xi \in \Xi_k\\
  \quad &  \gamma_k \geq \|q_k + 2Q_k \xi\|_* , \quad \forall k \in K,\xi \in \Xi_k\\
  \quad & \alpha_2 + \eta_j \geq \iprod{d, x- \xi_j}, \quad \forall j \in [m], d \in \{\pm 1\}^2\\
  \quad & Q_k \succeq 0, \quad \forall k \in K\\
  \quad & \gamma, \eta \geq 0,
\end{align*}
where $\{\xi_j\}_{j \in [m]}$ are samples generated from $\mathbb P$. The first two semi-infinite constraints can both be reformulated to contain the following minimization with $q:=q_k + d$ and $q:=q_k$, respectively.
$$ \min_{\xi \in \Xi_k, X \succeq \xi\xi^\intercal} ~\iprod{Q_k, X} + \iprod{q, \xi}. $$
Since $Q_k \succeq 0$, this reformulation is exact.
Applying the standard Schur complement and conic duality, we obtain the following dual problem for $\Xi_k := \{\xi \mid l_k \leq \xi \leq u_k\}$.
$$
\begin{aligned}
  \max_{\bar \tau, \underline \tau \geq 0, s} &~ \iprod{l_k, \underline \tau} - \iprod{u_k, \bar \tau} - s / 4\\
  \text{s.t.} &~ \begin{bmatrix}
    Q_k & q + \bar \tau - \underline \tau \\
    (q + \bar \tau - \underline \tau)^\intercal & s
  \end{bmatrix}
   \succeq 0.
\end{aligned}
$$
Since each constraint needs to be dualized independently, we use $\bar \tau^1_{kd}, \underline \tau^1_{kd}, s^1_{kd}$ and $\bar \tau^2_{k}, \underline \tau^2_{k}, s^2_{k}$ to denote the associated dual variables. For the third semi-infinite constraint, we take the 1-norm as the Lipschitz norm, thus the dual norm is $\|q_k + 2Q_k \xi\|_{\infty}= \max_{i \in [n]}\max_{\xi \in \Xi_k}|q_k^i + 2Q_k^i \xi|$.
Then, $\gamma_k$ can be further represented as
$$
\begin{aligned}
  \gamma_k &\geq q_k^i + 2\max_{\xi \in \Xi_k} \iprod{Q_k^i, \xi}, \quad \gamma_k \geq -q_k^i + 2\max_{\xi \in \Xi_k} \iprod{-Q_k^i, \xi}, \quad \forall k \in K, i \in [n],
\end{aligned}
$$
where $Q^i_k$ is the $i$th row of the matrix $Q_k$. Each linear program can be dualized using the definition of $\Xi_k$, where $\bar \pi$ and $\underline \pi$ denote the associated dual variables. Putting everything together, we obtain the following \emph{parametric reformulation}
\begin{align*}
  \inf_{\substack{x, \alpha, \gamma, s, \eta, \bar \tau, \underline \tau, \bar\pi, \underline\pi \\ \bar q=(q_{k0}, q_k, \mathrm{vec}(Q_k))_{k \in K}}} \quad &
  \alpha_1 + \alpha_2 
  +  \iprod{\bar \mu, \bar q}
  + \delta \|\Lambda \bar q\|_2
  + \sum_{k \in K} \epsilon_k \gamma_k
  + \frac{\beta}{(1-\beta)m} \sum_{j \in [m]} \eta_j \\[4pt]
  \text{s.t.} \quad &
  \sum_{i \in [2]}\alpha_i + q_{0k} + \iprod{l_k, \underline \tau^1_{kd}} - \iprod{u_k, \bar\tau^1_{kd}}-s^1_{kd}/4 \ge \iprod{d, x}, \ \forall k \in K, d \in \{\pm 1\}^2\\
  \quad & \begin{bmatrix}
    Q_k & q_k + d + \bar \tau^1_{kd} - \underline \tau^1_{kd} \\
    (q_k + d + \bar \tau^1_{kd} - \underline \tau^1_{kd})^\intercal & s^1_{kd}
  \end{bmatrix} \succeq 0, \ \forall k \in K, d \in \{\pm 1\}^2\\
  \quad & \alpha_1 + q_{0k} + \iprod{l_k, \underline \tau^2_{k}} - \iprod{u_k, \bar\tau^2_{k}}-s^2_{k}/4 \ge 0, \ \forall k \in K\\
  \quad & \begin{bmatrix}
    Q_k & q_k + \bar \tau^2_{k} - \underline \tau^2_{k} \\
    (q_k + \bar \tau^2_{k} - \underline \tau^2_{k})^\intercal & s^2_{k}
  \end{bmatrix} \succeq 0, \ \forall k \in K\\
  \quad &  \gamma_k \geq q_k^i + 2\left(\iprod{u_k, \bar \pi_{ki}^1} - \iprod{l_k, \underline\pi_{ki}^1}\right) , \quad \forall k \in K, i \in [n]\\
  \quad & Q^i_k = \bar \pi^1_{ki} - \underline\pi^1_{ki}, \quad \forall  k \in K, i \in [n]\\
  \quad &  \gamma_k \geq -q_k^i + 2\left(\iprod{u_k, \bar \pi_{ki}^2} - \iprod{l_k, \underline\pi_{ki}^2}\right) , \quad \forall k \in K, i \in [n]\\
  \quad & -Q^i_k = \bar \pi^2_{ki} - \underline\pi^2_{ki}, \quad \forall  k \in K, i \in [n]\\
  \quad & \alpha_2 + \eta_j \geq \iprod{d, x- \xi_j}, \quad \forall j \in [m], d \in \{\pm 1\}^2\\
  \quad & x \in [l, u], \gamma, \eta, \bar \tau, \underline \tau, \bar \pi, \underline\pi \geq 0.
\end{align*}
This yields a semidefinite program with solution robustness guaranteed by Theorem~\ref{thm:paradual}.

\subsubsection{Reformulation via Envelope Representation}
\label{sec:reformenv}
Since each function \(w \cdot \mathbb{I}_{\Xi_k}\) in \eqref{eq:casemain} admits its envelope representation on \(\Xi_k\), 
we perform the following reformulation,
$$
\begin{aligned}
  \inf_{ \gamma\geq 0, \alpha, w(\cdot)} \quad &
  \alpha_1 + \alpha_2 
  + \mathbb E_{\mathbb P}[w]
  + \delta \sqrt{\mathbb E_{\mathbb P}[w^2]}
  + \sum_{k \in K} \epsilon_k \gamma_k
  + \frac{\beta}{1-\beta}\,\mathbb E_{\mathbb P}[(f_x - \alpha_2)_+] \\[4pt]
  \text{s.t.} \quad &
  \alpha_1 + \inf_{\xi' \in \Xi_k} (w(\xi') + \gamma_k \|\xi - \xi'\|) \ge (\|x - \xi\|_1 - \alpha_2)_+, \quad \forall k \in K, \xi \in \Xi_k.
\end{aligned}
$$
For fixed \(\alpha\), this fits \eqref{eq:dualgen} with 
\(h(w) = (w,\,w^2)\) and \(g(a,b) = a + \delta\sqrt{b}\).
Hence, $g \circ h_{\mathbb Q}$ is \emph{not} monotone in general, and the assumptions of Theorem~\ref{thm:saa} and Corollary~\ref{coro:saa} are violated. To recover monotonicity, we restrict the feasible set by imposing the additional constraint $w \ge 0$ in the dual formulation. Under this restriction, $g \circ h_{\mathbb Q}$ becomes monotone, which allows Theorem~\ref{thm:saa} and Corollary~\ref{coro:saa} to be invoked.

This additional constraint admits a natural interpretation from the gauge perspective. Requiring $w \ge 0$ is equivalent to adding the indicator (or, equivalently, the gauge) of the nonnegative cone $\mathcal V_+ := \{v \in L^2(\mathbb P) \mid v \ge 0\}$ to the dual objective. By gauge algebra (Corollary~\ref{coro:multisum}), augmenting the dual problem with this cone gauge corresponds, in the primal, to taking the Minkowski sum of the original primal gauge set with the polar cone $\mathcal V_+^\circ = \{\nu \in L^2(\mathbb P) \mid \nu \le 0\}$. As a result, the admissible primal gauge set is enlarged, and the resulting formulation remains robust in the sense of admitting a superset of the original ambiguity.

With this modification, we can apply the envelope reformulation \eqref{eq:distdual_sp} either using samples drawn from any chosen nominal \(\mathbb P\) or taking the empirical measure \(\bar{\mathbb P}\) as the nominal: the former asymptotically converges to the optimal value of \eqref{eq:casemain} by Theorem~\ref{thm:saa}, while the latter is an exact reformulation of \eqref{eq:casemain} by Corollary~\ref{coro:saa}. For given samples $\{\xi_j\}_{j \in [m]}$, let $J_k:=\{j \in [m] \mid \xi_j \in \Xi_k\}$, we obtain the following semi-definite program
$$
\begin{aligned}
  \inf_{ \gamma \geq 0, s \geq 0, \alpha} \quad &
  \alpha_1 + \alpha_2
  + \frac{1}{m}\sum_{j \in [m]} s_j
  + \frac{\delta}{\sqrt m} \|s\|_2
  + \sum_{k \in K} \epsilon_k \gamma_k
  + \frac{\beta}{(1-\beta)m}\,\sum_{j \in [m]}(f_x(\xi_j) - \alpha_2)_+ \\[4pt]
  \text{s.t.} \quad &
  \alpha_1 + s_j + \gamma_k \|\xi - \xi_j\| \ge \iprod{d, x - \xi} - \alpha_2, \quad \forall k \in K, \xi \in \Xi_k, j \in J_k, d \in \{\pm 1\}^2\\
 \quad &  \alpha_1 + s_j + \gamma_k \|\xi - \xi_j\| \ge 0, \quad \forall k \in K, \xi \in \Xi_k, j \in J_k\\
\end{aligned}
$$
where we represent the 1-norm as the piecewise-affine function as before. Note that the optimization problem in the second constraint is $\min_{\xi \in \Xi_k}\|\xi - \xi_j\|$ with $\xi_j \in \Xi_k$ for each $j \in J_k$, the minimum is attained at $\xi_j$ with a value of $0$. Hence, this constraint reduces to $\alpha_1 + s_j \geq 0$ for all $j \in [m]$.
Then, the first semi-infinite constraint can be reformulated to contain the optimization $\min_{\xi \in \Xi_k}(\gamma_k\|\xi-\xi_j\| + \iprod{d, \xi})$ on the left-hand side, which is a convex minimization over a compact space with strong duality holds. Let $\Xi_k := \{\xi \mid l_k \leq \xi \leq u_k\}$, we obtain the following dual where $\|\cdot\|_*$ is the dual norm of the given norm $\|\cdot\|$ for Lipschitz.
$$
\begin{aligned}
  \max_{\bar \pi \geq 0, \underline \pi \geq 0} &~ \iprod{d, \xi_j} + \iprod{\xi_j - u_k, \bar \pi} + \iprod{l_k -\xi_j, \underline\pi} \\
                                                &~ \gamma_k \geq \|\bar \pi - \underline \pi + d\|_*.
\end{aligned}
$$
Since each dualization is independent, we use $\bar \pi_{kjd}$ and $\underline \pi_{kjd}$ to label the associated dual variables. Then, we obtain the following \emph{envelope reformulation}
\begin{align*}
  \min_{ x, \gamma, \alpha, s, \eta, \bar \pi, \underline \pi} \quad &
  \sum_{i \in [2]}{\alpha_i}
  + \frac{1}{m}\sum_{j \in [m]} s_j
  + \frac{\delta}{\sqrt m} \|s\|_2
  + \sum_{k \in K} \epsilon_k \gamma_k
  + \frac{\beta}{(1-\beta)m}\,\sum_{j \in [m]}\eta_j \\[4pt]
  \text{s.t.} \quad &
  \sum_{i \in [2]}\alpha_i + s_j + \iprod{\xi_j - u_k, \bar \pi_{kjd}} + \iprod{l_k - \xi_j, \underline \pi_{kjd}} \geq \iprod{d, x - \xi_j}, \forall k\in K, j \in J_k, d \in \{\pm 1\}^2\\
  \quad & \gamma_k \geq \|\bar \pi_{kjd} - \underline \pi_{kjd} + d\|_* , \quad  \forall k \in K, j \in J_k, d \in \{\pm 1\}^2\\
  \quad &  \alpha_1 + s_j  \ge 0, \quad \forall j \in [m] \\
 \quad & \alpha_2 + \eta_j \geq \iprod{d, x - \xi_j}, \quad \forall j \in [m], d \in \{\pm 1\}^2\\
 \quad & x \in [l,u], \gamma, s, \eta, \bar \pi, \underline \pi \geq 0.
\end{align*}
When the 1-norm is taken for the Lipschitz, the dual norm $\|\cdot\|_*$ is the infinite norm, leading to a linear program with one second-order conic term $\|s\|_2$ induced by the $\chi^2$-divergence gauge.

\subsection{Computational Study}

\begin{figure}[!tb]
  \centering
  \begin{tikzpicture}
    \node[anchor=south west, inner sep=0] (img) at (0,0)
      {\includegraphics[width=\linewidth]{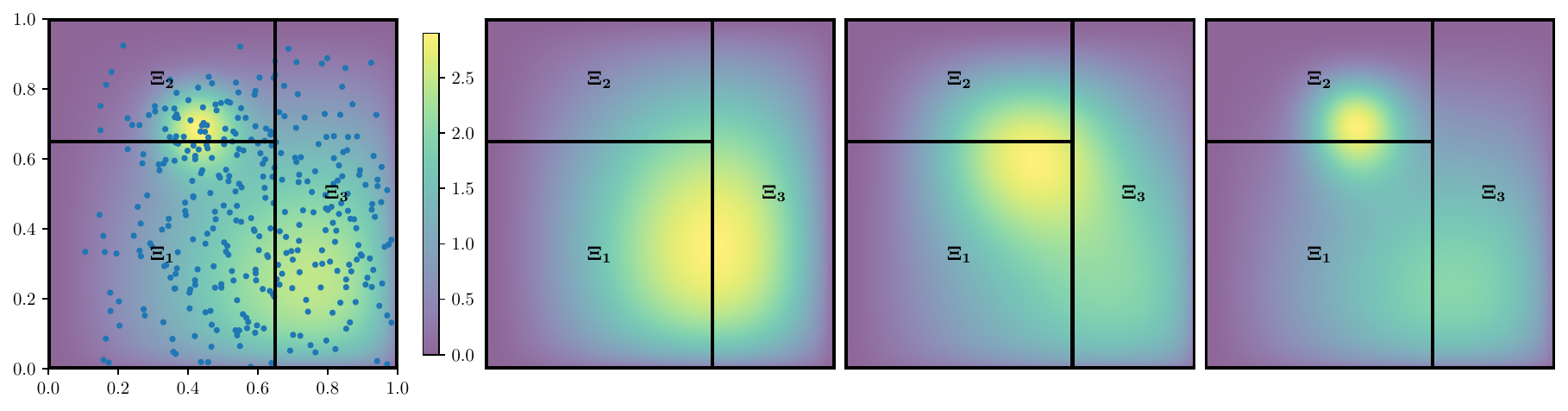}};

    \begin{scope}[x={(img.south east)}, y={(img.north west)}]
      \node[anchor=west, font=\scriptsize] (panelA) at (0.05,-0.02)
        {(a) True distribution $\tilde{\mathbb{P}}$};

      \node[anchor=west, font=\scriptsize] at (0.325,-0.02)
        {(b) Learned $\mathbb P$ ($m = 35$)};

      \node[anchor=west, font=\scriptsize] at (0.554,-0.02)
        {(c) Learned $\mathbb P$ ($m = 342$)};

      \node[anchor=west, font=\scriptsize] at (0.778,-0.02)
        {(d) Learned $\mathbb P$ ($m = 3411$)};
    \end{scope}
  \end{tikzpicture}
  \caption{Instance information. Panel (a) shows the three regions, the true distribution, and the retained $m=342$ observations. Panels (b)--(d) report Bayesian-learned distributions under different data retention levels. We fix $m=342$ and take the distribution in (c) as the nominal $\mathbb P$ for all models, except the data-driven Wasserstein CVaR baseline (WDRO), which uses the empirical measure $\bar{\mathbb P}_m=\tfrac{1}{m}\sum_{i\in[m]}\delta_{\xi_i}$ as the nominal.}
  \label{fig:city}
\end{figure}

\paragraph{Problem Instance.}
We normalize the city region to be $\Xi = [0,1]^2 \subseteq \mathbb R^2$, which is partitioned into three rectangular districts:
\[
\Xi_1 := [0,0.65]\times[0,0.65], \quad
\Xi_2 := [0,0.65]\times[0.65,1], \quad
\Xi_3 := [0.65,1]\times[0,1],
\]
as illustrated in Figure~\ref{fig:city}. The true incident distribution $\tilde{\mathbb P}$ is specified as a mixture of two entrywise independent beta distributions on $\Xi$, with density function given by
\[
\tilde p(\xi)
\;=\;
\sum_{i \in [2]} w_i \prod_{j\in [2]}
\mathrm{Beta}\!\left(\xi_j \mid a_{ij}, b_{ij}\right),
\ w:=(0.9, 0.1), 
\ a := \begin{bmatrix}
  2.5 & 1.5\\  
  15 & 20
\end{bmatrix},
\ b := \begin{bmatrix}
  1.5 & 2.5\\  
  20 & 9.2
\end{bmatrix},
\]
where $w_i$ are the mixture weights, and each $\mathrm{Beta}(\cdot \mid a_{ij},b_{ij})$ denotes the beta density on $[0,1]$ with shape parameters $(a_{ij},b_{ij})$. Figure~\ref{fig:city}a illustrates this distribution via a heatmap. We draw $\tilde m=500$ samples $\{\xi_j\}_{j=1}^{\tilde{m}}$ from $\tilde{\mathbb P}$ to represent historical incident occurrences across the city. 
Due to heterogeneous data retention rates, some observations may be lost. In particular, only $75\%$, $95\%$, and $55\%$ of the observations are preserved in $\Xi_1$, $\Xi_2$, and $\Xi_3$, respectively, resulting in $m_1=154$, $m_2 = 62$, and $m_3 = 126$ retained records, for a total of $m = 342$ observations. Figure~\ref{fig:city}a provides the spatial distribution of these incidents.

\paragraph{Nominal Distribution.}
We assume that the planner has prior knowledge of the parametric form of the incident density $\tilde p(\xi)$ but does not know the hyperparameters $(w_i,a_{ij},b_{ij})$. Given the available historical observations and the additional model-based information, we estimate a nominal distribution by fitting this mixture model in a Bayesian manner. Specifically, we place a Dirichlet prior on the mixture weights $w=(w_1,w_2)$, independent beta priors on the component-wise means, and gamma priors on the corresponding concentration parameters. Posterior inference is performed via automatic differentiation variational inference (ADVI) \citep{kucukelbir2017automatic} using the Python package PyMC, and the resulting nominal distribution $\mathbb P$ is formed by plugging in posterior means of $(w,a,b)$. Figure~\ref{fig:city}b--\ref{fig:city}d shows this learned distribution under different historical datasets. 
In particular, we assume that the retained dataset contains $m=342$ observations.
Accordingly, the fitted distribution shown in Figure~\ref{fig:city}c is taken as the nominal distribution $\mathbb P$ in our proposed reformulations. 
As a baseline, we also report a data-driven type-1 Wasserstein DRO formulation (WDRO). In contrast to the proposed approaches, this data-driven WDRO baseline is inherently tied to a discrete, sample-based center and therefore cannot directly incorporate model-based information by taking a continuous $\mathbb P$ as the nominal. It instead uses the standard empirical measure $\bar{\mathbb P}_m=\tfrac{1}{m}\sum_{i\in[m]}\delta_{\xi_i}$ as its nominal distribution with $m=342$.

\paragraph{Experimental Setting.}
We evaluate the parametric and envelope reformulations derived in Section~\ref{sec:casereform}, denoted by \textsc{PAR} and \textsc{ENV}, against two baselines. The first baseline is the following stochastic CVaR formulation with $\bar m$ i.i.d.\ samples generated from the trained posterior nominal $\mathbb P$:
$$
\begin{aligned}
  \mathrm{STO:} \quad \min_{x \in [l,u],\, \alpha \in \mathbb R,\, \eta \ge 0}
  &~\alpha + \frac{1}{(1-\beta)\bar m}\sum_{j \in [\bar m]} \eta_j \\
  \text{s.t.}\quad
  &~\alpha + \eta_j \ge \iprod{d, x-\xi_j},
    \quad \forall j \in [\bar m],\ \forall d \in \{\pm1\}^2,
\end{aligned}
$$
where $\beta$ denotes the confidence level in CVaR.
The second baseline is the data-driven WDRO model with a CVaR objective under the standard empirical measure $\bar{\mathbb P}_{m}$ as the nominal:
\[
\begin{aligned}
\min_{x \in [l,u],\,\alpha \in \mathbb{R}}
\ & \alpha + \frac{1}{1-\beta}
\inf_{w}\left\{
\mathbb{E}_{\bar{\mathbb{P}}_{m}}[w] + \epsilon\|w\|_{\mathrm{Lip}_1}
 ~\middle|~ \ w(\xi) \ge (\|x-\xi\|_1-\alpha)_+,\ \forall \xi\in\Xi
\right\}.
\end{aligned}
\]
This can be equivalently reformulated following the same steps as in Section~\ref{sec:reformenv}.
\begin{align*}
\mathrm{WDRO:}\quad 
\min_{x, \gamma, \alpha, s, \bar \pi, \underline{\pi}}
&~\alpha + \frac{1}{(1-\beta)m}\sum_{j \in [m]} s_j + \frac{\epsilon}{1-\beta} \gamma\\
  \text{s.t.} &~  \alpha + s_j + \iprod{\xi_j - u, \bar \pi_{jd}} + \iprod{l - \xi_j, \underline \pi_{jd}} \geq \iprod{d, x - \xi_j}, \forall j \in [m], d \in \{\pm 1\}^2\\
              &~ \gamma \geq \|\bar \pi_{jd} - \underline \pi_{jd} + d\|_{\infty} , \quad  \forall j \in [m], d \in \{\pm 1\}^2\\
              &~ x \in [l,u], \gamma, s, \bar \pi, \underline \pi \geq 0.
\end{align*}
We fix $\beta=0.8$ to target the upper $20\%$ tail. For model construction, we draw $\bar m:=2{,}000$ i.i.d.\ samples from the nominal distribution $\mathbb P$ for STO, PAR, and ENV. For PAR and ENV, we set the global $\chi^2$-divergence radius to $\delta:=\rho$, where we treat the scalar $\rho$ as the configuration parameter. We define the region-wise Wasserstein radii in PAR and ENV as $\epsilon_k:=\rho\sqrt{\mathrm{tr}(\Sigma_k)/m_k}$, where $m_k$ and $\Sigma_k$ are the retained sample size and empirical covariance in region $\Xi_k$, respectively.  Likewise, WDRO uses the global radius $\epsilon:=\rho\sqrt{\mathrm{tr}(\Sigma)/m}$ with $m=342$ and global empirical covariance $\Sigma$. Intuitively, smaller sample sizes or larger empirical variability indicate greater statistical uncertainty and therefore motivate larger radii.
To assess out-of-sample performance, we evaluate each computed in-sample solution $x$ on $20$ independently generated test datasets drawn from the true distribution $\tilde{\mathbb P}$, each containing $50{,}000$ observations. For each test set, we compute $\mathrm{CVaR}_{0.8}(x)$ and report the sample mean and standard error across the $20$ replications.
All algorithms are implemented in Python~3.10 using the Mosek~11.1 solver. Experiments are conducted on a MacBook Pro (2023) equipped with an Apple M2~Max processor and $64$~GB of memory.

\paragraph{Hyperparameter $\rho$ Tuning.}
In addition to illustrating how out-of-sample performance varies with $\rho$ across reformulations, we tune $\rho$ separately for each formulation to reflect practical use. Specifically, we select $\rho$ via cross-validation and then compare the resulting out-of-sample performance under the chosen values. The available samples 
are partitioned into five folds; in each split, one fold is held out for testing and the other four are used for training. 
We select the $\rho$ that minimizes the average out-of-sample $\mathrm{CVaR}_{0.8}$ across the five folds, with ties within numerical tolerance $(10^{-6})$ resolved by the smaller $\rho$.

\subsubsection{Performance Analysis}
\begin{figure}[t]
    \centering
    \begin{subfigure}[b]{0.42\textwidth}
        \centering
        \includegraphics[width=\textwidth]{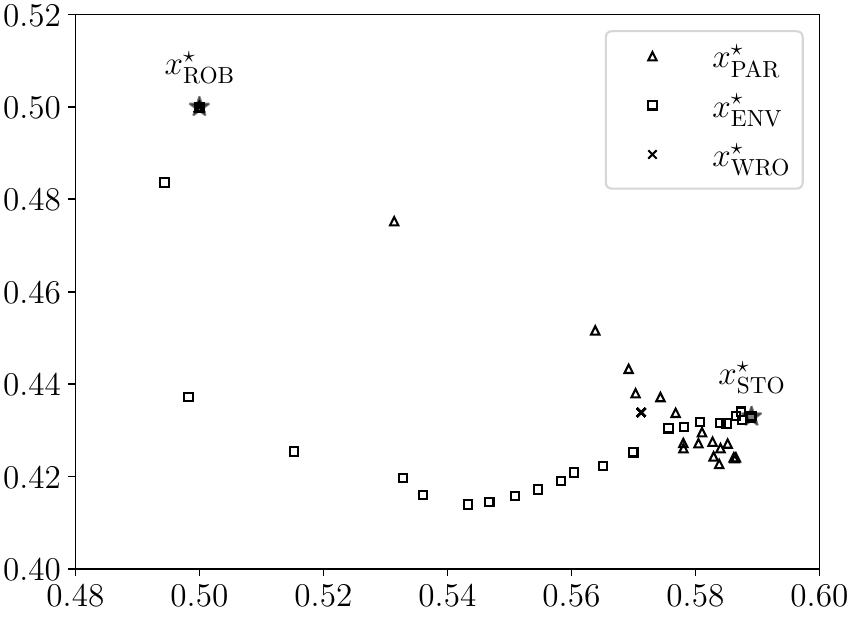}
        \caption{Solution trajectory in the city map.}
        \label{fig:sola}
    \end{subfigure}
    \hfill
    \begin{subfigure}[b]{0.57\textwidth}
        \centering
        \includegraphics[width=\textwidth]{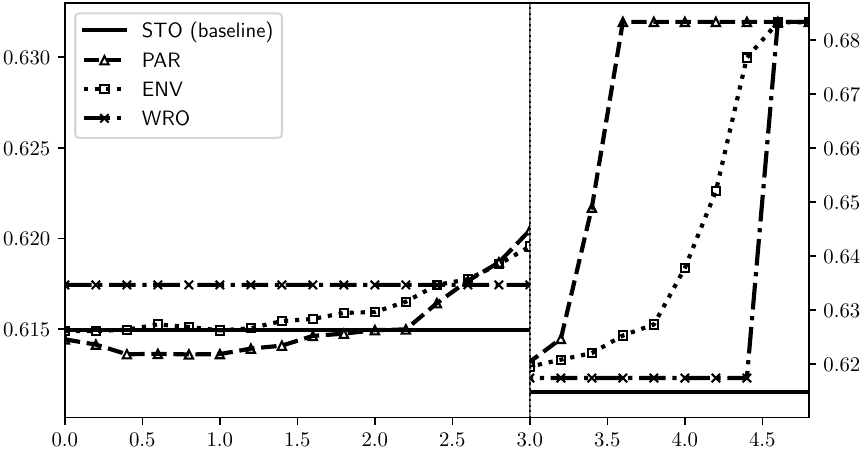}
        \caption{Out-of-sample $\mathrm{CVaR}_{0.8}$ versus $\rho$ (lower is better).}
        \label{fig:solb}
    \end{subfigure}
    \caption{PAR, ENV, and WDRO solution trajectories and out-of-sample $\mathrm{CVaR}_{0.8}$ across configurations. To improve readability, the $\mathrm{CVaR}_{0.8}$ axis is split at $\rho=3$, i.e., values over $\rho\in[0,3]$ are shown on a zoomed vertical scale, while $\rho\in[3,4.8]$ is shown on the full scale.}
    \label{fig:sol}
\end{figure}

Since the support of the true distribution $\tilde{\mathbb{P}}$ is the full square $\Xi=[0,1]^2$, the robust optimization attains its optimum at the center of $\Xi$, namely $x^\star_{\mathrm{ROB}}=(0.5,0.5)$. In contrast, the stochastic CVaR optimization (STO) solution is $x^\star_{\mathrm{STO}}=(0.589,0.433)$. 
Figure~\ref{fig:sola} marks both locations with stars and displays the solution trajectories of PAR, ENV, and WDRO.
As the primal gauge sets expand with increasing $\rho$ within the range $[0,5]$, the solutions produced by all three reformulations move from a neighborhood of $x^\star_{\mathrm{STO}}$ toward the robust solution $x^\star_{\mathrm{ROB}}$, with PAR reaching $x^\star_{\mathrm{ROB}}$ at smaller values of $\rho$.  
In contrast, $x^\star_{\mathrm{WDRO}}$ remains unchanged for most settings and then abruptly jumps to $x^\star_{\mathrm{ROB}}$ at the two largest configurations. 
This is consistent with the linear-program structure of WDRO: as $\epsilon$ varies, the optimal extreme point is stable over intervals and changes only at discrete breakpoints.

Figure~\ref{fig:solb} compares all methods in terms of out-of-sample $\mathrm{CVaR}_{0.8}$ across $\rho \in [0,5]$.
Overall, PAR attains the lowest risk for $\rho \in [0.0,2.4]$, while STO performs best for larger values of $\rho$. ENV slightly improves upon STO only at the first two configurations, which may be attributed to the additional constraint $w\ge 0$ in the envelope reformulation that can introduce extra conservatism.
By contrast, WDRO exhibits a noticeably different $\mathrm{CVaR}_{0.8}$ level from the other three methods, even at $\rho=0$, and remains less favorable than STO across the full range. A plausible explanation is that WDRO is optimized under a different nominal and sampling regime, using the empirical nominal $\bar{\mathbb P}_m$ with $m=342$ observations, whereas STO, PAR, and ENV are constructed using $\bar m=2{,}000$ samples drawn from the learned nominal $\mathbb P$.
At the two largest settings, PAR, ENV, and WDRO yield identical performance as their solutions converge to the robust optimizer $x_{\mathrm{ROB}}^\star$. 

As reported in Table~\ref{tab:expr}, cross-validation selects $\rho_{\mathrm{PAR}}=\rho_{\mathrm{ENV}}=0.8$ and $\rho_{\mathrm{WDRO}}=4.2$, under which both PAR and ENV achieve lower out-of-sample risk than WDRO in this setting. These results suggest that when the assumed distributional family is reasonably aligned with the data-generating process (here, a joint beta model), a fitted nominal distribution can improve out-of-sample performance. 
In terms of runtime, STO is fastest due to its smallest formulation size, whereas ENV is the slowest, consistent with its larger number of sample-dependent variables and constraints.

\begin{table}[t]
  \centering
 \scriptsize
\begin{tabular}{cccrr}
\toprule
Method & Best $\rho$ & Solution & $\mathrm{CVaR}_{0.8}$ (mean $\pm$ se) & Runtime (s) \\\midrule
PAR &	$0.8$ &	$(0.586,0.424)$ &	$0.61359 \pm 2.5 \times 10^{-4}$ &	$0.073$ \\
ENV &	$0.8$ &	$(0.587,0.433)$ &	$0.61511 \pm 2.4 \times 10^{-4}$ &	$1.454$ \\
WDRO &	$4.2$ &	$(0.571,0.434)$ &	$0.61743 \pm 2.4 \times 10^{-4}$ &	$0.073$ \\
STO &	-- &	$(0.589, 0.433)$ &	$0.61492 \pm 2.5 \times 10^{-4}$ &	$0.036$ \\
\bottomrule
\end{tabular}
   
 \caption{Performance comparison with the hyperparameter $\rho$ selected by cross-validation for each method, reporting mean $\pm$ standard error for $\mathrm{CVaR}_{0.8}$.}
  \label{tab:expr}
\end{table}

\section{Conclusion}
\label{sec:conclusion}
This paper introduced a gauge set framework for robustness design in optimization, offering a unified convex-analytic approach for modeling and analyzing robustness across stochastic, robust, and distributionally robust paradigms. 
By formulating robustness through the gauge set reweighting problem, we established quasi-strong duality and showed that the correspondence between primal and dual problems is governed by the geometry of gauge and polar gauge sets. 
This perspective recovers and extends classical results across existing robustness formulations, including moment-based, Wasserstein, and $\phi$-divergence ambiguity sets, while revealing a coherent structure for gauge manipulation through algebraic operations, decomposition, and composition principles. 
To enable computation under continuously supported uncertainty, we further develop two general reformulation schemes that decouple robustness design from reformulation choices, yielding flexible and problem-tailored solution strategies.

The connection between solution robustness and gauge set design opens several promising directions. 
Developing computationally efficient inner approximations of various polar gauges could yield new tractable DRO models with strong robustness guarantees. 
Incorporating structured constraints or hierarchical compositions into primal gauges may further enhance the expressiveness of robustness design in complex applications. 
Overall, shifting attention from ad hoc dual reformulations to the geometric design of gauge sets provides a more flexible framework for customizing robustness in optimization.

\bibliographystyle{plainnat}
\bibliography{bibi/refs}

\newpage
\appendix
\section{Gauge Set Application in Other Robustness Frameworks}
\label{sec:app2}

\subsection{DRO Chance Constraint}
Using gauge sets, we can model the general distributionally robust chance constraints as follows,
\begin{subequations}
\begin{align}
  \min_{x \in \mathcal X} ~& f(x) \\
  \text{s.t.} ~& \left\{
    \begin{aligned}
      \sup_{\nu \in \mathcal R(\mathbb P)} ~&\iprod{\nu, \mathbb I_{g^i_x > 0}}\\
      \text{s.t.} ~& \|\nu - 1\|_{\mathcal V_i} \leq \epsilon_i
    \end{aligned}
  \right\} \leq \beta, \quad \forall i \in [m].
\end{align}
\end{subequations}
where $\mathbb I$ denotes the indicator function of a set (so that $\mathbb I$ flags constraint violation) and $\beta$ is the tolerance level. Moreover, the indicator $\mathbb I_{\{g_x^i>0\}}$ is bounded below and lower semicontinuous whenever $g_x^i$ is, since $\{g_x^i>0\}$ is open. Although quasi-strong duality need not hold, the following dual formulation remains a valid approximation.
\begin{subequations}
\begin{align}
  \min_{x \in \mathcal X} ~& f(x) \\
  \text{s.t.} ~& \alpha_i + \mathbb E[w_i] + \epsilon_i\|w_i\|_{\mathcal V_i^\circ} \leq \beta, \quad \forall i \in [m]\\
  ~& \alpha_i + w_i \geq \mathbb I_{g^i_x > 0}, \quad \forall i \in [m].
\end{align}
\end{subequations}
Then, the gauge set $\mathcal V_i^\circ$ could be designed specifically to capture different types of robustness on the ambiguity of the probability.

\subsection{Robust Satisficing}
\label{sec:sat}
Robust satisficing is another paradigm that optimizes robustness without restricting the scope of ambiguity set \citep{long2023robust}. This method aims to minimize the ratio $(\mathbb E_{\tilde{\mathbb P}}[f_x] -\tau)/d(\tilde{\mathbb P}, \mathbb P)$ where $\tau$ is a given objective target and $d(\tilde{\mathbb P}, \mathbb P)$ signifies a general type of difference between the true probability measure $\tilde{\mathbb P}$ and the empirical measure $\mathbb P$. Since gauge sets provide a general way to specify such differences, we can formulate the general robust satisficing problem as follows.
\begin{subequations}
  \label{eq:sat}
\begin{align}
  \min_{x \in \mathcal X} \inf_{\gamma \geq 0} &~\gamma\\
  \text{s.t.} &~\iprod{f_x, \nu} - \tau \leq \gamma \|\nu - 1\|_{\mathcal V}, ~~\forall \nu \in \mathcal R(\mathbb P). \label{eq:sat1}
\end{align}
\end{subequations}
Using the similar derivation as in Theorem \ref{thm:dist}, we can derive the following reformulation results by rewriting \eqref{eq:sat1} as $\sup_{\nu \in \mathcal R(\mathbb P)} \iprod{f_x, \nu}-\gamma \|\nu - 1\|_{\mathcal V} \leq \tau$.
\begin{align*}
  \min_{x \in \mathcal X}\inf_{w(\cdot)}~& \|w\|_{\mathcal V^\circ} \\
  \text{s.t.} ~& \alpha + \mathbb E[w] \leq \tau,\\
              ~& \alpha + w \geq f_x.
\end{align*}

This reformulation provides a neat dual interpretation for robust satisficing. We again use $\alpha + w$ to upper approximate $f_x$, but with an additional upper bound $\tau$ on the expectation of this approximator. Then, the objective is to minimize the gauge of $w$ under these two constraints. All the previous results regarding different designs of $\mathcal V$ can be carried over to study this robust satisficing problem, facilitating various robustness requirements under this setting. For instance, using Proposition~\ref{prop:phigauge} and Theorem~\ref{thm:gauge}, we can obtain the following robust satisficing dual problem with respect to $\phi$-divergence.
\begin{align*}
  \min_{x \in \mathcal X}\inf_{\gamma \geq 0, w(\cdot)}~& \gamma \\
  \text{s.t.} ~& \alpha + \gamma \mathbb E[\phi^*(w/\gamma)] \leq \tau,\\
              ~& \alpha + w \geq f_x.
\end{align*}

\section{Mathematical Proofs}
\label{sec:app1}
\bipolar*
\pfstart
For Statement 1, when $\epsilon = 0$, $0\nu = 0 \in \mathcal V$, implying $\|0\nu\|_{\mathcal V} = 0$. Otherwise $\epsilon > 0$, then 
$$\|\epsilon \nu\|_{\mathcal V} =  \inf\{t > 0 \mid \nu \in (t/\epsilon)\mathcal V\}=\inf\{\epsilon\beta > 0 \mid \nu \in \beta \mathcal V\} = \epsilon\inf\{\beta > 0 \mid \nu \in \beta \mathcal V\},$$
where the second equality is due to the substitution $\beta:=t / \epsilon$.

For Statement 2, take any $\alpha > \|\nu\|_{\mathcal V}$ and $\beta > \|w\|_{\mathcal V}$. By definition of gauge, $\nu = \alpha \nu_0$ and $w = \beta w_0$ for some $\nu_0, w_0 \in \mathcal V$. Then,
$$\nu + w = (\alpha + \beta) \left(\frac{\alpha}{\alpha + \beta}\nu_0+ \frac{\beta}{\alpha + \beta}w_0\right) \in (\alpha + \beta) \mathcal V.$$
where the membership is due to the convexity of $\mathcal V$. Thus, $\|\nu + w\|_{\mathcal V} \leq \alpha + \beta$. Since this holds for every $\alpha$ and $\beta$, the triangle inequality holds at the infimum.

For Statement 3, since $0$ is an interior of $\mathcal V \cap \mathcal U$ relative to the subspace $\mathcal U$. There exists some $\alpha > 0$ such that the subspace $\alpha$-ball $B_{\mathcal U}(0, \alpha) \subseteq \mathcal V \cap \mathcal U$, i.e., $\alpha u / \|u\| \in \mathcal V$ for every $u \in \mathcal U$. Equivalently, $\|u\|_{\mathcal V} \leq \|u\|/\alpha$ for all $u \in \mathcal U$. Then, the triangle inequality of the gauge function gives
$$\left|\|\nu\|_{\mathcal V} - \|w\|_{\mathcal V}\right| \leq \max\{\|\nu - w\|_{\mathcal V}, \|w - \nu\|_{\mathcal V}\}\leq \frac{\|\nu - w\|}{\alpha}$$
for every $\nu, w \in \mathcal U$ ($(\nu - w) \in \mathcal U$ due to $\mathcal U$ is a subspace), which shows that $\|\cdot\|_{\mathcal V}$ is Lipschitz on $\mathcal U$. Then, for every $\nu \in \mathcal U$, we have 
$\|\nu\|_{\mathcal V} = \|\nu\|_{\mathcal V \cap \mathcal U} = \|\nu\|_{\cl_{\mathcal U}(\mathcal V \cap \mathcal U)}$, since every element in the subspace closure can be approximated by a sequence in $\mathcal V \cap \mathcal U$ with $\|\cdot\|_{\mathcal V}$ acting continuously.

For Statement 4, every $\nu \in \epsilon \mathcal V$ satisfies $\|\nu\|_{\mathcal V} \leq \epsilon$ by definition, proving the first inclusion. For the second inclusion, $\|\nu\|_{\mathcal V} \leq \epsilon$ implies that, for every $\delta > 0$, $\nu/(\epsilon + \delta) \in \mathcal V$ since $\mathcal V$ is convex and contains zero. Then, $\delta \to 0$ provides a sequence in $\mathcal V$ approaching $\nu/\epsilon$, implying $\nu/\epsilon \in \overline{\mathcal V}$.

For Statement 5, the convexity of $\|\cdot\|_{\mathcal V}$ is a direct consequence of Statements 1 and 2. Suppose $\mathcal V$ is further closed, Statement 4 implies $\{\nu \mid \|\nu\|_{\mathcal V} \leq \epsilon\} = \epsilon \mathcal V$ for every $\epsilon > 0$, which is a closed set. For $\epsilon = 0$, the level set becomes $\bigcap_{\epsilon > 0} \epsilon \mathcal V$, where the closedness is preserved through intersection. Since every level set is closed, the function $\|\cdot\|_{\mathcal V}$ is a closed function.

For Statement 6, since $\mathcal V^{\circ\circ}$ contains $\mathcal V$ and is closed, the direction $\overline{\mathcal V}\subseteq \mathcal V^{\circ\circ}$ holds by the definition of closure. For the other direction, suppose $\bar\nu \notin \overline{\mathcal V}$. Since $\overline{\mathcal V}$ is convex-closed and $\{\bar\nu\}$ is compact and convex, the (strong) Hahn-Banach separation theorem along with Riesz representation theorem on $\mathcal L^2(\mathbb P)$ provides some $w$ such that
$$s: = \sup_{\nu \in \overline{\mathcal V}} \iprod{w, \nu} < \iprod{w, \bar\nu}.$$
Pick $\alpha$ such that $s < \alpha < \iprod{w, \bar\nu}$ and define $w':= w/ \alpha$, then
$$\sup_{\nu \in \overline{\mathcal V}} \iprod{w', \nu} = s / \alpha < 1,$$
implying $w' \in \mathcal V^\circ$. However, $\iprod{w', \bar \nu} = \iprod{w, \bar\nu}/ \alpha > 1$, which shows $\bar \nu \notin \mathcal V^{\circ\circ}$.

For Statement 7, the definitions are
$$
\begin{aligned}
  \ker \|\cdot\|_{\mathcal V} &:= \{w \mid w \in \gamma \mathcal V,~\forall \gamma > 0\}\\
  \rec(\mathcal V) &:= \{w \mid \alpha w \in \mathcal V, ~\forall \alpha \geq 0\} = \{w \mid \alpha w \in \mathcal V, ~\forall \alpha > 0\} = \{w \mid w \in (1/\alpha)\mathcal V, ~\forall \alpha > 0\},
\end{aligned}
$$
where we can relax the requirement $\alpha =0$ in $\rec(\mathcal V)$ due to $0w = 0 \in \mathcal V$ by definition.

For Statement 8, we unpack the definitions for every $t > 0$ as
$$w \in t \mathcal V^\circ \Longleftrightarrow \iprod{w/t, \nu} \leq 1 ~\forall \nu \in \mathcal V \Longleftrightarrow \sup_{\nu \in \mathcal V}\iprod{w, \nu} \leq t.$$
Then, we obtain
$$\|w\|_{\mathcal V^\circ} = \inf\{t > 0 \mid w \in t{\mathcal V}^\circ\} = \inf\left\{t > 0 ~\middle|~ \sup_{\nu \in \mathcal V} \iprod{w, \nu} \leq t\right\} = \sup_{\nu \in \mathcal V} \iprod{w, \nu}.$$

For Statement 9, we expand the definitions as follows
$$
\begin{aligned}
  \|w\|_{\mathcal V} &= \inf\{\gamma > 0 \mid w \in \gamma \mathcal V\},\\
  \|w\|_{\mathcal V^\circ} &= \sup_{\nu \in \mathcal V} \iprod{w, \nu}.
\end{aligned}
$$
Then, the first quantity is $0$ whenever $w/\gamma \in \mathcal V$ for every $\gamma > 0$, which implies $\|w\|_{\mathcal V^\circ} \geq \sup_{\gamma > 0}\iprod{w, w}/ \gamma = \infty$.

For Statement 10, $w \in \mathcal V^\perp$ entails $\iprod{\nu, w}= 0$ for all $\nu \in \mathcal V$, implying $\|w\|_{\mathcal V^\circ} = 0$ through the previous statement. On the other hand, if $w \in \mathcal V^\perp$ and $w \in \gamma \mathcal V$ for some $\gamma > 0$, then $\iprod{w, w} = 0$ forcing $w = 0$.
\pfend

\extball*
\pfstart
To show (i), note that every element in $L^2(\mathbb P)$ can be uniquely represented as $\nu + w$ for $\nu \in \mathcal R_0$ and $w \in \mathcal R_0^\perp$ due to the orthogonal decomposition theorem of Hilbert space. Moreover, $\mathcal R_0^\perp$ can be directly computed as the constant functionals $\{\alpha 1 \mid \alpha \in \mathbb R\}$. Then, consider the unit open ball $\mathcal B:=\{\nu + \alpha 1 \mid \|\nu + \alpha 1 \|_2 < 1\}$ under this representation, there exists some sufficiently small $\delta > 0$ so that $\delta \mathcal B \cap \mathcal R_0 \subseteq \mathcal V \cap \mathcal R_0$ since $\mathcal V$ contains an open neighborhood of $0$ inside $\mathcal R_0$. Thus, every $\delta \nu + \delta \alpha 1 \in \delta \mathcal B$ is in $\tilde{\mathcal V}$, i.e., $\tilde{\mathcal V}$ contains an open neighborhood of $0$ in $L^2(\mathbb P)$. For (ii), 
under $\tilde{\mathcal V}$, the gauge constraint in \eqref{eq:gdist} is
$\nu-1 \in \bigcap_{t>\epsilon} t\tilde{\mathcal V}$; since
$\langle 1,\nu\rangle=1$ implies $\nu-1\in\mathcal R_0$ and
$\big((\mathcal V\cap\mathcal R_0)+\mathcal R_0^\perp\big)\cap\mathcal R_0
=\mathcal V\cap\mathcal R_0$, this constraint is equivalent to
$\nu-1 \in \bigcap_{t>\epsilon} t\mathcal V$.
\pfend

\unitight*
\pfstart
Without loss of generality, we assume $\epsilon' = 1$ throughout the proof by replacing $\mathcal V$ by $\epsilon' \mathcal V$ and $\epsilon$ by $\epsilon / \epsilon'$.
For Type-I regularity, we first prove (i) and (ii) for the case $\epsilon = 1$. (i) is provided by Assumption~\ref{asm:reg}. To show uniformly tightness,
for every $\delta > 0$, choose $M > \sup_{\mathbb Q \in \overline{\mathcal P}_{\mathcal V}} \mathbb E_{\mathbb Q}[\Phi(\xi)] / \delta$, and define $\Xi_M:=\{\xi \in \Xi \mid \Phi(\xi) \leq M\}$ to be the level set, which is closed and bounded due to the closedness and coerciveness of $\Phi$. Since $\Xi \subseteq \mathbb R^n$, $\Xi_M$ is compact. For any $\mathbb Q \in \overline{\mathcal P}_{\mathcal V}$, we have
$$\mathbb Q(\Xi \setminus \Xi_M) = \mathbb Q(\Phi > M) \leq \mathbb E_{\mathbb Q}[\Phi(\xi)]/M < \delta.$$
Thus, uniform tightness follows since this is valid for every $\mathbb Q$. Then, for every $\epsilon < 1$, (i) and (ii) trivially hold since $\overline{\mathcal P}_{\epsilon \mathcal V} \subseteq \overline{\mathcal P}_{\mathcal V}$ and both objective finiteness and uniformly tightness are preserved under subsets. 
(iii) is true due to the boundedness of $|f_x|$.

For (iv), by the triangle inequality of the gauge function (Proposition~\ref{prop:bipolar}), we have
$$\|\nu - 1\|_{\tilde{\mathcal V}} = \|\nu - w + (w-1)\|_{\tilde{\mathcal V}} \leq \|\nu -w\|_{\tilde{\mathcal V}} + \|w-1\|_{\tilde{\mathcal V}}.$$
Thus, every $\nu$ such that $\|\nu - w\|_{\tilde{\mathcal V}} \leq \epsilon$ also satisfies $\|\nu - 1\|_{\tilde{\mathcal V}} \leq \epsilon + \|w -1\|_{\tilde{\mathcal V}}$. Since the extended gauge $\tilde {\mathcal V}$ contains a neighborhood of $0$ in $L^2(\mathbb P)$, the function $\|\cdot\|_{\tilde{\mathcal V}}$ is Lipschitz by Statement 3 of Proposition~\ref{prop:bipolar}. Then, when $\|w-1\| < \delta$ for some $\delta > 0$, it ensures $\epsilon + \|w -1\|_{\tilde{\mathcal V}} \leq 1$ since $\epsilon < \epsilon' = 1$ is assumed.
This means, when $\|w-1\| \leq \delta$, the set $\{\nu \mid \|\nu -1 \|_{\tilde {\mathcal V}} \leq 1\}$ fully contains $\{\nu \mid \|\nu - w\|_{\tilde{\mathcal V}} \leq \epsilon\}$. Since the measure closure of the former is $\overline {\mathcal P}_{\tilde {\mathcal V}}$ while the closure of the latter is $\overline{\mathcal P}_{\epsilon \tilde{\mathcal V}, w}$, we also have
$$\overline {\mathcal P}_{\tilde {\mathcal V}} \supseteq \overline{\mathcal P}_{\epsilon \tilde{\mathcal V}, w}.$$
Because (ii) proves the first set is uniformly tight, and uniformly tightness is preserved in subsets, we conclude the proof of (iv).

Since Type-II regularity provides a stronger light-tail condition on $\overline{\mathcal P}_{\mathcal V}$, (i), (ii), and (iv) still hold. For (iii), we have
$$\sup_{\mathbb Q \in \overline{\mathcal P}_{\epsilon\mathcal V}} \mathbb E_{\mathbb Q}[f_x] \leq \alpha + \beta\sup_{\mathbb Q \in \overline{\mathcal P}_{\epsilon\mathcal V}}\mathbb E_{\mathbb Q}[\Phi] \leq \alpha + \beta\sup_{\mathbb Q \in \overline{\mathcal P}_{\epsilon\mathcal V}}\mathbb E_{\mathbb Q}[\Phi^{1+\eta}] < \infty.$$
For (v), we prove the case $\epsilon = 1$. For any $M > 0$, on $\{\xi \mid \Phi(\xi) > M\}$ we have $\Phi \leq \Phi^{1+\eta}/ M^{\eta}$. Hence
$$\sup_{\mathbb Q \in \overline{\mathcal P}_{\mathcal V}} \E_{\mathbb Q}[\Phi \mathbb I_{\Phi > M}] \leq \frac{1}{M^\eta} \sup_{\mathbb Q \in \overline{\mathcal P}_{\mathcal V}} \E_{\mathbb Q}[\Phi^{1+\eta}] \to 0, \text{ as } M \to \infty,$$
where the limit is induced by $\sup_{\mathbb Q \in \overline{\mathcal P}_{\mathcal V}} \E_{\mathbb Q}[\Phi^{1+\eta}] < \infty$ from Assumption~\ref{asm:reg}.
Then, the case $\epsilon < \epsilon' = 1$ directly follows.
\pfend

\valeq*
\pfstart
Under Type-I regularity, $f_x\in C_b(\Xi)$, hence $\mathbb Q\mapsto \langle f_x,\mathbb Q\rangle$
is weak$^\ast$-continuous by the definition of weak convergence.
Under the Type-II regularity, we only prove the weak$^\ast$-lower-semicontinuity, then the weak$^\ast$-upper-semicontinuity follows by a symmetric argument.
Define the truncation $f_x^{M}:=\max\{f_x, -M\}$ for every $M > \max\{\alpha, 0\}$.
Each $f_x^{M}$ is closed and bounded below, hence
\[
\liminf_{n\to\infty}\iprod{f_x^{M},\mathbb Q_n}
\;\ge\;
\iprod{f_x^{M},\mathbb Q}.
\]
Moreover, $f_x^M\downarrow f_x$ pointwise.
To quantify the truncation error, obtain
\[
0\le f_x^{M}(\xi)-f_x(\xi)=(-M-f_x(\xi))_+\le |f_x(\xi)|\,\mathbb I\{f_x(\xi)<-M\}.
\]
Hence, for any $\mathbb Q\in\overline{\mathcal P}_{\epsilon\mathcal V}$, we have
\[
0\le \iprod{f_x^{M},\mathbb Q}-\iprod{f_x,\mathbb Q}
\le \E_{\mathbb Q}\!\left[|f_x(\xi)|\,\mathbb I\{f_x(\xi)<-M\}\right].
\]
Using the Type-II growth bound $|f_x|\le \alpha+\beta\Phi$ (with $\beta > 0$), we obtain
\[
\E_{\mathbb Q}\!\left[|f_x|\,\mathbb I\{f_x<-M\}\right]
\le
\E_{\mathbb Q}\!\left[(\alpha+\beta\Phi)\,\mathbb I\{\alpha+\beta\Phi>M\}\right]
\le
\alpha\,\mathbb Q\!\left(\Phi>\tfrac{M-\alpha}{\beta}\right)
+\beta\,\E_{\mathbb Q}\!\left[\Phi\,\mathbb I\!\left\{\Phi>\tfrac{M-\alpha}{\beta}\right\}\right],
\]
where the first inequality holds due to $|f_x| \leq \alpha + \beta \Phi$ and the indicated set becomes larger.
By Markov's inequality,
\[
\mathbb Q\!\left(\Phi>\tfrac{M-\alpha}{\beta}\right)
\le
\frac{\beta}{M-\alpha}\,\E_{\mathbb Q}[\Phi],
\]
and $\sup_{\mathbb Q\in\overline{\mathcal P}_{\epsilon\mathcal V}}\E_{\mathbb Q}[\Phi]<\infty$
(Lemma~\ref{lem:unitight}) implies that the first term vanishes as $M\to\infty$ uniformly over $\mathbb Q$.
The second term vanishes uniformly by (v) of Lemma~\ref{lem:unitight}.
Combining Portmanteau for $f_x^{M}$ with the uniform truncation bound,
\[
\sup_{\mathbb Q\in\overline{\mathcal P}_{\epsilon\mathcal V}}
\bigl(\langle f_x^{M},\mathbb Q\rangle-\langle f_x,\mathbb Q\rangle\bigr)\to 0,
\]
and letting $M\to\infty$ yields
$\liminf_{n\to\infty}\langle f_x,\mathbb Q_n\rangle\ge \langle f_x,\mathbb Q\rangle$,
i.e., the functional $\mathbb Q \mapsto \iprod{f_x, \mathbb Q}$ is weak$^\ast$-lower-semicontinuous. Weak$^\ast$-upper-semicontinuity follows analogously by applying Portmanteau to the truncations
$g_x^{M}:=\min\{f_x,M\}$ (upper semicontinuous and bounded above).

Finally, since $z_{\epsilon\mathcal V}$ is equivalent to $\sup_{\mathbb Q \in \mathcal P_{\epsilon\mathcal V}} \iprod{f_x, \mathbb Q}$, the ``$\leq$'' direction is trivial. For the other direction, by definition of closure, for every $\mathbb Q$ in the closure, there is a sequence $\mathbb Q_n \in \mathcal P_{\epsilon\mathcal V}$ that weak$^\ast$ converges to $\mathbb Q$. Clearly, $\sup_{\mathbb Q' \in \mathcal P_{\epsilon\mathcal V}} \iprod{f_x, \mathbb Q'} \geq \iprod{f_x, \mathbb Q_n}$ for every $n$, implying 
$$\sup_{\mathbb Q' \in \mathcal P_{\epsilon\mathcal V}} \iprod{f_x, \mathbb Q'} \geq \liminf_{n\rightarrow \infty}\iprod{f_x, \mathbb Q_n} \geq \iprod{f_x, \mathbb Q},$$
where the second inequality is due to the weak$^\ast$-lower-semicontinuity of $\mathbb Q \mapsto \iprod{f_x, \mathbb Q}$.
\pfend

\dist*
\pfstart
  Adopting the conjugate duality framework \citep{rockafellar1974conjugate,bot2009conjugate}, we define the following perturbation function where $h(\nu)$ denotes the function $\|\nu - 1\|_{\mathcal V}$.
  $$F(\nu, u, z) := \begin{cases}
  \iprod{-f_x, \nu}, & \text{if } \nu \geq 0, \iprod{1, \nu} = 1, \text{and } h(\nu - z) - \epsilon \leq u\\
  \infty, & \text{otherwise.}
\end{cases}
$$
Then, the corresponding dual problem can be computed as
\begin{align*}
  \inf\limits_{\gamma, w}F^*(0, -\gamma, -w) & = 
  \inf\limits_{\gamma,w}\sup\limits_{u, z, \nu \geq 0, \iprod{1,\nu} = 1}\left\{ -\gamma u - \iprod{w,z} + \iprod{f_x, \nu} \mid h(\nu-z)-\epsilon \leq u\right\}\\
                      &= \inf\limits_{\gamma\geq 0,w}
  \sup\limits_{z, \nu \geq 0, \iprod{1,\nu}=1} \left\{-\gamma(h(\nu-z)-\epsilon) - \iprod{w, z} + \iprod{f_x, \nu}\right\}\\
                      &= \inf\limits_{\gamma\geq 0,w} \epsilon\gamma + \sup\limits_{\nu \geq 0, \iprod{1,\nu}=1} \left\{ \iprod{f_x,\nu} + \sup\limits_{z}\left\{ - \iprod{w, z} -\gamma h(\nu-z)\right\}\right\}\\
                      &= \inf_{\gamma \geq 0, w}  \epsilon \gamma + \sup\limits_{\nu \geq 0, \iprod{1,\nu}=1} \left\{\iprod{f_x,\nu} + \sup\limits_{z'}\left\{ - \iprod{w, \nu-z'} -\gamma h(z') \right\}\right\}\\
                      &= \inf_{\gamma \geq 0, w} \epsilon \gamma + \sup\limits_{\nu \geq 0, \iprod{1,\nu}=1}\left\{\iprod{f_x-w,\nu} + \sup\limits_{z'} \left\{ \iprod{w, z'} -\gamma h(z') \right\}\right\}\\
                      &= \inf_{\gamma \geq 0, w} \epsilon \gamma + (\gamma h)^*(w) + \sup\limits_{\nu \geq 0, \iprod{1,\nu}=1} \iprod{f_x - w, \nu}\\
                      &= \inf_{\gamma \geq 0, w} \epsilon \gamma + \gamma h^*(w/\gamma) + \sup\limits_{\nu \geq 0} \inf_{\alpha} \left\{\iprod{f_x - w, \nu} + \alpha(1-\iprod{1, \nu})\right\}\\
                      &\leq \inf_{\gamma \geq 0, w} \epsilon \gamma + \gamma h^*(w/\gamma) + \inf_{\alpha} \left\{\alpha + \sup\limits_{\nu \geq 0} \iprod{f_x -\alpha - w, \nu}\right\}\\
                      &= \inf\limits_{\alpha,w, \gamma \geq 0}\left\{\alpha  + \gamma h^*(w/\gamma)+ \epsilon \gamma \mid \alpha + w \geq f_x \right\}.
\end{align*}
Note that the seventh equality holds for the case $\gamma = 0$ under the definition $(0 h)^*(w) = \delta_0(w)$. Then, we compute $h^*(w)$ explicitly as follows.
  \begin{align*}
    h^*(w) &= \sup_{\nu}~\iprod{w, \nu} - \|\nu -1\|_{\mathcal V}\\
           &= \sup_{\nu'}~\iprod{w, \nu'+1} - \|\nu'\|_{\mathcal V}\\
           &= \iprod{1, w} + \sup_{\nu'}~\iprod{w,\nu'} - \|\nu'\|_{\mathcal V}\\
           &= \mathbb E[w] + \delta^{**}_{\mathcal V^\circ}(w)\\
           &= \mathbb E[w] + \delta_{\mathcal V^\circ}(w),
  \end{align*}
  where the fourth equality is by the identity $\|\cdot\|_{\mathcal V} = \delta^*_{\mathcal V^\circ}(\cdot)$ whenever $\mathcal V$ is convex and closed. Then, the dual problem becomes
\begin{align*}
  \inf_{\alpha, w(\cdot)}~& \alpha + \mathbb E[w] + \epsilon \inf \{\gamma \geq 0 \mid w \in \gamma \mathcal V^\circ\}\\
  \text{s.t.} ~& \alpha + w \geq f_x,
\end{align*}
which gives the desired dual formulation by the definition of gauge function. 

According to the Fenchel-Young inequality, the weak duality always holds. For quasi-strong duality, we verify the conditions in Proposition \ref{prop:qduality}. Since the primal \eqref{eq:gdist} is always feasible under $\nu = 1$ with the value $\mathbb E[f_x]$ finite according to Assumption~\ref{asm:reg} and the fact $\mathbb P \in \overline{\mathcal P}_{\epsilon\mathcal V}$, the infimum value function of $F$ is finite at $0$. Moreover, if $F(\nu, u, z) = -\infty$ at some $\nu, u, z$, i.e., $\iprod{f_x, \nu}=+\infty$ at some $\nu, u, z$, it contradicts to Lemma~\ref{lem:unitight}. Hence, $F$ is proper. The convexity of $F$ is also straightforward by our perturbation scheme and the convexity of $h$.

Therefore, it suffices to verify that $\phi(u, z) = \inf_{\nu} F(\nu, u, z)$ is lower semicontinuous at $(0,0)$, i.e., every $(0,0,t)$ that arises as a limit of points from $\epi \phi$ remains in $\epi \phi$. We first note that the parameters $z$ and $u$ are essentially designed to perturb the center $1$ and radius $\epsilon$ of the gauge function. Thus, for every $(u,z)$, we have
$$
\begin{aligned}
  F(\nu, u,z) &= \{\iprod{-f_x, \nu} \mid \nu \in \mathcal R(\mathbb P) \cap \left((1+z) + (\epsilon + u) \mathcal V\right) \}\\
  \phi(u,z)&= \inf_{\nu} F(\nu, u, z)= \inf_{\nu \in  \mathcal R(\mathbb P) \cap \left((1+z) + (\epsilon + u) \mathcal V\right)} \iprod{-f_x, \nu}.
\end{aligned}
$$
Due to Lemma~\ref{lem:extball}, we can safely replace $\mathcal V$ with its extended gauge $\tilde {\mathcal V}$ to preserve the same value. Then, the associated measure closure is $\overline{\mathcal P}_{(\epsilon + u)\tilde{\mathcal V}, 1 + z}$. We define the associated optimization in the measure space as
$$
\begin{aligned}
  \hat F(\mathbb Q, u,z) &= \{\iprod{-f_x, \mathbb Q} \mid \mathbb Q \in \overline{\mathcal P}_{(\epsilon + u)\tilde{\mathcal V}, 1 + z}\}\\
  \hat\phi(u,z)&= \inf_{\mathbb Q} \hat F(\mathbb Q, u, z)= \inf_{\mathbb Q \in \overline{\mathcal P}_{(\epsilon + u)\tilde{\mathcal V}, 1 + z}} \iprod{-f_x, \mathbb Q}.
\end{aligned}
$$
Now, take any convergence sequence $(u_n, z_n, t_n)\to (0,0, t)$ where $(u_n, z_n, t_n) \in \epi \phi$ for every $n$. Since $\phi$ is the infimum of $F$ over $\nu$, $\epi \phi$ is the projection of $\epi F$ onto the space of $(u, z, t)$. By the definition of projection, there exists a sequence $(\nu_n, u_n, z_n, t_n)$ in $\epi F$. By the choice of $\nu_n$, the lifted measures $\nu_n \mathbb P$ belongs to $\overline{\mathcal P}_{(\epsilon + u)\tilde{\mathcal V}, 1 + z}$. Hence, we obtain a sequence $(\nu_n \mathbb P, u_n, z_n, t_n)$ in $\epi \hat F$. According to Lemma~\ref{lem:unitight}, the set $\overline{\mathcal P}_{(\epsilon + u)\tilde{\mathcal V}, 1 + z}$ becomes uniformly tight for all sufficiently small $\epsilon>0$ and all $z$ lying in a sufficiently small neighborhood of $0$.
By Prokhorov’s theorem, uniform tightness ensures precompactness. Since $\overline{\mathcal P}_{(\epsilon + u)\tilde{\mathcal V}, 1 + z}$ is weak$^\ast$-closed by definition, it is weak$^\ast$-compact.
Thus, there is a convergent subsequence in $(\nu_n \mathbb P, u_n, z_n, t_n)$. Passing to this subsequence, we have $(\nu_n \mathbb P, u_n, z_n, t_n) \to (\mathbb Q, 0, 0, t)$ for some $\mathbb Q \in \overline{\mathcal P}_{\epsilon\tilde{\mathcal V}, 1}$. Since $\epi \hat F$ is a closed set as $\iprod{f_x, \cdot}$ is weak$^\ast$-continuous by Lemma~\ref{lem:valeq}, we have $(\mathbb Q, 0, 0, t) \in \epi \hat F$, which implies $(0,0,t) \in \epi \hat \phi$. Finally, notice that when $u=0$ and $z=0$, $-\phi(0,0)$ is the original problem \eqref{eq:gdist}, while $-\hat \phi(0,0) = \sup_{\mathbb Q \in \overline{\mathcal P}_{\epsilon \tilde{\mathcal V}}}\iprod{f_x, \mathbb Q}$. By Lemma~\ref{lem:valeq}, two problems have the same value, i.e., $\phi$ and $\hat \phi$ coincide at $(0,0)$. This shows $(0,0,t) \in \epi \phi$, which proves the lower semicontinuity of $\phi$ at $(0,0)$, and concludes the quasi-strong duality.
\pfend

\crm*
\pfstart
By the definition of this gauge set, \eqref{eq:gdist02} is satisfied for $\epsilon = 1$ if and only if $\nu - 1 \in \mathcal V$, which is equivalent to $\nu \in \tilde{\mathcal Q}$ by the definition of $\mathcal V:=\tilde{\mathcal Q} -1$. Thus, the equivalence holds. When $\tilde{\mathcal Q}$ has the assumed explicit representation, we have
$$\mathcal V^\circ = (\tilde{\mathcal Q} -1)^\circ = \left\{w ~\middle|~ \sup_{g(\nu + 1) \leq 0}\iprod{w, \nu}\leq 1\right\},$$
where $g(\nu + 1)\leq 0$ comes from the shift by $1$. Then, the claimed result follows a direct computation of conjugate duality, and the quasi-strong duality holds by the same proof as in Theorem \ref{thm:dist}.
\pfend

\cmgauge*
We note that the following proof requires a later result Theorem~\ref{thm:gauge}.
\pfstart
  According to the gauge set dual formulation \eqref{eq:gdist}, Proposition~\ref{prop:crm}, and Theorem \ref{thm:gauge}, it suffices to show that the following function 
  $$h(w):=\inf_{\gamma \geq 0} \iprod{\gamma, g(\cdot)}^*(w) - \iprod{1, w}= \sup_{g(\nu+1)\leq 0} \iprod{w, \nu}$$
is positively homogeneous and non-negative. Both are trivially true from the above supremum form and the assumption $g(1) \leq 0$.
\pfend

\cvar*
\pfstart
Since $\tilde{\mathcal Q} = \{\nu \mid \nu \leq (1-\beta)^{-1}\}$, the corresponding $\mathcal V_\beta = \tilde{\mathcal Q} - 1$ has the claimed definition by Proposition~\ref{prop:crm}. To determine the polar set $\mathcal V^\circ_\beta$, we directly compute the following for some input $w$.
  \begin{align*}
    \inf_{\gamma \geq 0}\iprod{\gamma, g(\cdot)}^*(w) &= \inf_{\gamma \geq 0}\left\{\sup_{\nu} \iprod{w, \nu} - \iprod{\gamma, \nu - (1-\beta)^{-1}}\right\}\\
                                                      &= \inf_{\gamma \geq 0}\left\{(1-\beta)^{-1}\iprod{1, \gamma} + \sup_{\nu} \iprod{w -\gamma, \nu}\right\}\\
                                                      &= \begin{cases}
                                                        (1 - \beta)^{-1} \iprod{1, w}, & \text{if } w \geq 0\\
                                                        +\infty, & \text{otherwise.}
                                                      \end{cases}
  \end{align*}
  This proves the definition of the polar set. Then, the definition of the gauge function $\|\cdot\|_{\mathcal V^\circ_\beta}$ follows Theorem \ref{thm:gauge} directly. Hence, problem \eqref{eq:distdual} becomes
\begin{align*}
  \inf_{\alpha, w(\cdot) \geq 0}~& \alpha + (1-\beta)^{-1}\mathbb E[w]\\
  \text{s.t.} ~& \alpha + w \geq f_x.
\end{align*}
Then, $w = (f_x - \alpha)_+$ is an optimal functional for every $\alpha$, which reduces the above formulation to the familiar CVaR optimization.
\pfend

\bdab*
\pfstart
  By definition, $\nu$ is in the kernel if and only if $\nu \in \bigcap_{\epsilon > 0}\epsilon\mathcal V$. When $\mathcal V$ is bounded, every nonzero $\nu$ will be excluded for some sufficiently small $\epsilon$, hence the kernel is $\{0\}$. For the second statement, if $\mathcal V$ is absorbing, then there exists $\epsilon > 0$ such that the open $\epsilon$-$L_2$-ball is contained within $\mathcal V$. Then, every $\nu \in L^2(\mathbb P)$ is contained in the scaled set $(\|\nu\|/\epsilon) \mathcal V$.
\pfend

\absorb*
\pfstart
Recall that $\mathcal V^\circ=\{w\mid \sup_{v\in\mathcal V}\langle w,v\rangle\le 1\}$. Thus $w\in\lambda\mathcal V^\circ$ if and only if $\sup_{v\in\mathcal V}\langle w,v\rangle\le\lambda$. Hence, $\mathcal V^\circ$ is absorbing if and only if for every $w\in L^2(\mathbb P)$ the quantity $\sup_{v\in\mathcal V}\langle w,v\rangle$ is finite.
If $\mathcal V$ is bounded, let $R:=\sup_{v\in\mathcal V}\|v\|_2<\infty$. Then by the Cauchy--Schwarz inequality, $\sup_{v\in\mathcal V}\langle w,v\rangle\le\sup_{v\in\mathcal V}\|w\|_2\|v\|_2\le R\|w\|_2<\infty$, which shows that $\mathcal V^\circ$ is absorbing.
Conversely, assume that $\mathcal V^\circ$ is absorbing. For each $v\in\mathcal V$, define the linear functional $l_v(w):=\langle w,v\rangle$. Absorbingness implies that $\sup_{v\in\mathcal V}|l_v(w)|=\sup_{v\in\mathcal V}|\langle w,v\rangle|<\infty$ for every $w\in L^2(\mathbb P)$. By the Uniform Boundedness Principle, it follows that $\sup_{v\in\mathcal V}\|l_v\|<\infty$. Since $L^2(\mathbb P)$ is a Hilbert space, the operator norm satisfies $\|l_v\|=\|v\|_2$, and hence $\sup_{v\in\mathcal V}\|v\|_2<\infty$, i.e., $\mathcal V$ is bounded.
\pfend

\mdro*
\pfstart
 For the first moment constraint, we have
 \begin{align*}
   &~ (\mathbb E_{\nu \mathbb P}[\xi] - \mu)^\intercal \Sigma^{-1}(\mathbb E_{\nu \mathbb P}[\xi] - \mu) \\
   =&~ \mathbb E[(\nu-1)\cdot \xi]^\intercal Q^\intercal \Lambda^{-1} Q\mathbb E[(\nu-1)\cdot \xi]\\
   =&~ \|\Lambda^{-1/2} Q\mathbb E[(\nu-1)\cdot \id]\|^2_2\\
   =&~ \|\mathbb E[(\nu-1)\cdot \Sigma^{-1/2}]\|^2_2\\
   =&~ \|\nu-1\|^2_{\mathcal V_1}.
 \end{align*}
 The third equality is because $\nu -1$ is a reweighting function and $\Sigma^{-1/2}(\cdot)$ is a random vector (we consider it as a function with input $\xi$). For the second moment constraint, we first subtract $\Sigma$ on both sides then multiply by $\Sigma^{-1/2}$ and $(\Sigma^{-1/2})^{\intercal}$ on the left and right of both sides. Both operations are compatible with the semi-definite inequality given $\Sigma$ is positive-definite. Then, we have
 \begin{align*}
   &~ \Sigma^{-1/2}\left(\mathbb E_{\nu \mathbb P}[(\xi - \mu)(\xi - \mu)^\intercal] - \Sigma\right)(\Sigma^{-1/2})^\intercal\\
   = &~ \mathbb E_{\nu \mathbb P}\left[\Sigma^{-1/2}(\id - \mu)(\id - \mu)^\intercal(\Sigma^{-1/2})^\intercal\right] - \mathbb E\left[\Sigma^{-1/2}(\id-\mu)(\id - \mu)^\intercal(\Sigma^{-1/2})^\intercal\right]\\
   = &~ \mathbb E[(\nu - 1)T_2(\Sigma^{-1/2}(\id-\mu))]\\
   = &~ \mathbb E[(\nu - 1)T_2\circ \Omega_2].
 \end{align*}
 By the same operations, the right-hand side becomes $(\gamma_2 -1)I$. Hence, the semi-definite inequality holds if and only if the largest eigenvalue of the above matrix is bounded by $\gamma_2  - 1$, i.e., the corresponding spectral norm is bounded by $\gamma_2 - 1$, which completes the proof.
\pfend

\mdroi*
\pfstart
To compute the explicit description of $\mathcal V_m^\circ$, we have
 \begin{align*}
   \|\nu\|_{\mathcal V_m} &= \left\|\mathbb E_{\nu \mathbb P}[T_m \circ \Omega]\right\|_{\mathcal N} \\
   &= \sup_{X \in \mathcal N^\circ} \iprod{X, \mathbb E_{\nu \mathbb P}[T_m \circ \Omega]}\\
                                                                           &= \sup_{X \in \mathcal N^\circ} \mathbb E_{\nu \mathbb P}[\iprod{X, T_m \circ \Omega}]\\
                                                                           &= \sup_{w \in \{\iprod{X, T_m \circ \Omega} \mid X \in \mathcal N^\circ\}} \iprod{\nu, w},
 \end{align*}
 where the first equality is the definition of $\mathcal V_m$, the second is by the relationship between gauge set and support function, along with the fact that $\mathcal N$ is convex-closed, the third is due to the linearity of expectation, and the last one is by the definition of expectation in Hilbert space. We also note that the first two inner products are equipped with the corresponding tensor space, and the last one is from the Hilbert space. Since $\|\nu\|_{\mathcal V_m} = \delta^*_{\mathcal V_m^\circ}(\nu)$ whenever $\mathcal V_m$ is convex-closed, we proved the description of $\mathcal V_m^\circ$. 

 Hence, $\mathcal V_m^\circ$ can be considered as the lifting of the norm ball $\mathcal N^\circ$ into the functional space using functions in $T_m \circ \Omega$. Moreover, since $0$ is an interior point of $\mathcal N^\circ$ (since it is a norm ball) and the functional lifting is a surjection onto its image, the zero function $0 = \iprod{0, T_m \circ \Omega}$ is a relative interior of $\mathcal V_m^\circ$, which implies $\cone (\mathcal V_m^\circ) = \mspan(T_m \circ \Omega)$. Given any $w \in L^2(\mathbb P)$, by definition of gauge function, $\|w\|_{\mathcal V_m^\circ} = +\infty$ if $w$ is not within $\cone (\mathcal V_m^\circ)=\mspan(T_m \circ \Omega)$. Otherwise, 
 let $A := [w]_{T_m \circ \Omega}$ denote a coefficient tensor of minimal
$\|\cdot\|_{\mathcal N^\circ}$-norm among all tensors representing $w$.
Then $w \in t\,\mathcal V_m^\circ$ if and only if there exists
$X \in \mathcal N^\circ$ such that
$w = \iprod{tX,\, T_m \circ \Omega}$.
By minimality of $A$, this is equivalent to $A \in t\,\mathcal N^\circ$. Thus, we have
 \begin{align*}
   \| \iprod{A, T_m\circ \Omega}\|_{\mathcal V_m^\circ} &= \inf\left\{t \mid \iprod{A, T_m\circ \Omega} \in t \mathcal V_m^\circ\right\}\\
                                                         &= \inf\left\{t \mid A \in t \mathcal N^\circ\right\} = \|A\|_{\mathcal N^\circ}.
 \end{align*}
To show that $\mathcal V_m^\circ$ induces a pseudonorm, observe that
$\mspan(\mathcal V_m^\circ)=\mspan(T_m\circ \Omega)$ is finite-dimensional.
Since $\mathcal V_m^\circ=\{\iprod{X,T_m\circ \Omega}\mid X\in\mathcal N^\circ\}$ is bounded,
its gauge $\|\cdot\|_{\mathcal V_m^\circ}$ is a norm on $\mspan(T_m\circ \Omega)$ (hence has
trivial kernel there), and equals $+\infty$ outside this subspace.

 Then, the decomposition of the primal gauge set $\mathcal V_m$ is a direct consequence of the later proved gauge set decomposition theorem (Theorem~\ref{thm:decomp}), where the essential part $\mathcal V_m':=\mathcal V_m^\dagger$ (see Theorem~\ref{thm:decomp}) is the polar set of $\mathcal V_m^\circ$ relative to the subspace spanned by $T_m \circ \Omega$. Specifically, we have $\iprod{X', T_m \circ \Omega} \in \mathcal V_m^\dagger$ if and only if $X'$ belongs the following set
 \begin{align*}
   &~\left\{X' ~\middle|~ \sup_{X \in \mathcal N^\circ} \iprod{\iprod{X, T_m\circ \Omega}, \iprod{X', T_m\circ \Omega}}_{\mathbb P} \leq 1\right\}\\
                      =& \left\{X' ~\middle|~ \sup_{X \in \mathcal N^\circ} \sum_{J,J'}X_J X'_{J'}\iprod{[T_m\circ \Omega]_J, [T_m\circ \Omega]_{J'}}_{\mathbb P} \leq 1\right\}\\
                      =& \left\{X' ~\middle|~ \sup_{X \in \mathcal N^\circ} \iprod{X \otimes X', \mathfrak C} \leq 1\right\} = \left\{X' ~\middle|~ \sup_{X \in \mathcal N^\circ} \iprod{\mathfrak C X', X} \leq 1\right\}\\
                      =& \left\{X' ~\middle|~ \mathfrak C X' \in \mathcal N^{\circ \circ} = \mathcal N \right\} = \mathfrak C^{-1}\mathcal N,
 \end{align*}
 where the first equality is by expressing the two functions as linear combinations of basis in $T_m \circ \Omega$; the second and third are by the algebra of tensor product and the fact that $\mathfrak C$ is symmetric; the fourth one is due to $\mathcal N$ is convex-closed; the last one is by the definition of the set inverse operator.
\pfend

\mdroform*
\pfstart
Having moment constraints up to degree $m$ is equivalent to using the intersection of the associated gauge sets. By later proved Corollary \ref{coro:multicon} regarding gauge set intersection, the dual problem immediately becomes \eqref{eq:pweighted3}. By Theorem \ref{thm:mdro1}, each $w_i$ is a function from $\mspan(T_i \circ \Omega_i)$ where $\Omega_i$ is injective. Hence, $w:=\alpha + \sum_{i \in [m]}w_i$ is a polynomial of degree at most $m$. Thus, the constraint \eqref{eq:pweighted3} essentially says using an arbitrary $m$-degree polynomial to upper approximate $f_x$. Then, the first part of the objective penalizes the expectation of this upper approximation, and the second part penalizes the coefficient tensor $[w]_{T_m \circ \Omega_m}$ using the corresponding dual norm induced by $\mathcal N_m^\circ$ according to Theorem \ref{thm:mdro1}.
\pfend

\wdroi*
\pfstart
  By definition, $\text{Lip}_1^\circ = \{\nu \mid \sup_{w \in \text{Lip}_1}\iprod{\nu , w} \leq 1\}$. Hence,
  $$ \text{Lip}_1^\circ + 1 =  \left\{\nu \middle| \sup_{w \in \text{Lip}_1}\iprod{\nu-1 , w} \leq 1\right\} = \{\nu \in \mathcal R(\mathbb P) \mid W_1(\nu \mathbb P, \mathbb P) \leq 1\},$$
  according to Proposition \ref{prop:lip}. Hence, $\text{Lip}^\circ_1$ is the $W_1$ ball centered at $1$ shifted to the center by the translation vector $1$. Since it is known that $W_1$ distance is a metric on the probability simplex, then the shifted set is also a full-dimensional metric ball (for $\epsilon > 0$) restricted to the shifted probability simplex centered at zero. Consequently, every $\nu \in \text{Lip}^\circ_1$ must have a total measure of zero. Then, for every constant function $\alpha \in \mspan(1)$ and every $\nu \in \text{Lip}^\circ_1$, we have $\iprod{\nu, \alpha} = \alpha \iprod{\nu, 1} = 0$, which shows that $\mspan(1)$ is the orthogonal subspace. Hence, $\text{Lip}^\circ_1$ induces a pseudonorm. By Theorem \ref{thm:decomp}, $\text{Lip}_1$ induces a seminorm with $\mspan(1)$ as its kernel.
\pfend

\droiform*
\pfstart
  A direct application of the dual problem \eqref{eq:distdual} gives
\begin{align*}
  \inf_{\alpha, w(\cdot)}~& \mathbb E[\alpha + w] + \epsilon\|w\|_{\text{Lip}_1}\\
  \text{s.t.} ~& \alpha + w \geq f_x.
\end{align*}
By Proposition~\ref{prop:w1}, $\|w+\alpha\|_{\text{Lip}_1} = \|w\|_{\text{Lip}_1}$. Then, replacing $\alpha + w$ with $w$ gives the result.
\pfend

\wpdro*
\pfstart
 By definition, we have
 $$\mathcal V^\circ_{p,\epsilon} = \left\{w \in L^2(\mathbb P)~\middle|~\left\{\begin{array}{rl}
     \sup\limits_{\nu(\cdot), \pi(\cdot,\cdot) \geq 0} & \iprod{w, \nu}\\
     \text{s.t.}  & \iprod{d(\xi, \xi')^p, \pi(\xi, \xi')} \leq \epsilon^p\\
                  & \iprod{1, \pi(\cdot, \xi')} = \nu + 1\\
                  & \iprod{1, \pi(\xi, \cdot)} = 1\\
                  & \nu + 1 \geq 0\\
                  & \iprod{1, \nu} =0.
 \end{array}\right\} \leq 1
 \right\},$$
 where the inner part is a linear program in the Hilbert space $L^2(\mathbb P)$. This problem is clearly feasible by letting $\nu = 0$ and $\pi = 1$, and it is also bounded for every $w \in \mathcal V^\circ_{p,\epsilon}$. In this case, the quasi-strong duality holds by Proposition \ref{prop:qduality} under the standard RHS perturbation. Let $\beta \geq 0, s(\xi), t(\xi'), r(\xi) \geq 0, z$ be the corresponding dual variables in order, the dual problem can be computed as
 \begin{align*}
   \inf\limits_{\beta \geq 0, r(\xi) \geq 0, s(\xi), t(\xi'), z} &~~\epsilon^p \beta - \iprod{1, s} - \iprod{1, t} + \iprod{1, w-s + z}\\
     \text{s.t.}  &~~w-s + z = r\\
                  &~~s(\xi) + t(\xi') \leq \beta d(\xi, \xi')^p.
 \end{align*}
 Note that the first constraint can be reduced to $w+z \geq s$ by eliminating $r\geq 0$. Moreover, the only term of $w+z$ is in the last inner product in the objective. Then, we can set $w(\xi) + z = s(\xi)$, which makes the last inner product equal to zero. We can further replace $s$ by $w + z$ in all occurrences, which gives the following
 \begin{align*}
   \inf\limits_{\beta \geq 0, t(\xi'), z} &~~\epsilon^p \beta - \iprod{1, w} - \iprod{1,t} + z\\
     \text{s.t.}  &~~w(\xi) + z + t(\xi') \leq \beta d(\xi, \xi')^p.
 \end{align*}
 Finally, setting $t(\xi') = \inf_{\xi} \left\{\beta d(\xi, \xi')^p - w(\xi)\right\} - z$ gives us the dual formulation as
 \begin{align*}
 \inf\limits_{\beta \geq 0} &~~\epsilon^p \beta - \iprod{1, w} - \iprod{ 1,\inf_{\xi} \left\{\beta d(\xi, \cdot)^p - w(\xi)\right\}},
 \end{align*}
which proves the claimed polar set definition.
\pfend

\wpdroform*
\pfstart
By Proposition~\ref{prop:wp}, we have
$$\|w\|_{\mathcal V^\circ_{p, \epsilon}}= \inf\limits_{\beta \geq 0} \iprod{1,-w(\cdot) - \inf_{\xi} \left\{\beta (d(\xi, \cdot)^p-\epsilon^p) - w(\xi)\right\}}.$$
Similar to the type-1 case, $\|w\|_{\mathcal V^\circ_{p, \epsilon}}$ is invariant under constant addition. Then, \eqref{eq:distdual} becomes
\[
\inf_{w\ge f_x}\ \inf_{\beta\ge 0}
\left\{\epsilon^p\beta - \iprod{1,\inf_{\xi}\left\{\beta d(\xi,\cdot)^p- w(\xi)\right\}}\right\}.
\]
For each fixed $\beta\ge 0$, the map
$w \mapsto -\iprod{1,\inf_{\xi}\{\beta d(\xi,\cdot)^p-w(\xi)\}}$
is nondecreasing in $w$ (pointwise). Hence, the infimum over $w\ge f_x$
is attained at $w=f_x$, and we obtain
\[
\inf_{\beta\ge 0}\left\{\epsilon^p\beta - \iprod{1,\inf_{\xi}\left\{\beta d(\xi,\cdot)^p-f_x(\xi)\right\}}\right\},
\]
which proves the claim.
\pfend

\phigauge*
\pfstart
  A direct verification shows that $\|\nu - 1\|_{\mathcal V_{\phi,\epsilon}} \leq 1$ if and only if $\mathbb E[\phi(\nu)] \leq \epsilon$, which proves the equivalence. We can compute the polar set using the definition
 $$\mathcal V^\circ_{\phi,\epsilon} = \left\{w ~\middle|~\left\{\begin{array}{rl}
     \sup\limits_{\nu(\cdot)} & \iprod{w, \nu}\\
     \text{s.t.}  & \iprod{1, \phi(\nu + 1)} \leq \epsilon.
 \end{array}\right\} \leq 1
 \right\}.$$
Since $\phi$ is convex-closed, we can use the following perturbation function to compute the dual of the inner optimization.
  $$F(\nu, u, z) := \begin{cases}
    \iprod{-w, \nu}, & \text{if } \iprod{1, \phi(\nu + 1 - z)} - \epsilon \leq u\\
  \infty, & \text{otherwise.}
\end{cases}
$$
Using a similar conjugate duality computation as in Theorem \ref{thm:dist}, we obtain the dual as
$$g(w):=\inf_{\gamma \geq 0} \iprod{1, \gamma(\phi^*(w/\gamma) + \epsilon) - w}.$$
Moreover, by the same argument as in Theorem \ref{thm:dist}, the quasi-strong duality holds since $\phi$ is convex and closed, which concludes the description of the polar set $\mathcal V^\circ_{\phi, \epsilon}$.
\pfend

\divdual*
\pfstart
  We note that the function $g(w)$ that defines $\mathcal V^\circ_{\phi, \epsilon}$ is positively homogeneous 
  since the quasi-strong duality holds in the computation of $\mathcal V^\circ_{\phi, \epsilon}$ and $\sup_{\nu \in \mathcal U} \iprod{\alpha w, \nu} = \alpha \sup_{\nu \in \mathcal U} \iprod{w, \nu}$ for every $\alpha \geq 0$, $w \in L^2(\mathbb P)$, and nonempty $\mathcal U$. It is also non-negative since we have the following when the optimal $\gamma > 0$.
  \begin{align*}
    \inf_{w} g(w) &= \inf_{\gamma \geq 0, w(\cdot)} \gamma\iprod{1, \phi^*(w/\gamma) - w/\gamma + \epsilon}\\
                  &= \inf_{\gamma \geq 0} \epsilon\gamma - \gamma \iprod{1,  \sup_{w} w/\gamma - \phi^*(w/\gamma)}\\
                  &= \inf_{\gamma \geq 0} \epsilon\gamma - \gamma \iprod{1,  \phi^{**}(1)}\\
                  &= \inf_{\gamma \geq 0} \epsilon\gamma - \gamma \iprod{1,  0} = 0.
  \end{align*}
  The third equality is due to $\phi$ being convex and closed, and the fourth is by the property that $\phi(1) = 0$. In the case the optimal $\gamma = 0$, we have $\inf_w g(w) = \inf_w \iprod{1, \delta_0(w) - w} = 0$ by the definition of $0 \phi^*(w/0)$. Hence, $g(w)$ is non-negative. By Theorem \ref{thm:gauge}, $\|w\|_{\mathcal V^\circ_{\phi, \epsilon}} = g(w)$, which gives the claimed reformulation \eqref{eq:phidual} by plugging $g(w)$ into \eqref{eq:distdual}.
  
When $\phi$ is continuously differentiable, the gradient of the objective function with respect to $\nu$ can be computed directly as $w - \phi'(\nu)$, which gives the optimal solution $\nu = (\phi')^{-1}(w)$. This inverse is well-defined since $\phi$ is strictly convex, implying that $\phi'$ is strictly increasing. 
\pfend

\alg*
\pfstart
The first statement is directly from the definition. For the second statement, we first show the ``$\supseteq$'' direction. It suffices to verify every $w \in \conv\left(\bigcup_{i \in I} \mathcal V_i^\circ\right)$ since the set on the left-hand side is closed. Such a $w$ can be represented as some convex combination $w=\sum_{i \in I_n} \lambda_i w_i$ for some finite index subset $I_n \subseteq I$ with $w_i \in \mathcal V_i^\circ$ for every $i \in I_n$. Take an arbitrary $\nu \in \bigcap_{i \in I}\mathcal V_i$, we have 
  $$\iprod{w, \nu} = \sum_{i \in I_n}\lambda_i\iprod{w_i, \nu} \leq \sum_{i \in I_n}\lambda_i = 1,$$
  where the inequality is due to $w_i \in \mathcal V_i^\circ$ and $\nu \in \mathcal V_i$ for every $i$. This completes the proof of this direction. For the other direction, since both sides are convex-closed and contain the origin, we can prove the following equivalent statement invoking Proposition \ref{prop:bipolar}.
  $$\left(\bigcap_{i \in I} \mathcal V_i\right)^{\circ\circ} = \bigcap_{i \in I} \mathcal V_i \supseteq \left(\cl\conv\left(\bigcup_{i \in I} \mathcal V_i^\circ\right)\right)^\circ = \left(\conv\left(\bigcup_{i \in I} \mathcal V_i^\circ\right)\right)^\circ,$$
where the last equality is due to the fact that the polar set automatically ensures the closure property using the intersection of half-spaces. Take $\nu$ from the set on the right, we have $\iprod{\nu, \sum_{i \in I_n}\lambda_i w_i} \leq 1$ for every $I_n \subseteq I$, every convex combination coefficients $\lambda$, and every $w_i \in \mathcal V_i^\circ$. In particular, for every $i \in I$, taking $\lambda_i = 1$ implies $\iprod{\nu, w_i} \leq 1$ for every $w_i \in \mathcal V_i^\circ$, which means $\nu \in \mathcal V_i^{\circ\circ} = \mathcal V_i$ by Proposition \ref{prop:bipolar}. This shows $\nu$ belongs to the intersection of $\mathcal V_i$'s.

  For the third statement, we note that $w \in (\bigoplus_{i \in I} \mathcal V_i)^\circ$ if and only if 
  $$\sup_{\nu \in \bigoplus_{i \in I} \mathcal V_i}\iprod{w, \nu} = \sup_{I_n\subseteq I}\sum_{i \in I_n}\sup_{\nu_i \in \mathcal V_i}\iprod{w, \nu_i} \leq 1,$$
  where $I_n$ is any finite subset of $I$ by the definition of direct sum. Since $0 \in \mathcal V_i$, each summation term is non-negative. Hence, the above inequality is satisfied if and only if $\sup_{\nu_i \in \mathcal V_i} \iprod{w, \nu_i} \leq \lambda_i$ for some $\lambda =(\lambda_i)_{i \in I} \in \Delta$, which is equivalent to $w \in \lambda_i \mathcal V_i^\circ$ for every $i \in I$, i.e., $w \in \bigcap_{i \in I}\lambda_i \mathcal V_i^\circ$ for some $\lambda \in \Delta$. This concludes the proof of this statement.

  The fourth statement is trivial. For the fifth, since $\bigcap_{i \in I} \mathcal V_i\subseteq \mathcal V_i$ for every $i$, we have
  $$\|\nu\|_{\bigcap_{i \in I} \mathcal V_i} \geq \|\nu\|_{\mathcal V_i},~ \forall i \in I $$
  due to the gauge function value being larger for a smaller gauge set. For the other direction, we have $\nu \in \gamma_i\mathcal V_i \text{ for every } \gamma_i > \|\nu\|_{\mathcal V_i}$ by the definition of gauge function. This implies $\nu / \gamma \in \mathcal V_i$ for every $i \in I$ given that $\gamma \geq \sup_{i \in I} \|\nu\|_{\mathcal V_i}$, i.e., $\nu \in \left(\sup_{i \in I} \|\nu\|_{\mathcal V_i}\right) \bigcap_{i \in I} \mathcal V_i$. This concludes this statement.

  For Statement 6, the ``$\leq$'' direction is obvious since $\bigcup_{i \in I} \mathcal V_i \supseteq \mathcal V_i$ for every $i$. We left to show that this inequality cannot be strict. Suppose otherwise $\|\nu\|_{\bigcup_{i \in I}\mathcal V_i} < \gamma' < \gamma:= \inf_{i \in I}\|\nu\|_{\mathcal V_i}$, then $\nu \in \gamma'\bigcup_{i \in I} \mathcal V_i$. That is, there exists some $i \in I$ such that $\nu \in \gamma' \mathcal V_i$, i.e., $\gamma' \geq \|\nu\|_{\mathcal V_i}$. This contradicts that $\gamma$ is the infimum. We note that, in this case, the union is not necessarily convex-closed anymore, but still contains the origin.

 For Statement 7, we have the following 
 \begin{align*}
   \|\nu\|_{\conv\left(\bigcup_{i \in I} \mathcal V_i\right)} &= \inf\left\{\gamma>0 ~\middle|~ \nu \in \gamma \conv\left(\bigcup_{i \in I}\mathcal V_i\right)\right\}\\
                                                              &=\inf\left\{\gamma>0 ~\middle|~ \begin{array}{rl}
                                                                  I_n &\subseteq I\\
                                                                  \nu &= \gamma\sum_{i \in I_n}\lambda_i\nu_i\\
                                                                  \nu_i & \in \mathcal V_i,\quad\forall i \in I_n\\
                                                                  \lambda_i &\geq 0,\quad\forall i \in I_n\\
                                                                  \sum_{i \in I_n}\lambda_i &= 1.
\end{array}\right\}\\
                                                              &=\inf\left\{\gamma>0 ~\middle|~ \begin{array}{rl}
                                                                  I_n &\subseteq I\\
                                                                  \nu &= \sum_{i \in I_n}\gamma\lambda_i\nu_i\\
                                                                  \gamma\lambda_i\nu_i & \in \gamma\lambda_i\mathcal V_i,\quad\forall i \in I_n\\
                                                                  \gamma\lambda_i &\geq 0,\quad\forall i \in I_n\\
                                                                  \sum_{i \in I_n}\gamma\lambda_i &= \gamma.
\end{array}\right\},
 \end{align*}
 where the third equality is obtained by multiplying $\gamma > 0$ on both sides of the constraints. We then substitute $\gamma_i=\gamma\lambda_i$ and $\nu_i' = \gamma_i \nu_i$ to simplify the above formula, which gives
 \begin{align*}
   \|\nu\|_{\conv\left(\bigcup_{i \in I} \mathcal V_i\right)} &= \inf_{I_n \in I,\nu = \sum_{i \in I_n}\nu_i'}\left\{\sum_{i \in I_n}\gamma_i ~\middle|~ \begin{array}{rl}
          \nu_i' & \in \gamma_i\mathcal V_i,\quad\forall i \in I_n\\
          \gamma_i&\geq 0,\quad\forall i \in I_n.
\end{array}\right\}\\
                                                              &= \inf_{I_n \in I,\nu = \sum_{i \in I_n}\nu_i'} \sum_{i \in I_n}\inf\{\gamma_i\geq 0 \mid \nu_i' \in \gamma_i \mathcal V_i\},
 \end{align*}
 where each summand is exactly $\|\nu_i'\|_{\mathcal V_i}$ by definition. This finishes the proof of this statement. 
 Similarly, for the eighth statement, we have
 \begin{align*}
   \|\nu\|_{\bigoplus_{i \in I} \mathcal V_i} &=\inf\left\{\gamma>0 ~\middle|~ \begin{array}{rl}
                                                                  I_n &\subseteq I\\
                                                                  \nu &= \sum_{i \in I_n}\gamma\nu_i\\
                                                                  \gamma\nu_i & \in \gamma\mathcal V_i,\quad\forall i \in I_n
\end{array}\right\}\\
                                              &=\inf_{I_n \subseteq I,\nu=\sum_{i \in I_n}\nu_i'}\left\{\gamma>0 ~\middle|~ \begin{array}{rl}
                                                                  \nu_i' & \in \gamma_i\mathcal V_i,\quad\forall i \in I_n\\
                                                                  \gamma & \geq \gamma_i,\quad \forall i \in I_n
\end{array}\right\},
 \end{align*}
 where we substitute $\nu_i'=\gamma\nu_i$ to obtain the second equality. According to this form, the infimum of $\gamma$  equals $\max_{i \in I_n}\|\nu_i\|_{\mathcal V_i}$, which proves the desired result.

 For the last statement, Since $I$ is finite, we have
 $$\|w\|_{\bigcup_{\lambda \in \Delta} \bigcap_{i \in I}\lambda_i \mathcal V_i} = \inf_{\lambda \in \Delta} \max_{i \in I}\|w\|_{\mathcal V_i}/\lambda_i$$
 by Statements 4--6. We can also safely assume that $\|w\|_{\mathcal V_i} > 0$ for every $i \in I$, since otherwise we can remove the corresponding terms on both sides. Then, the optimal $\lambda$ would make $\|w\|_{\mathcal V_i}/\lambda_i$ equal for every $i \in I$. Otherwise, changing any value $\lambda_i$ would increase the maximum due to $\Delta$ is a simplex. Then, we have
 $$\lambda_j = \frac{\|w\|_{\mathcal V_j}}{\|w\|_{\mathcal V_i}} \lambda_i \Longrightarrow \lambda_i = \frac{\|w\|_{\mathcal V_i}}{\sum_{i \in I}\|w\|_{\mathcal V_i}} \Longrightarrow \frac{\|w\|_{\mathcal V_i}}{\lambda_i} = \sum_{i \in I}\|w\|_{\mathcal V_i}$$
for every $i \in I$, which concludes the proof.
\pfend

\gauge*
\pfstart
  By definition, we have the following thanks to positive homogeneity.
  $$\|w\|_{\mathcal V} = \inf\{\gamma > 0 \mid w = \gamma w', g(w') \leq \epsilon \}=\inf\{\gamma > 0 \mid g(w)/\epsilon \leq \gamma\}.$$
  Then, the non-negativity ensures $\gamma = g(w)/\epsilon = \|w\|_{\mathcal V}$.
\pfend

\multicon*
\pfstart
  In this case, the constraint set \eqref{eq:gdist02} is equivalent to $\|\nu-1\|_{\epsilon_i \mathcal V_i} \leq 1$ for all $i \in [m]$, and is the same as $\|\nu- 1\|_{\bigcap_{i \in [m]} \epsilon_i \mathcal V_i} \leq 1$ by the definition of gauge function. By Theorem \ref{thm:alg} Statement 1 and 2, the polar set is $\conv\left(\bigcup_{i \in [m]}\mathcal V_i^\circ/\epsilon_i\right)$. Then, the claimed result follows the statements 4 and 7 in Theorem \ref{thm:alg}.
\pfend

\multisum*
\pfstart
  In this finite summation case, we have $\mathcal V = \bigoplus_{i \in [m]}\beta_i\mathcal V_i$. Then, by the first and third statements of Theorem \ref{thm:alg}, we have
  $\mathcal V^\circ = \cl\left(\bigcup_{\lambda \in \Delta} \bigcap_{i \in [m]} \lambda_i\mathcal (\beta_iV_i)^\circ\right)$. Let $\mathcal W:= \bigcup_{\lambda \in \Delta} \bigcap_{i \in [m]} \lambda_i\mathcal (\beta_iV_i)^\circ$, by Statements 3 and 9 of Theorem \ref{thm:alg}, we obtain 
  $$\|w\|_{\mathcal W} = \sum_{i \in [m]}\beta_i \|w\|_{\mathcal V_i^\circ}.$$
  This is a convex-closed function due to each $\mathcal V_i^\circ$ is, thus $\|w\|_{\mathcal V^\circ}=\|w\|_{\cl \mathcal W} = \|w\|_{\mathcal W}$, which concludes the proof.
\pfend

\decomp*
\pfstart
Since both $\lin(\mathcal V)$ and $\mathcal V^\perp$ are closed and orthogonal to each other by definition, the subspace $\lin(\mathcal V) \oplus \mathcal V^\perp$ is also closed in $L^2(\mathbb P)$. Then, the decomposition follows the orthogonal decomposition theorem in Hilbert space. To show the decomposition of $\mathcal V$, we write any $\nu \in \mathcal V$ as $\nu_1 + \nu_2 + \nu_3$ from the space decomposition. Then, 
\begin{equation}
  0 = \iprod{\nu, \nu_2} = \iprod{\nu_1, \nu_2} + \iprod{\nu_2, \nu_2} + \iprod{\nu_3, \nu_2} = \|\nu_2\|_2^2 \tag{a}\label{eq:decompa}
\end{equation}
where the first equality is due to $\nu \in \mathcal V$ and $\nu_2$ is from $\mathcal V^\perp$, the second is due to the decomposition, and the third is due to orthogonality between the three spaces. Hence, $\nu_2 = 0$.
Then, $\nu = \nu_1 + \nu_3 \in \mathcal V$ implies $\nu_3 \in \mathcal V - \nu_1 = \mathcal V$, where the equality is due to $\mathcal V$ is invariant under translation of any $\nu_1 \in \lin(\mathcal V)$ by definition of $\lin(\mathcal V)$. Thus, $\nu_3 \in \ess(\mathcal V) \cap \mathcal V$. For the other direction, every $\nu \in \lin(\mathcal V) + \mathcal V^{\dagger}$ can be written as $\nu = \nu_1 + \nu_2$ for some $\nu_1 \in \lin(\mathcal V)$ and $\nu_2 \in \ess(\mathcal V) \cap \mathcal V$, which implies $\nu_1 + \nu_2 \in \mathcal V$ as $\mathcal V$ is invariant under addition of any element in $\lin(\mathcal V)$.

For the decomposition of $\mathcal V^\circ$, any $w \in \mathcal V^\circ$ decomposed as $w=w_1 + w_2 + w_3$ must satisfy
$$\iprod{w_1, \nu_1} + \iprod{w_2, \nu_2} + \iprod{w_3, \nu_3} \leq 1$$
for every $\nu= \nu_1 + \nu_2 + \nu_3 \in \mathcal V$. Then, $w_1$ must be $0$ to ensure boundedness (otherwise, take $\nu_1 = \gamma w_1 \in \lin(\mathcal V)$ and increase $\gamma \to \infty$), and $w_2$ will never affect the summation value since $\nu_2 = 0$ by \eqref{eq:decompa}. The above criterion is then reduced to  
$$\iprod{w_3, \nu_3} \leq 1,~\forall \nu_3 \in \mathcal V^\dagger,$$
implying $w_3 \in \left(\mathcal V^\dagger\right)^\circ_{\ess(\mathcal V)}$ by definition. For the other direction, adding such $w_3$ with any $w_2 \in \mathcal V^\perp$ remains the inequality, which implies any element $w \in \mathcal V^\perp + \left(\mathcal V^\dagger\right)^\circ_{\ess(\mathcal V)}$ also belongs to $\mathcal V^\circ$.

Finally, the subsequent three identities in the statements follow directly by the uniqueness and orthogonality of the two decompositions, and the convexity, closedness, and containing zero are preserved by the intersection that defines $\mathcal V^\dagger$.
\pfend

\gcomp*
\pfstart
We prove by induction. The basic case $m=1$ is true by Theorem~\ref{thm:dist}. Then, the case of $m$ can be written as follows by the induction hypothesis:
$$
\begin{aligned}
  \sup_{\substack{\nu_1 \geq 0 \\ \iprod{1, \nu_1}_{\mathbb P}=1 \\ \|\nu_1 - 1\|_{\mathcal V_1} \leq \epsilon_1}} \inf_{\{\alpha_i, w_i(\cdot)\}_{i=2}^m} &~ \sum_{i=2}^m \left(\alpha_i + \epsilon_i \|w_i\|_{\mathcal V_i^\circ}\right) + \mathbb E_{\nu_1 \mathbb P}[w_2]\\
  \text{s.t.} &~ \alpha_i + w_i \geq w_{i+1}, \quad \forall i \in \{2,3,\dots, m\},
\end{aligned}
$$
where the nominal measure of the inner problem is $\nu_1 \mathbb P$. By Assumption~\ref{asm:reg} and Lemma~\ref{lem:unitight}, $\nu_1 \mathbb P$ still satisfies $\mathbb E_{\nu_1 \mathbb P}[\Phi] < \infty$, thus regularity is preserved, making the above induction step valid. Then, we swap the supremum and infimum to obtain the following with a potential minimax gap,
$$
\begin{aligned}
  \inf_{\{\alpha_i, w_i(\cdot)\}_{i=2}^m} &~ \sum_{i=2}^m \left(\alpha_i + \epsilon_i \|w_i\|_{\mathcal V_i^\circ}\right) + \sup_{\substack{\nu_1 \geq 0 \\ \iprod{1, \nu_1}_{\mathbb P}=1 \\ \|\nu_1 - 1\|_{\mathcal V_1} \leq \epsilon_1}}\iprod{\nu_1, w_2}_{\mathbb P}\\
  \text{s.t.} &~ \alpha_i + w_i \geq w_{i+1}, \quad \forall i \in \{2,3,\dots, m\}.
\end{aligned}
$$
By the definition of $C_{\Phi}(\Xi)$, let $M:=\sup_{\xi\in \Xi} |w_2(\xi)|/(1+\Phi(\xi))< \infty$, then $|w_2| \leq M + M\Phi$, which implies it is Type-II regular with respect to $\Phi$. Thus, Assumption~\ref{asmp} and \ref{asm:reg} are both satisfied.
Thus, applying Theorem~\ref{thm:dist} on the inner supremum finishes the induction step and obtains the dual form.

It remains to verify that the minimax equality holds. Under Assumption~\ref{asm:reg}, the feasible set $\{\nu_1\mathbb P  \mid \nu_1 \geq 0, \iprod{1,\nu_1}_{\mathbb P}=1,\, \|\nu_1-1\|_{\mathcal V_1}\le \epsilon_1\}$ has a weak$^\ast$-compact closure in the measure space by uniform tightness and Prokhorov’s theorem. Furthermore, for each fixed $(\alpha_i, w_i)_{i=2}^m$, the mapping $\nu_1 \mapsto \iprod{ w_2,\nu_1}_{\mathbb P}$ is affine and weak$^\ast$-continuous on the feasible set by Assumption~\ref{asm:reg} and Lemma~\ref{lem:valeq}. Thus, the inner value function is concave and upper semicontinuous in $\nu_1$. It is also convex and lower semicontinuous with respect to $(\alpha_i, w_i)_{i=2}^m$ since it consists of a sum of convex-closed gauge penalties and closed linear inequalities. Therefore, the hypotheses of Sion’s minimax theorem are satisfied, and the supremum and infimum can be interchanged, which completes the proof.
\pfend

\paradual*
\pfstart
By the second statement of Theorem~\ref{thm:alg}, we have
$$\left(\conv\left(\mathcal V \cup \Lambda_\phi^\circ\right)\right)^\circ = \mathcal V^\circ \cap (\Lambda_\phi^\circ)^\circ = \mathcal V^\circ \cap \overline{\Lambda}_\phi.$$
Moreover, $\inf\{t> 0 : w \in t(\mathcal V^\circ \cap \overline{\Lambda}_\phi) = t\mathcal V^\circ \cap t\overline{\Lambda}_\phi = t\mathcal V^\circ \cap \overline{\Lambda}_\phi\}$ since $\overline{\Lambda}_\phi$ is a cone induced by $\Lambda$. Thus, when $w \in \overline{\Lambda}_\phi$, its gauge equals $\|w\|_{\mathcal V^\circ}$, and when $w \notin \overline{\Lambda}_\phi$ it becomes infinite, i.e., we have
$$
\|w\|_{\mathcal V^\circ \cap \overline{\Lambda}_\phi} = 
\begin{cases}
  \|w\|_{\mathcal V^\circ}, & \text{if } w \in \overline{\Lambda}_\phi\\
  \infty, & \text{otherwise.}
\end{cases}
$$
Therefore, when $\Lambda_\phi$ is closed, the problem \eqref{eq:distdual} only allows functions from $\overline{\Lambda}_\phi = \Lambda_\phi$ to be upper approximators, which proves \eqref{eq:paradual} is exactly the dual of \eqref{eq:paraprime}. Otherwise, \eqref{eq:paradual} is the dual problem restricted to the smaller decision space $\Lambda_\phi \subsetneq \overline{\Lambda}_{\phi}$, inducing an upper bound.
\pfend

\lipg*
\pfstart
Convexity can be verified directly. $\mathcal V_c$ contains zero due to $c \geq 0$. Given $\gamma > \|w\|_{\mathcal V_c}$, we have $w \in \gamma \mathcal V_c$, implying $w(\xi) \leq \inf_{\xi'} w(\xi') + \gamma c(\xi, \xi')$.
On the otherhand, taking $\xi' := \xi$ gives
$$\inf_{\xi'} w(\xi') + \gamma c(\xi, \xi') \leq w(\xi) + \gamma c(\xi, \xi) = w(\xi),$$
which proves Statement 2. The third statement can be verified directly by definition. 
For Statement 4, fix $\xi_0\in\Xi$. By the triangle inequality, $c(\xi,\xi_0)\le c(\xi,\xi')+c(\xi',\xi_0)$ for all $\xi,\xi'$, hence $c(\xi,\xi_0)-c(\xi',\xi_0)\le c(\xi,\xi')$, i.e., $c(\cdot,\xi_0)\in\mathcal V_c$ and thus $\|c(\cdot,\xi_0)\|_{\mathcal V_c}\le 1$. Conversely, if $c(\cdot,\xi_0)\in t\mathcal V_c$, then setting $\xi'=\xi_0$ gives $c(\xi,\xi_0)\le t\,c(\xi,\xi_0)$ for all $\xi$, so $t\ge 1$ whenever $c(\xi,\xi_0)>0$ for some $\xi$, and therefore $\|c(\cdot,\xi_0)\|_{\mathcal V_c}=1$.

For Statement 5, the premise implies
$$
  s_i + \gamma c(\xi_j, \xi_i) \leq s_j + \gamma c(\xi_j, \xi_j) \Longrightarrow s_j - s_i \geq \gamma c(\xi_j, \xi_i).
$$
Then, for every $\xi \in \Xi$, we have
$$
\theta_{\gamma, s_j, \xi_j}(\xi) - \theta_{\gamma, s_i, \xi_i}(\xi) = s_j - s_i + \gamma(c(\xi, \xi_j) - c(\xi, \xi_i)) \geq \gamma (c(\xi_j, \xi_i) + c(\xi, \xi_j) - c(\xi, \xi_i)) \geq 0,
$$
where the last is the triangle inequality. Statement 6 can be proved by the same argument. For Statement 7, one direction is trivial. For the other, suppose $\theta_{\gamma, s_i, \xi_i}$ is not active at $\xi_i$, then Statement 6 shows that it is fully dominated in the entire domain by some other atomic envelope. For Statement 8, since $\theta_{\gamma, s_i, \xi_i}(\xi_i) = s_i$ and $\hat w$ is the minimum over all atomic envelopes, the first inequality holds. Suppose $\theta_{\gamma, s_i, \xi_i}$ is active, then it must be active at $\xi_i$ according to Statement 7, which proves the equality. For Statement 9, take any $\xi, \xi' \in \Xi$ and let $\theta_{\gamma, s_i, \xi_i}$ be an active envelope at $\xi'$. We have
$$\hat w_{\gamma, s}(\xi') = s_i + \gamma c(\xi', \xi_i), \quad  \hat w_{\gamma, s}(\xi) \leq  s_i + \gamma c(\xi, \xi_i),$$
which implies
$$\hat w_{\gamma, s}(\xi) - \hat w_{\gamma, s}(\xi') \leq \gamma \left(c(\xi, \xi_i) - c(\xi', \xi_i)\right)\leq \gamma c(\xi, \xi').$$
Since $\xi, \xi'$ are arbitrarily chosen, we have $\|\hat w_{\gamma, s}\|_{\mathcal V_c} \leq \gamma$. For the equality claim, suppose some $\theta_{\gamma, s_i, \xi_i}$ is active at multiple points, then one of them is $\xi_i$ by Statement 7. Take another active point $\xi_j \neq \xi_i$, we have
$$\hat w_{\gamma, s}(\xi_j) - \hat w_{\gamma, s}(\xi_i) = s_i + \gamma c(\xi_j, \xi_i) - s_i = \gamma c(\xi_j, \xi_i),$$
which attains the upper bound $\gamma$. Finally, suppose $\Xi$ has a larger cardinality than $I$, then there exists some $\xi \in \Xi$ does not associated with any atomic envelope. Then, the active atomic at this point, say $\theta_{\gamma, s, \xi_i}$, must be active at both $\xi$ and $\xi_i$.
\pfend

\saacomp*
\pfstart
Since \eqref{eq:distdual_sp02} is equivalently to $\alpha + \min_{i \in [m]}\theta_{\gamma, s_i, \xi_i} = \alpha + \hat w_{\gamma, s} \geq f_x$, the functional $\alpha + \hat w_{\gamma, s}$ is feasible to \eqref{eq:dualgen}. By the definition of envelope, we always have $\hat w_{\gamma, s}(\xi_i) \leq s_i$. 
For the last statement, suppose $\hat w_{\gamma, s}(\xi_i) < s_i$, Statements 6 and 8 in Proposition~\ref{prop:lipg} imply that $\theta_1 := \theta_{\gamma, s_i, \xi_i}$ is strictly dominated by another atomic envelope. Let $\theta_{\gamma, s_j, \xi_j}$ be the active atomic envelope at $\xi_i$. Then, we can modify $\theta_1$ by reducing the constant $s_i \gets \hat w_{\gamma, s}(\xi_i) = \theta_{\gamma, s_j, \xi_j}(\xi_i)$ so that $\theta_1$ becomes active. Moreover, due to Statement 5 in Proposition~\ref{prop:lipg}, this operation will not change the function value of $\hat w_{\gamma, s}$. Since the function $g\circ h_m$ is non-decreasing on $s$, the objective value will not increase by this operation of reducing $s_i$, which retains the optimality of the solution.
\pfend

\saa*
\pfstart
Define $\hat w_{\gamma, s}:=\min_{i \in [m]}\theta_{\gamma, s_i, \xi_i}$ relative to the given samples $S$, we focus on the relationship between the following problem and \eqref{eq:dualgen}.\\
\begin{subequations}
  \label{eq:pf01}
  \begin{align}
    \inf_{\gamma, \alpha, s} &~ \alpha + g\left(\sum_{i \in [m]}h(s_i)/m\right) + \epsilon\gamma\\
                           &~ \hat w_{\gamma, s}\geq f_x - \alpha.
  \end{align}
\end{subequations}
Clearly, \eqref{eq:pf01} is equivalent to \eqref{eq:distdual_sp} by the definition of $\hat w$. For a given $S$ with size $m$, let $z^\star, z_{m}$ be the optimal value of \eqref{eq:dualgen} and \eqref{eq:pf01}, respectively. Let $(\alpha^\star, w^\star)$ be the optimal solution of \eqref{eq:dualgen} (which may be obtained via a weak$^\ast$ convergent sequence), we construct a feasible solution for \eqref{eq:pf01} as $\gamma:=\|w^\star\|_{\mathcal V^\circ}$, $\alpha:= \alpha^\star$, and $s_i:=w^\star(\xi_i)$, which induces the finite-envelope $\hat w_m:=\hat w_{\gamma, s}$. Then, $\hat w_m$ is feasible to \eqref{eq:pf01} due to the envelope property, and its objective value is exactly $\alpha^\star + g\left(\sum_{i \in [m]}h\circ w^\star(\xi_i)/m\right) + \epsilon \|w^\star\|_{\mathcal V^\circ}$ by construction. Then, we have
$$
\begin{aligned}
  z_{m} - z^\star &\leq \left(\alpha^\star + g\left(\sum_{i \in [m]}h\circ w^\star(\xi_i)/m\right) + \epsilon \|w^\star\|_{\mathcal V^\circ}\right) - z^\star \\
                    &= g\left(\sum_{i \in [m]}h\circ w^\star(\xi_i)/m\right) - g\left(\mathbb E[h\circ w^\star]\right) \\
                    & \leq L_g\left|\sum_{i \in [m]}h\circ w^\star(\xi_i)/m - \mathbb E[h\circ w^\star]\right|\\
                    & = L_g \left|\iprod{h \circ w^\star, \bar{\mathbb P}_m - \mathbb P}\right|.
\end{aligned}
$$
Since gauge compatibility cancels the gauge term, we obtain the first equality. The second inequality applies the Lipschitz inequality of $g$. Since this is the SAA estimation error of the random variable $h\circ w^\star$, which is independent of the i.i.d.\ samples, the strong-law-of-large-numbers guarantees $\limsup z_m \leq z^\star$ almost surely.

For the other bound, take any optimal solution $(\alpha, s, \gamma_m)$ from \eqref{eq:pf01} that is non-redundant as in Lemma~\ref{lem:saacomp}, the associated $\alpha + \hat w_{\gamma, s}$ is feasible to \eqref{eq:dualgen}. We have the following
$$
\begin{aligned}
  z^\star - z_m & = \left(\alpha + g\left(\mathbb E[h\circ \hat w_{\gamma_m, s}]\right) + \epsilon \|\hat w_{\gamma_m, s}\|_{\mathcal V^\circ}\right) - \left(\alpha + g\left(\sum_{i \in [m]} \frac{h(s_i)}{m}\right) + \epsilon\gamma_m\right)\\
  & \leq L_g \left|\mathbb E[h\circ \hat w_{\gamma_m, s}] - \sum_{i \in [m]} \frac{h(s_i)}{m}\right|\\
  & = L_g \left|\mathbb E_{\mathbb P}[h \circ \hat w_{\gamma_m, s}] - \mathbb E_{\bar {\mathbb P}_m}[h \circ \hat w_{\gamma_m, s}]\right|\\
  & = L_g \left|\int h\circ \hat w_{\gamma_m, s}(\xi) - h\circ \hat w_{\gamma_m, s}(\xi') ~d \pi(\xi, \xi')\right|, \quad \forall \pi \in \Pi(\bar{\mathbb P}_m, \mathbb P)\\
  & \leq L_g L_h \int \left|\hat w_{\gamma_m, s}(\xi) - \hat w_{\gamma_m, s}(\xi')\right| ~ d\pi(\xi, \xi'), \quad \forall \pi \in \Pi(\bar{\mathbb P}_m, \mathbb P)\\
  & \leq L_g L_h \int \gamma_m \max\{c(\xi, \xi'), c(\xi', \xi)\} ~ d\pi(\xi, \xi'), \quad \forall \pi \in \Pi(\bar{\mathbb P}_m, \mathbb P)\\
  & = L_g L_h \gamma_m \inf_{\pi \in \Pi(\bar{\mathbb P}_m, \mathbb P)}\int \bar c(\xi, \xi') ~ d\pi(\xi, \xi')\\
  & = L_g L_h \gamma_m W_1^{\bar c}(\bar {\mathbb P}_m, \mathbb P).
\end{aligned}
$$
Since gauge compatibility ensures $\|\hat w_{\gamma_m, s}\|_{\mathcal V^\circ} = \gamma_m$, we obtain the second inequality using Lipschitz inequality of $g$.
For the third inequality above, the following holds by the definition of Lipschitz gauge
$$\|\hat w_{\gamma_m, s}\|_{\mathcal V^\circ} = \gamma_m \geq \max\left\{\frac{\hat w_{\gamma_m, s}(\xi) - \hat w_{\gamma_m, s}(\xi')}{c(\xi, \xi')},\frac{\hat w_{\gamma_m, s}(\xi') - \hat w_{\gamma_m, s}(\xi)}{c(\xi', \xi)}\right\},$$
which implies $\max\{c(\xi, \xi'), c(\xi', \xi)\} \gamma_m \geq |\hat w_{\gamma_m, s}(\xi) - \hat w_{\gamma_m, s}(\xi')|$.
Then, the convergence $z_m \to z^\star$ follows when $W_1^{\bar c}(\bar{\mathbb P}_m, \mathbb P) \to 0$ and $\limsup_{m \to \infty} \gamma_m < \infty$.
\pfend

\discssa*
\pfstart
Let $z^\star$ and $z_m$ denote the optimal objective values of 
\eqref{eq:dualgen} and \eqref{eq:pf01} (equivalently, \eqref{eq:distdual_sp}), respectively. 
Consider any feasible pair $(\alpha_0, w)$ of \eqref{eq:dualgen} under the 
discrete nominal distribution $\bar{\mathbb P}_m$. 
Let $z$ denote its objective value, and define
\[
\gamma := \|w\|_{\mathcal V^\circ}, \quad
\alpha := \alpha_0
\quad 
s_i := w(\xi_i), \quad i \in [m].
\]
Then, the induced envelope function 
$\hat w_{\gamma,s}(\xi) := \min_{i \in [m]} \theta_{\gamma,s_i,\xi_i}(\xi)$
satisfies $\hat w_{\gamma,s} \ge f_x - \alpha_0$, and hence $(\gamma, s)$ is feasible for 
\eqref{eq:distdual_sp}. 
The corresponding objective value is
\[
  \alpha + g\left(\frac{1}{m}\sum_{i \in [m]} h(s_i)\right) + \epsilon\gamma
  = \alpha_0 + g\left(\frac{1}{m}\sum_{i\in [m]} h\circ w(\xi_i)\right) + \epsilon\|w\|_{\mathcal V^\circ}
= z,
\]
which coincides with the value of $(\alpha,w)$ in \eqref{eq:dualgen} under 
$\bar{\mathbb P}_m$. 
Since this holds for every feasible $(\alpha,w)$, we obtain $z_m \le z^\star$.
Conversely, take any optimal solution $(\gamma, \alpha, s)$ of \eqref{eq:pf01} that is non-redundant as in Lemma~\ref{lem:saacomp}, 
and let $\hat w_{\gamma,s}$ be the corresponding finite envelope function. 
By construction, $\hat w_{\gamma,s} \geq f_x - \alpha$, so $(\alpha, \hat w_{\gamma,s})$ is 
feasible for \eqref{eq:dualgen}. 
The difference between the two objective values becomes
$$
\begin{aligned}
  z^\star - z_m &= \left(\alpha + g\left(\mathbb E_{\bar{\mathbb P}_m}[h\circ \hat w_{\gamma,s}]\right)
  + \epsilon \|\hat w_{\gamma,s}\|_{\mathcal V^\circ}\right) - \left(\alpha  + g\left(\frac{1}{m}\sum_{i \in [m]} h(s_i)\right) + \epsilon\gamma\right)\\
                         & = g\left(\frac{1}{m}\sum_{i \in [m]} h\circ \hat w_{\gamma, s}(\xi_i)\right) - g\left(\frac{1}{m}\sum_{i\in [m]} h(s_i)\right) \\
                         & = g\left(\frac{1}{m}\sum_{i \in [m]} h(s_i)\right) - g\left(\frac{1}{m}\sum_{i\in [m]} h(s_i)\right) = 0.
\end{aligned}
$$
The second equality is due to gauge compatibility $\|\hat w_{\gamma, s}\|_{\mathcal V^\circ} = \gamma$, and the thrid equality is by the SAA compatibility shown in Lemma~\ref{lem:saacomp}.
\pfend

\end{document}